\documentclass[letterpaper,11pt,
]{article}
\usepackage{amsmath,amsthm,amsfonts,epsfig,setspace}
\usepackage{grffile}

\newcommand{\R}{\mathbb{R}}

\newcommand{\be}[1]{\begin{equation}\label{#1}}
\newcommand{\ee}{\end{equation}}
\renewcommand{\d}{\mathrm{d}}
\DeclareMathOperator{\supp}{supp}
\DeclareMathOperator{\dist}{dist}
\newcommand{\bo}{\partial \Omega}
\newcommand{\Id}{\mbox{I}}
\renewcommand{\r}[1]{(\ref{#1})}

\DeclareMathOperator{\WF}{WF}
\newtheorem{theorem}{Theorem}

\newtheorem{corollary}{Corollary}

\title{A New Numerical Algorithm for Thermoacoustic and Photoacoustic Tomography with Variable Sound Speed}

\author{Jianliang Qian\thanks{Department of Mathematics, Michigan State University, East Lansing, MI 48824. 
Email: {\bf qian@math.msu.edu}} 
\and Plamen Stefanov\thanks{Mathematics Department, Purdue University, West Lafayette, IN 47907. 
Email: {\bf stefanov@math.purdue.edu}}  
\and Gunther Uhlmann\thanks{Department of Mathematics, University of Washington, Seattle, WA 98195 and University of California, Irvine, CA92697-3875. 
Email: {\bf gunther@math.washington.edu}} 
\and Hongkai Zhao\thanks{Department of Mathematics, University of California, Irvine, CA92697-3875. 
Email: {\bf zhao@math.uci.edu}}
}

\usepackage{graphicx}
\graphicspath{%
    {converted_graphics/}
    {/}
    {Fig/}
}
\begin{document}
\thispagestyle{plain} \maketitle \vspace{0.5cm}

\begin{abstract}
We present a new algorithm for reconstructing an unknown source in Thermoacoustic and Photoacoustic Tomography based on the recent advances in understanding the theoretical nature of the problem. We work with variable sound speeds that might be also discontinuous across some surface. The latter problem arises in brain imaging. The new algorithm is based on an explicit formula in the form of a Neumann series. 
We present numerical examples with non-trapping, trapping and piecewise smooth speeds, as well as examples with data on a part of the boundary. These numerical examples demonstrate the robust performance of the new algorithm. 
\end{abstract}

\section{Introduction}  
Thermoacoustic (TAT) and Photoacoustic (PAT) Tomography are emerging medical imaging modalities \cite{XuWang06,wanwu07}. 
These are hybrid medical imaging methods that combine the high resolution of acoustic waves with the large 
contrast of optical waves. TAT and PAT have been developed to overcome the limitations of both conventional 
ultrasound and microwave imaging. The physical principle underlying TAT and PAT is the photoacoustic effect 
which can be roughly described as follows. A short impulse of electromagnetic microwaves or light is sent through a patient's body. The tissue heats up slightly and the heat expansion generates weak acoustic waves. These waves are measured away from the patient's body and one tries to recover the acoustic source that gives us information about the rate of absorption at each point in the body, thus creating an image. One of the potential applications for TAT 
is early breast cancer detection. The American Cancer Society reports that breast cancer is 
the second overall leading cause of death among women in the United States. The mortality rate 
from breast cancer has declined in recent years due to progress in both early detection and more 
effective treatment. Better detection techniques, however, are still needed. The significance of TAT and PAT 
is that they yield images of high electromagnetic contrast at high ultrasonic resolution 
in relatively large volumes of biological tissues \cite{XuWang06,wanwu07}.  


A first step in TAT and PAT is to reconstruct the amount of deposited energy from time-dependent 
boundary measurement of acoustic signals. To start with, we describe here the widely accepted mathematical model of TAT and PAT \cite{XuWang06,wanwu07,finchPR}. Let $\Omega\subset\R^n$ be an open set with a smooth strictly convex boundary. 
Let $c(x)>0$ be the sound speed that is either  smooth or piecewise smooth. Assume that $c(x)=1$ outside $\Omega$. 
The acoustic pressure $u(t,x)$ then solves the wave equation 
\begin{equation}   \label{1}
\left\{
\begin{array}{rcll}
(\partial_t^2 -c^2\Delta)u &=&0 &  \mbox{in $(0,T)\times \R^n$},\\
u|_{t=0} &=& f,\\ \quad \partial_t u|_{t=0}& =&0, 
\end{array}
\right.               
\end{equation}
where $T>0$ is fixed, and $f(x)$ is a source that we want to recover supported in $\bar\Omega$. The measurements are modeled by the operator
\begin{equation}\label{1b}
\Lambda f : = u|_{[0,T]\times\partial\Omega}.
\end{equation}
Given $\Lambda f$, the problem is to reconstruct the unknown $f$ that is related to the absorbing properties of the body at any point. 


There have been significant progresses in both mathematical theories and medical applications; see \cite{AgrKuchKun2008,finchPR,FinchRakesh08,Haltmeier04,Haltmeier05,Hristova08,HristovaKu08, Kruger03, Kruger99, KuchmentKun08,Patch04,SU-thermo, SU-thermo_brain, XuWang06, XuWang2006,YangWang2008,jinwan06} and references therein.
 Theoretically, one is interested in uniqueness and stability of the solution for 
the inverse problem; numerically, one is interested in designing efficient numerical algorithms to recover the solution of the inverse problem. Naturally, the above two aspects have been well studied in the case of the sound speed being  constant. In fact, if the sound speed is constant and the observation surface $\partial\Omega$ is of some special  geometry, such as planar, spherical or cylindrical surface, there are explicit closed-form inversion formulas;  see \cite{finchPR,xuwan05,Haltmeier04,Haltmeier05,finhalrak07} and references therein. 
In practice there are many cases when the constant sound speed model is inaccurate \cite{XuWang06,jinwan06,xuwan03,kuforjinxuhunwan05}. For instance in breast imaging, the different components of the breast, such as the glandular tissues, stromal tissues, cancerous tissues and other fatty issues, have different acoustic properties. The variations between their acoustic speeds can be as great as 10\% \cite{jinwan06}.

To tackle variable sound speeds in TAT and PAT, the time-reversal method has been suggested in \cite{finchPR} and used in \cite{xuwan04,gruetal07,Hristova08,HristovaKu08}. However, only when the dimension is odd and the sound speed is constant, the time-reversal method gives an exact reconstruction for $T$ large enough by the Huygens principle. 
When the dimension is even or the sound speed is not constant, the time-reversal method only yields an approximate 
reconstruction as $T\gg1$.  A natural question arises immediately: is there any exact reconstruction formula which is able 
to handle both a variable sound speed and an irregular observation geometry, for a fixed $T$? It turns out that the work  \cite{SU-thermo} summarized below does provide such a formula and this formula is the foundation for our algorithmic 
development. Some results about variable sound speeds are contained in the references above but a more complete analysis in this case has been given in \cite{SU-thermo}. This analysis includes if-and-only-if conditions for uniqueness and stability, including the case with observations on a part of the boundary.  An explicit recovery formula of a type of convergent Neumann series is derived in \cite{SU-thermo} when $\Lambda f$ is known on the whole $\bo$ and $T$ is greater than the stability threshold. This formula is the base for our numerical reconstruction, and we also apply it if the speed is trapping.

There is another situation related to a variable discontinuous sound speed which arises in brain imaging \cite{XuWang2006,wanetal08}. The skull has a discontinuous sound speed which is piecewise smooth with jump-discontinuities across the boundary of the skull.  A typical situation is that the speed in the skull is about twice as large as that in the soft tissue; see e.g., \cite{XuWang2006, YangWang2008}. Such speeds change drastically the way that singularities propagate; see Figure~\ref{fig:skull2}. In \cite{SU-thermo_brain}, the second and the third author studied this model and proved that the Neumann series expansion works as well, provided that all singularities issued from $\supp f$ have a path of non-diffractive segments that reaches $\bo$ until time $T$.  Based on this work, we develop an efficient numerical algorithm as well 
to handle the situation of variable discontinuous sound speeds in PAT and TAT. 

The rest of the paper is organized as follows. In Section \ref{prep}, we provide some preliminary materials for 
understanding theories summarized in the following sections. In Section \ref{smoothwhole} we summarize the 
theoretical results in \cite{SU-thermo} and explain consequences of the theory.  In Section \ref{smoothpart} 
we summarize the theoretical results in \cite{SU-thermo} for the case of partial data. In Section \ref{discont} 
we summarize the theoretical results in \cite{SU-thermo_brain} when the sound speed is discontinuous. 
In Section \ref{smoothvsnonsmooth} we make comparison between smooth and nonsmooth sound speeds in terms 
of uniqueness and stability. In Section \ref{algorithm} we summarize the iterative algorithm for constructing 
the Neumann series so that the inversion formula in \cite{SU-thermo} can be implemented.  In Section \ref{numerics} 
we give an extensive examples to demonstrate the robustness of the new algorithm for PAT and TAT.  

\section{Preliminaries}
\label{prep}
Assume for now that $c>0$ is smooth. 
The speed $c$ defines a Riemannian metric $c^{-2}\d x^2$. For any piecewise smooth curve $[a,b]\mapsto \gamma\in \R^n$, the length of $c$ in that metric is given by
\[
\textrm{length}(c) = \int_a^b \frac{|\dot\gamma(t)|}{c(\gamma(t))}\d t.
\]
The so defined length is independent of the parameterization of $\gamma$. 
The distance function $\dist(x,y)$ is then defined as the infimum of the lengths of all such curves connecting $x$ and $y$. 

For any $(x,\theta)\in \R^n \times S^{n-1}$ 
we denote by $\gamma_{x,\theta}(t)$ the unit speed (i.e., $|\dot\gamma|=c(\gamma)$) geodesics issued at $x$ in the  direction $\theta$.


Recall that the energy of $u(t,x)$ in a domain $U$ is given by
\[
E(u(t)) = \int_U \left(  |\nabla_x u|^2+c^{-2}|u_t|^2      \right) \d x,
\]
where $u(t)=u(t,\cdot)$. The energy of any Cauchy data $(f,g)$ for equation \r{1} is given by the same integral with $\nabla_x u=\nabla_x f$ and $u_t=g$. In particular, the energy of $(f,0)$ in $U$ is given by the square of the Dirichlet norm
\[
\|f\|_{H_D(U)}^2 := \int_U |\nabla_x f|^2\,\d x.
\]
We always assume below that $f\in H_D(\Omega)$, where the latter is the Hilbert space defined by the norm above. We will denote by $\|\cdot\|$ the norm in $H_D(\Omega)$, and in the same way we denote the operator norm in that space. 

There are two main geometric quantities that  are crucial for the results below. First we set
\be{T0}
 T_0 := \max\{\dist(x,\bo);\;x\in\bar\Omega\}.
\ee
Let $T_1\le\infty$ be the supremum of the lengths of all maximal geodesics lying in $\bar\Omega$. Clearly, $T_0< T_1$ but while the first number is always finite, the second one can be infinite. It can be shown actually that
\be{T1}
T_0\le T_1/2.
\ee

One of the main ingredients of our approach in \cite{SU-thermo, SU-thermo_brain} is understanding the microlocal nature of the problem. We recall the definition of a wave front set of a function, or more generally, a distribution; see \cite{hor71a}. The definition is based on the known property of the Fourier transform: one can tell whether a compactly supported function $f$ is smooth by looking at the decay of the Fourier transform $\hat f(\xi)$ as $|\xi|\to\infty$:   $f\in C_0^\infty$ if and only if  $|\hat f(\xi)|\le C_N(1+|\xi|)^{-N}$ for any $N$. The idea behind the wave front set is to localize this near a fixed $x_0$ and in a conic neighborhood of a fixed $\xi_0\not=0$. Conic neighborhoods are defined as open conic sets, i.e., sets of the type $\Gamma = \{r\theta; \; r>0,\, \theta\in V\}$, where $V$ is an open subset of $ S^{n-1}$. 
We say that $(x_0,\xi_0)\not\in \WF(f)$, $\xi_0\not=0$, if there exist $\phi\in C_0^\infty$ with $\phi(x_0)\not=0$ and a conical neighborhood $\Gamma$ of $\xi_0$, 
so that
\[
|\widehat{\phi f}(\xi)| \le C_N(1+|\xi|)^{-N}, \quad \forall \xi\in \Gamma, \forall N. 
\]
If $(x,\xi)\in \WF(f)$, we say that $(x,\xi)$ is a singularity of $f$, or that $f$ is singular at $(x,\xi)$. Since singularities are defined by conic sets, we can restrict $\xi$ to unit vectors. For example, the Dirac Delta function $\delta(x)$ has wave front at $x=0$ and all directions, i.e., $\WF(\delta) = \{(0,\xi);\; \xi\not=0\}$, and a piecewise smooth function $f$ that has a jump across some smooth surface $S$ (and nowhere else) is singular at all points of $S$ in (co)normal directions, i.e., $\WF(f) = \{(x,\xi);\; x\in S, \; 0\not =\xi\perp S\ \text{at $x$}\}$. 

The propagation of singularities for the wave equation \r{1} can be described as follows. If $(x,\theta)\in \WF(f)$, then at time $t$, both $(\gamma_{x,\theta}(t), \dot \gamma_{x,\theta}(t))$ and $(\gamma_{x,\theta}(-t), \dot \gamma_{x,\theta}(-t))$ are in $\WF(u(t,\cdot))$, where $u$ is the solution of \r{1}. This is due to the fact that the symbol of the  wave operator has two sound speeds, $\pm c(x)|\xi|$, and that the initial velocity on \r{1} is zero, therefore each singularity splits into two equal parts starting to propagate in opposite directions. While this is a classical result in the linear PDE theory, we refer to \cite{SU-thermo} for more details in this specific case.

\section{Smooth speed and data on the whole $\bo$} 
\label{smoothwhole}
We will describe below the theoretical results in \cite{SU-thermo}.

\subsection{Uniqueness} If $T\gg1$, $\Lambda f$ recovers $f$ uniquely. We  have the following sharp result based on the unique continuation theorem by Tataru \cite{Tataru99}. 

\begin{theorem}\label{thm_uniq1}
Let $\Lambda f=0$. Then $f(x)=0$ for $\dist(x,\bo)\le T$.  Moreover, $f$ can be arbitrary in the set $\dist(x,\bo)> T$, if the latter set is non-empty. 
\end{theorem}

\begin{corollary}\label{cor_uniq1}

$\Lambda$ is injective on $H_D(\Omega)$ if and only if  $T\ge T_0$.
\end{corollary}

We refer to \cite{SU-thermo} for proofs.

\subsection{Stability} We showed in \cite{SU-thermo} that  $\Lambda f$ recovers $f$ in a stable way if each singularity $(x,\xi)$, i.e., each element of the wave front set $\WF(f)$, reaches $\bo$ for time $t$ (positive or negative) such that $|t|<T$. 
In other words, if functions $f$ are a priori supported in a fixed compact $\mathcal{K}\subset\bar\Omega$, then we have 
the following equivalent statement: 
\be{cond_stab}
\begin{split}
& \text{For any $(x,\theta)\in \mathcal{K} \times S^{n-1}$, the unit speed geodesic through $(x,\theta)$ at $t=0$}\\& \text{reaches $\bo$ at time $|t|<T$.}
\end{split}
\ee
Moreover, the following condition is sufficient regardless of the choice of $\mathcal{K}$:
\be{cond_stab2}
T_1/2<T.
\ee
This condition is equivalent to \r{cond_stab} if $\mathcal{K}=\bar\Omega$.

Furthermore, to formulate a condition equivalent to \r{cond_stab} by taking $\mathcal{K}$ into account, we define $T_1=T_1(\mathcal{K})$ as above which is explicitly related to geodesics that pass through $\mathcal{K}$. 

If $T_1=\infty$, then the sound speed $c$ is called trapping (in $\Omega$). In this case, there is no stability regardless of the choice of $T$. 

We summarize the above discussion into the following.

\begin{theorem}\label{thm_stab} Let $\mathcal{K}\subset \Omega$ be compact. 

(a) 
Let $T> T_1(\mathcal{K})/2$. Then there exists a constant $C>0$ so that
\[
\|f\|\le C\|\Lambda f\|_{H^1([0,T]\times\bo)}.
\]

(b) Let $T< T_1(\mathcal{K})/2$. Then for any $C>0$, $s_1$ and $s_2$, there is $f\in C^\infty$ supported in $\mathcal{K}$ so that
\[
\|f\|_{H^{s_1}}\ge  C\|\Lambda f\|_{H^{s_2}([0,T]\times\bo)}.
\]
\end{theorem} 

In other words,  $T> T_1(\mathcal{K})/2$  is a sufficient and necessary condition for stability for functions $f$ supported in $\mathcal{K}$ up to replacing the $<$ sign by $\le$. 

\medskip 
\textbf{Visible and Invisible Singularities.} The condition \r{cond_stab} can be explained in the following way. As a general principle, for a stable recovery in such a linear inverse problem, we need to detect all singularities. We refer to \cite{SU-JFA09} for more details. As explained above, each singularity starts to travel in a positive and a negative direction because $u_t=0$ at $t=0$, so it can leave two traces on $\bo$. It is enough to detect one of them for stability. On the other hand, if we can detect both, we can expect better numerical results.
Condition \r{cond_stab} then says that all singularities are visible at the boundary. 

We want to emphasize that $T_1$ can be much larger that $\text{diam}(\Omega) :=\max\{\dist(x,y);\; (x,y)\in \bo\times\bo\}$. If there are no conjugate points in $\bar\Omega$, then those two quantities coincide. If there are conjugate points however, this is not necessarily so. If there are closed geodesics, then $T_1=\infty$, while  $\text{diam}(\Omega)$ is finite. In this case, although we always have uniqueness for some $T\gg1$, we never have stability. Finally, we notice that while $T_0$ and $\text{diam}(\Omega)$ can always be  estimated analytically or numerically, $T_1$ is much harder to estimate or even to tell whether it is finite or not. In any case,
\be{T}
\text{diam}(\Omega)\le T_1. 
\ee
If $\gamma_{x,\theta}$ does not hit the boundary at time $|t|\le T$, we call $(x,\theta)$  an invisible (possible) singularity. It is easy to show that if one such pair exists, then there is a non-empty open set of invisible singularities.

\subsection{Reconstruction}
 The reconstruction method in \cite{SU-thermo} is based on the following ideas. If we knew the Cauchy data $(u,u_t)$ on $\{T\}\times \Omega$, we could just solve a mixed problem like the one below with that Cauchy data on $t=T$ and boundary data given by $\Lambda f$. Although we do not know $(u,u_t)$ on $\{T\}\times \Omega$, we do know the boundary values of $u$ on $\partial\Omega$ for $t=T$, i.e., $u$ on $\{T\}\times\partial\Omega$. Assuming that $[f,0]$ is in the energy space, i.e., $f\in H_D(\Omega)$, we can only say that $u_t(T,\cdot)$ is in $L^2(\Omega)$, and its boundary values might not be well defined. Now from all possible functions with prescribed boundary values on $\{T\}\times \Omega$, we choose the one that minimizes the energy norm $\|\cdot\|_{H_D(\Omega)}$. By the Dirichlet principle, it is given by the harmonic extension of $u(T,\cdot)|_{\partial\Omega}$. 

Consequently, given $h$ defined on $[0,T]\times \bo$ which eventually will be replaced by $\Lambda f$, we first solve the elliptic boundary value problem
\be{3}
\Delta\phi=0, \quad 
\phi|_{\partial\Omega} = h(T,\cdot),
\ee 
and we introduce the notation $P_\Omega$ for the Poisson operator of harmonic extension: $P_\Omega h(T,\cdot) := \phi$. 
In fact, the equation is $c^2\Delta u=0$ but $c^2$ cancels out. Then we perform the modified back projection:
\begin{equation}   \label{2}
\left\{
\begin{array}{rcll}
(\partial_t^2 -c^2\Delta )v &=&0 &  \mbox{in $(0,T)\times \Omega$},\\
v|_{[0,T]\times\partial\Omega}&= &h,\\
v|_{t=T} &=& P_\Omega h(T,\cdot),\\ \quad   v_t|_{t=T}& =&0. \\
\end{array}
\right.               
\end{equation}
 Note that the initial data at $t=T$ satisfy compatibility conditions of first order (no jump at $\{T\}\times\bo$). 
Then we define the following left pseudo-inverse of $\Lambda$
\be{4}
A h := v(0,\cdot) \quad \mbox{in $\bar\Omega$}.
\ee
The operator $A$ is not an actual inverse (unless $n$ is odd, $c$ is constant,  and $T$ is greater than the diameter) and we have
\[
A\Lambda = \Id -K,
\]
where $K$ is an ``error'' operator. We showed in \cite{SU-thermo} that 
\be{4K}
\|Kf\|_{H_D(\Omega)}\le \|f\|_{H_D(\Omega)}, \quad\forall f\in H_D(\Omega),
\ee
for any smooth speed, trapping or not, and for any time $T>0$. If $T>T_0$, the inequality is strict, see \cite{SU-thermo}. To show that $K$ is a contraction requires some assumptions, however.


\begin{theorem}  \label{thm2.1} \ 

(a) Let $c$ be non-trapping, and let $T>T_1/2$. Then 
$A\Lambda=\Id-K$, where   $\|K\|_{H_{D}(\Omega)\to {H_{D}(\Omega)}}<1$. 
In particular, $\Id-K$ is invertible on $H_{D}(\Omega)$, and the inverse thermoacoustic problem has an explicit solution of the form
\be{2.2}
f = \sum_{m=0}^\infty K^m A h, \quad h:= \Lambda f.
\ee

(b) Let $T>T_1$. Then in addition to the conclusions above, $K$ is compact in $H_{D}(\Omega)$.
\end{theorem}
For a proof, see \cite[Theorem~1]{SU-thermo}  and \cite[Remark~2.2]{SU-thermo_brain}.

In other words, under condition \r{cond_stab2}, we have not only stability but an explicit solution in a form of a convergent Neumann series as well. On the other hand, if \r{cond_stab2} fails, i.e., if $T<T_1/2$, there is no stability by Theorem~\ref{thm_stab}. This does not mean that the series \r{2.2} would not converge  in this case. 
If it does, it will converge to $f$ for $T>T_0$. Indeed, then it is easy to see that the limit $g$ solves $(\Id-K)(g-f)=0$, and since \r{4K} is strict in this case, $f=g$.

Based on this theorem and its proof, we can expect good convergence when $T>T_1$, and the first term would already be a good approximation of the high frequency part of $f$. 

If $T_1/2<T<T_1$, we can expect the first term in the series to recover only a fraction of the high frequency part of $f$, and the successful terms to improve this gradually; the series \r{2.2} would still converge but that convergence would be slower. 

If $T<T_1$, then $\|K\|=1$. Indeed, if we assume that $\|K\|<1$, we would get uniform convergence of the Neumann series, and stability; on the other hand, there is no stability in this case. 

\medskip

\subsection{Summary: Dependence on $T$}

\begin{itemize}
  \item[(i)]{$T<T_0$} \  \par 
$\Lambda f$ does not recover $f$ uniquely; see Theorem~\ref{thm_uniq1}. Then $\| K\|=1$, and for any $f$ supported in the inaccessible region, $Kf=f$. 

 \item[(ii)]{$T_0< T<T_1/2$} \  \par 
This can happen only if there is a strict inequality in \r{T1}. Then we have uniqueness but not stability. In this case, $\|K\|=1$, $\|Kf\|<\|f\|$, and we do not know if the Neumann series \r{2.2} converges. If it does, it converges to $f$. 
  \item[(iii)]{$T_1/2<T<T_1$} \  \par 
This assumes that $\Omega$ is non-trapping for $c$.  
The Neumann series \r{2.2} converges exponentially but maybe not as fast as in the next case. There is stability, and $\|K\|<1$. 
  \item[(iv)]{$T_1<T$} \  \par 
This also assumes that $\Omega$ is non-trapping for $c$. The Neumann series \r{2.2} converges exponentially. There is  stability, $\|K\|<1$, and $K$ is compact. 
\end{itemize}

\subsection{Recovery of singularities} 
In cases (iii) and (iv) above, we recover explicitly the whole $f$, including its singularities. If the goal is to recover only the singularities of $f$ (the wave front set $\WF(f)$) with less computation, then one can perform the classical back-projection (time reversal) as follows; see \cite{xuwan04,gruetal07,Hristova08,HristovaKu08,KuchmentKun08}. Assuming the non-trapping condition, $T_1<\infty$, let $\chi\in C_0^\infty(\R)$ be such that $\chi(t)=1$ for $t\in [0,T_1]$. 
Set
\be{BP}
Rf = A\chi \Lambda f.
\ee
In general, $Rf$ is not close to $f$ unless $T\to\infty$; see \cite{Hristova08}. 

In case (iv), $R$ is a parametrix of infinite order; see \cite{SU-thermo}. Therefore, it recovers correctly all singularities, including jumps across smooth surfaces --- it will recover correctly the location and the size of the jump. 

In case (iii), $R$ is elliptic but not a parametrix itself; see \cite[Theorem~3]{SU-thermo}. Then $Rf$ will have the singularities at the right places but the amplitudes will be  in general between $1/2$ and $1$.  

In cases (i) and (ii), only singularities ``close enough to the boundary'' will be recovered with amplitudes between $1/2$ and $1$. 

Finally, if the speed is trapping in $\Omega$, i.e., if $T_1=\infty$, then $Rf$ recovers the visible singularities up to time $T$ 
by choosing $\chi$ appropriately. 

These comments are just another way to formulate \cite[Theorem~3]{SU-thermo}.

\subsection{Comparison with the Time Reversal Method} 
Let $T_1<\infty$ first, i.e., assume that $c$ is non-trapping. Let the time reversal approximation be defined as in \r{BP} with $\chi=1$ in a neighborhood of  $[0,T_1]$. Then $R$ is a parametrix of infinite order, i.e., $Rf=f-Qf$, where $Qf\in C^\infty$ for any $f$. On the other hand, $\|Q\|$ is not necessarily small for any fixed $T$. Assuming that $\chi$ is properly chosen, as $T\to\infty$, that norm gets smaller at a rate dictated by the local energy decay for the wave equation: $\|Q\|=O(t^{1-n})$ for $n$ even 
and $\|Q\|=O(e^{-Ct})$ for $n$ odd; see \cite{Hristova08}. Therefore, for $T\gg1$ so that $\|Q\|<1$, one can write
\[
R\Lambda f = (\Id-Q)f,
\]
and one solve this equation by Neumann series as above. However, it is not straightforward to impose sharp conditions on $\chi$ and $T$ to guarantee $\|Q\|<1$. On the other hand, for the method that we propose, that condition is $T>T_1$ and it is sharp for a stable inversion. Also, the proposed method minimizes the norm of the ``error'' operator, and when both Neumann expansions converge uniformly, the one in Theorem~\ref{thm2.1} will converge faster in the uniform topology. Numerical experiments not shown here confirm that. 

If $T_1=\infty$ ($c$ is trapping), then the error in the time reversal method decays like $O(1/\log T)$ if $f\in H_0^2(\Omega)$, and the error decays even slower if $f\in H_D(\Omega)$ only. The first term of the Neumann series inversion has an error no less than that, and numerically the error improves with a few more terms. We do not know whether it converges or not, however. 

\section{Smooth speed and data on a part of $\bo$} 
\label{smoothpart}
Let $\Gamma\subset\bo$ be a relatively open set of $\bo$, and assume that we only have data available on $[0,T]\times\Gamma$. We will suppose that $f$ is supported in some compact $\mathcal{K}\subset\Omega$. 
\subsection{Uniqueness} 

As in \r{T0}, set
\be{T0p}
 T_0 := T_0(\mathcal{K},\Gamma)=  \max\{\dist(x,\Gamma);\;x\in\mathcal{K}\}.
\ee
Theorem~\ref{thm_uniq1} has the following analog in this case. 

\begin{theorem}\label{thm_pd_uniq}
Let $\Lambda f=0$ on $[0,T]\times\Gamma$. If $T\ge T_0$, then $f=0$. If $T<T_0$, then $f=0$ on $\mathcal{K}\cap \{x;\; \dist(x,\Gamma)<T\}$ and can be arbitrary in the complement of this set in $\mathcal{K}$. 
\end{theorem}

The proof of this theorem is not easy. It combines Tataru's uniqueness theorem with arguments that first appeared in \cite{finchPR} in the case of constant speed and were extended later in \cite{SU-thermo} to variable speeds. 

As a corollary,  $\Lambda f|_{[0,T]\times\Gamma}$ determines uniquely $f$ if and only if $T\ge T_0$. 

\subsection{Stability} 
The following condition guarantees that we can detect all singularities originating from $\mathcal{K}$ as singularities of our data: 
\be{cond_stab_pd}
\begin{split}
& \text{For any $(x,\theta)\in \mathcal{K} \times S^{n-1}$, the unit speed geodesic through $(x,\theta)$ at $t=0$}\\& \text{reaches $\Gamma$ at time $|t|<T$.}
\end{split}
\ee
Let $T_1/2  \le\infty$ be the maximum    of all such times. 
 The following condition then is equivalent to \r{cond_stab_pd} 
\be{cond_stab2a}
T_1/2<T.
\ee

\begin{theorem}\label{thm_stab2} Let $\mathcal{K}\subset \Omega$ be compact. 

(a) Let $T> T_1/2$. Then there exists a constant $C>0$ so that
\[
\|f\|\le C\|\Lambda f\|_{H^1([0,T]\times\Gamma)}.
\]

(b) Let $T< T_1/2$. Then for any $C>0$, $s_1$ and $s_2$, there is $f\in C^\infty$ supported in $\mathcal{K}$ so that
\[
\|f\|_{H^{s_1}}\ge  C\|\Lambda f\|_{H^{s_2}([0,T]\times\Gamma)}.
\]
\end{theorem} 

In other words,  $T> T_1/2$  is a sufficient and necessary condition for stability for functions $f$ supported in $\mathcal{K}$ up to replacing the $<$ sign by $\le$. 

\subsection{Reconstruction} 
An explicit formula of the type \r{2.2} is not available in this case but one can show that the recovery is reduced to a Fredholm equation if \r{cond_stab2a} holds. 

Let $\chi\in C^\infty([0,T]\times\Gamma)$ be a cutoff function supported in $[0,T]\times\Gamma$ so that $\chi=1$ on a slightly smaller set $[0,T']\times\Gamma'$ that still satisfies the condition. Assume also that $0\le\chi\le1$. 
Then we know $\chi\Lambda f$.  Apply the time reversal operator $A$ to that to get $A\chi\Lambda f$. Since $\chi(\cdot,T)=0$, the harmonic extension is zero, so this is the classical time reversal. Let $K$ be the ``error operator''
\[
K = \Id-A\chi\Lambda.
\]
To find $f$, we need to solve
\be{pd1}
(\Id-K)f=h, \quad h := A\chi\Lambda h
\ee
with $h$ determined by the data.  By \cite[Theorem~3]{SU-thermo}, $A\chi\Lambda$ is a pseudo-differential operator that is elliptic under the stability condition \r{cond_stab_pd}. Its principal symbol is given by 
\[
\frac12 \chi(\gamma_{x,\xi}(\tau_+(x,\xi))) + \frac12 \chi(\gamma_{x,\xi}(\tau_-(x,\xi))) ,
\]
where $\tau_\pm(x,\xi) $ is the positive/negative time of the (unit speed) geodesic $\gamma_{x,\xi}$ through $(x,\xi)$ to reach $\bo$. Then $K$ is also a pseudo-differential operator with principal symbol
\be{sk}
\sigma_p(K) = 1-\frac12 \chi(\gamma_{x,\xi}(\tau_+(x,\xi))) - \frac12 \chi(\gamma_{x,\xi}(\tau_-(x,\xi))).
\ee
Since $0\le\chi\le1$, and by the stability condition, at least one of the $\chi$ terms above is equal to $1/2$. Therefore,
\[
0\le \sigma_p(K)\le \frac12,
\]
and  in the set where $\chi$ is not  $0$ nor $1$ (that is relatively small in applications, but not too small to keep $|\partial_t\chi|+|\partial_x \chi|$ under control),   $\sigma_p(K)$ is either  $0$ or $1/2$. By the G\.{a}rding inequality, one can express $K$ in the form
\[
K = K_1+K_2,\quad \|K_1\|\le\frac12, \quad \text{$K_2$ compact},
\]
and if $\Gamma\not=\emptyset$, one can show that actually $\|K_1||=1/2$. This does not allow us to claim that the Neumann series \r{2.2} converges. However, by \cite{Medkova}, since the distance from $K$ to the compact operators is less than one, for a fixed $f$,
\be{NS2}
\sum_{m=0}^\infty K^mA\chi \Lambda f
\ee
converges if and only if $K^mA\chi \Lambda f\to 0$ as $m\to\infty$. It is trivial to see that the limit $g$ solves
\be{pd_A}
A\chi \Lambda g = A\chi \Lambda f.
\ee
While we know that $\chi\Lambda$ is injective for $T>T_0$, this is not true for $A$, so we cannot draw the conclusion that $g=f$. The original equation $\chi\Lambda f=h$ with $h$ in the range of $\Lambda$ is Fredholm with a trivial kernel, but the equation \r{pd_A} for $g$ is Fredholm which may have a non-trivial kernel. 

On the other hand, the partial sums
\be{pd_P}
g_{N} := \sum_{m=0}^N K^mA\chi \Lambda f
\ee
used in the reconstruction below, recover singularities of $f$ asymptotically as $N\to\infty$. This follows from the standard pseudo-differential calculus.

For the purpose of the numerical reconstruction, we take $\chi$ to be a function of $x\in\bo$ only, independent of $t$. In particular, we no longer have $\chi=0$ near $t=T$. Then we define $A$ as before but this time the harmonic extension of $\chi\Lambda f(T,\cdot)$ is not trivial in general. Then we use partial sums as in \r{pd_P}. 

If the stability condition \r{cond_stab2a} is not satisfied, then inverting $\chi\Lambda$ is not equivalent to solving an Fredholm equation anymore and it is unstable. One can show that $K$ is the sum of an operator with norm not exceeding $1$ and a compact one. Of course, this does not guarantee convergence of the Neumann series and the partial sums recover asymptotically the visible singularities only. In the numerical example below, we check the norm of each successive term in the partial sum \r{pd_P} and stop when that norm starts to increase.

\section{Discontinuous sound speed; modeling brain imaging} 
\label{discont}
We will briefly review the results in \cite{SU-thermo_brain}. Let $c(x)$ be piece-wise smooth with a non-zero jump across one or several smooth closed non-intersecting surfaces that we call $S$ in $\Omega$. We still assume that $c=1$ outside $\Omega$. The formulation of the problem is still the same and classical solutions $u(t,x)$ are assumed to be in $C^1$, and this implies the following transmission conditions across $S$ where $c(x)$ jumps: the limits of $u(t,x)$ and its normal derivative match as $x$ approaches $S$ from either side. 

\subsection{Propagation of singularities} 
It is well known that this case is quite different from the previous one from the viewpoint of propagation of singularities. Let us see what happens when a ``ray'' (a geodesic) approaches $S$ from one of its sides, that we call an interior one. Let $c_-$ be the limit 
of $c$ on $S$ from inside, and let $c_+$ be that from outside. 

If $c_->c_+$, this ray splits into two parts when hitting $S$. One of them reflects according to the usual laws of reflection and goes back to the interior. Another one goes into the exterior and refracts by changing its angle with respect to the boundary. The incoming and the outgoing angles $\alpha_\pm$ with the normal satisfy Snell's Law: 
\be{snell}
\frac{\sin\alpha_-}{\sin\alpha_+} = \frac{c_-}{c_+}.
\ee
Assume now that $c_-<c_+$. Then there is a critical angle $0<\alpha_0<\pi$ with the normal at any point so that if $\alpha_-<\alpha_0$, there are still a reflected and a transmitted (refracted) ray as above satisfying Snell's law. If $\alpha_->\alpha_0$, then there is no refracted ray, while the reflected one still exists. This is known as a full internal reflection. This critical angle can be read off from the Snell Law: it is the value of $\alpha_-$ that forces $\sin\alpha_+$ to be equal to $1$, and therefore to be greater than $1$ when $\alpha_- >\alpha_0$. Therefore, $\alpha_0=\sin^{-1}(c_-/c_+)$. Propagation of singularities when $\alpha_-=\alpha_0$ is more delicate and will not be analyzed here. 

After splitting or not into two rays, each segment may split into two, etc. Assuming that we never get rays tangent to $S$, reflection and refraction occur in the same way as above. 
Figure~\ref{fig:skull2} below shows a possible time evolution of a single singularity, both for positive and  negative time. The sound speed in the ``skull'' is higher than that on either side.  

\begin{figure}  
  \centering
  \includegraphics[bb=112 37 797 514,width=2.98in,height=2.08in,keepaspectratio]{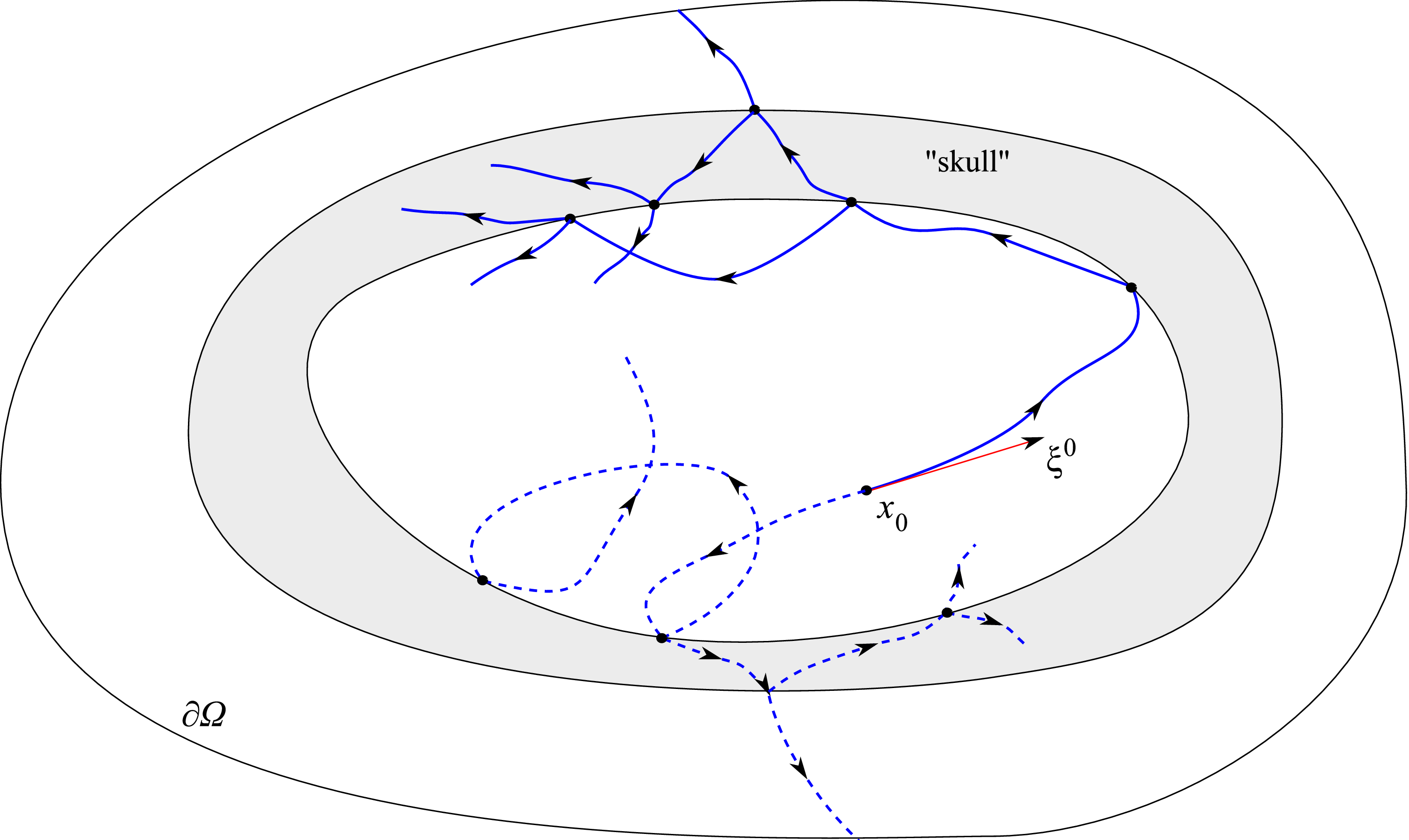}
  \caption{Rays in the case of a discontinuous sound speed. Dotted lines represent negative times.}
  \label{fig:skull2}
\end{figure}

We define the distance function $\dist(x,y)$ as the infimum of the length (in the metric $c^{-2}\d x$) of all piecewise smooth curves connecting $x$ and $y$ that intersect $S$ transversely. 

\subsection{Uniqueness} 
Theorem~\ref{thm_uniq1} and Corollary~\ref{cor_uniq1} still hold in this case; see \cite[Proposition~5.1]{SU-thermo_brain} and 
its proof. The definition of $T_0$ in this case is the same.  

\subsection{Stability. Visible and invisible singularities} 
The general principle that we should be able to detect each singularity for stable recovery still applies. However, in the case of full internal reflection, there might be singularities that never reach $\bo$. For example, let  $\Omega$ be the ball $B(0,R)$ (centered at $0$ with radius $R$), and let $S$ be the sphere $|x|=R_1$, $R_1<R$. Assume also that inside the sphere $S$, the speed is a constant less than outside $S$, where it is also a constant. Then all rays coming from the interior of $S$, hitting $S$ at an angle smaller than the critical one, will reflect and give rise to no transmitted rays. By the rotational symmetry, those reflected rays will hit $S$ again at the same angle, fully reflect, etc. Such  singularities will be invisible, no matter what $T$ we choose. This is similar to the trapping situation in the case of  a smooth sound speed. Of course, we may also have rays that are trapped but never reach the  interface $S$. 
Therefore the singularities that are certain to be visible up to time $T$ consist of the following set: 
\be{2.1a'}
\begin{split}
\mathcal{U} = &\big\{(x,\theta)\in (\Omega\setminus S) \times S^{n-1} ;\;  \text{there is a path  of the  ``geodesic'' issued from either}\\
& \text{ $(x,\theta)$ or $(x,-\theta)$  at $t=0$ never tangent to $S$, that is outside $\bar\Omega$ at time $t=T$} \big\}.
\end{split}
\ee
If we restrict our attention to $f$ a priori supported in some compact $\mathcal{K}\subset \Omega\setminus S$, then the stability condition can be formulated as follows:
\be{2.1ab}
\begin{split}
&\forall (x,\theta)\in \mathcal{K} \times S^{n-1} ;\;  \text{there is a path  of the  ``geodesic'' issued from either}\\
& \text{ $(x,\theta)$ or $(x,-\theta)$  at $t=0$ never tangent to $S$, that is outside $\bar\Omega$ at time $t=T$} .
\end{split}
\ee
The reconstruction formula below shows that this condition is sufficient not only for stability but also for an explicit convergent formula of the type \r{2.2}.

Similar to the smooth case, we can define $T_1/2$ (instead of defining $T_1$ directly) to be infimum of all $T$ satisfying \r{2.1ab} when $\mathcal{K}=\bar\Omega$. Then \r{cond_stab2} remains a necessary and sufficient condition (up to replacing $<$ by $\le$) for stability in this case as well.

\subsection{Reconstruction} 
In \cite{SU-thermo_brain}, we analyzed also the energy at high frequencies that is carried by the reflected and the transmitted rays. We showed that if none of the incident, reflected and transmitted rays are tangent to $S$, then a positive fraction of the energy reflects and transmits at high frequencies. Therefore, even though under the condition \r{2.1ab} the corresponding singularity is visible at $\bo$, only a fraction of its would be measured. This fraction could be quite small, depending on the speed,  the number of reflections and refractions before the ray hits $\bo$, and the angles there. So we should not expect the first term in the series \r{2.2}, even if it converges, to be a good approximation even at high frequencies. The theorem below, proved in \cite{SU-thermo_brain}, shows that an analog of \r{2.2} still converges but it needs to be modified first.

The needed modification is connected to the following problem. When applying $A$ to $\Lambda f$ with $\supp f \subset\mathcal{K}$, we get $A\Lambda f$ that may be supported everywhere in $\Omega$. The successful terms in \r{2.2} would then force $A\Lambda$ to be applied to the result. To ``restrict'' $A\Lambda f$ to $\mathcal{K}$, we cannot just restrict in the usual sense, because this may take the results out of the energy space (no zero trace on $\partial\mathcal{K}$). For this reason, we project orthogonally on $H_D(\mathcal{K})$. It turns out that the orthogonal projection is given by $\Pi_\mathcal{K}f := f -P_\mathcal{K} (f|{\partial\mathcal{K}})$, where $P_\mathcal{K}$ is the Poisson operator of harmonic extension from $\partial\mathcal{K}$ to $\mathcal{K}$ defined in \r{3}. 

\begin{theorem}  \label{thm2.1a} 
Let $\mathcal{K}$ satisfy \r{2.1ab}. Then  ${\Pi}_{\mathcal{K}}{A}\Lambda_1=\Id-{K}$ in $H_D(\mathcal{K)}$, with  
 $\|{K}\|_{H_D(\mathcal{K})}<1$. In particular, $\Id-{K}$ is invertible on $H_D(\mathcal{K})$, and $\Lambda$ restricted to $H_D(\mathcal{K})$  has an explicit left inverse of the form
\be{2.2'}
{f} = \sum_{m=0}^\infty  {K}^m  {\Pi}_{\mathcal{K}} A h, \quad h= \Lambda f.
\ee
\end{theorem}

Note that the theorem does not say how to reconstruct $f$ when $\supp f$ is not in a set $\mathcal{K}$ satisfying \r{2.1ab}. Estimate \r{4K} still holds but we do not know if $\Id-K$ can be inverted by a Neumann series. 
In the numerical examples, however, we work with $f$ supported everywhere in $\Omega$, and then $\Pi_\mathcal{K}=\Pi_{\bar\Omega}=\Id$. 

\section{Comparison between the smooth and the non-smooth case} 
\label{smoothvsnonsmooth}
If the sound speed $c$ is smooth, then there is always uniqueness for large enough $T$ because $T_0<\infty$; however, there is stability and the Neumann series converges to a solution only when $c$ is non-trapping and $T>T_1$. In the trapping case, one still has stability and an explicit solution of the type in Theorem~\ref{thm2.1} under the a priori assumptions that $\supp f\subset\mathcal{K}$ and all singularities with a base point $x$ in $\mathcal{K}$ are visible. The assumption $\supp f\subset\mathcal{K}$ means that we need to know $f$ outside $\mathcal{K}$; then we can subtract that part from $f$ and apply the reconstruction procedure. 

When $c$ is of jump type, we may still have rays that are trapped even if they never reach $S$. On the other hand, $S$ may split some rays. Under the stability condition \r{2.1ab} which is equivalent to \r{cond_stab2} with the modified definition of $T_1$, we still have stability and an explicit solution for $\supp f\subset\mathcal{K}$. We still need to know  $f$ outside $\mathcal{K}$ for a full reconstruction. The essential difference is that even though all singularities except those (of measure zero) contributing to tangent rays on $S$ will be detected in the best possible case, $T_1<\infty$, they might be detected with a loss of energy. In the example presented in Figure~\ref{fig:skull2}, the singularity that exits at the top has been split twice before that, and each such event takes away a fraction of the energy. In contrast, in the case of a smooth sound speed, such a singularity will exit for $t>0$  with half of the energy of $f$  while the other half will exist for negative time; as such, for example, we will detect both in case (iv). Going back to the non-smooth case, the Neumann series \r{2.2'} will gradually restore the right strength of the singularity, and the whole $f$, actually, but we can expect this to be more sensitive to noise and computational errors. 

\section{Algorithmic formulations}
\label{algorithm}
The main issue in implementing the above formulation is how to compute the operator $K$. In the case of data available on the whole $\bo$, by definition, $K=Id-A\Lambda$. Thus, given a function $\psi$, 
\[K\psi = \psi-A\Lambda\psi.\] 
Thus we need an efficient computational algorithm to carry out the actions of time-reversal 
operator $A$ and the measurement operator $\Lambda$. Since the measurement 
operator $\Lambda$ can be simulated by forward modeling, we will detail the implementation 
of forward modeling first. 

\subsection{Complex Scaling/Perfectly matched layers (PML) for the acoustic wave equation}
Since the forward model equation \eqref{1} is formulated as a pure initial value problem, we have to truncate the computational domain to be finite. 
The truncated domain has to be large enough to enclose the domain where the measurements are taken along its boundary. On the other hand, we have to impose some artificial boundary conditions on the boundary of the truncated domain. To perform accurate 
long-time simulation, we will adopt the complex scaling method (see \cite{Sjostrand_Zw_91} and the references therein), 
also known as perfectly matched layer (PML) as an absorbing boundary condition for acoustic waves \cite{ber94,chewee94,liutao97}. 
This is done so that waves will not be reflected into the computational domain.  

To derive the PML for the equation, we follow \cite{liutao97}. We rewrite the 2nd-order 
equation as a first-order system by introducing the two-component velocity vector 
${\bf v}(x,y,t)=(v_x(x,y,t),v_y(x,y,t))$:  
\begin{eqnarray}
\frac{\partial {\bf v}}{\partial t}&=& - \nabla u, \\
\frac{\partial u}{\partial t}&=& - c^2\nabla \cdot {\bf v}.    
\end{eqnarray}

By using the complex coordinate stretching \cite{ber94,chewee94,liutao97}, we have the following 
equations: 
\begin{eqnarray}
\frac{\partial v_{\eta}}{\partial t}+\omega_{\eta}v_{\eta} &=& - \frac{\partial u}{\partial \eta}, \label{pmlv}\\
\frac{\partial u^{(\eta)}}{\partial t}+\omega_{\eta}u^{(\eta)} &=& - c^2\frac{\partial v_{\eta}}{\partial \eta}, \label{pmlu}   
\end{eqnarray}
where $\eta=x,y$, $u=u^{(x)}+u^{(y)}$, and $\omega_{\eta}$ represents a loss in the PML and is zero 
in a regular non-PML region. 

In our simulation, we will take the following $\omega_{\eta}: [0,1]\rightarrow R^{+}$:
\begin{eqnarray}
\omega_{\eta}(s)&=& \left\{\begin{array}{cc}
\frac{b}{\sigma}\left(\frac{s-\sigma}{\sigma}\right)^2,\quad\quad & s\in [0,\sigma];  \\
 0,\quad\quad & s\in [\sigma, 1-\sigma]; \\
\frac{b}{\sigma}\left(\frac{s-1+\sigma}{\sigma}\right)^2,\quad\quad & s\in [1-\sigma,1],
\end{array}\right.
\nonumber
\end{eqnarray}
where $\sigma>0$ and $b>0$ are some appropriate constants. Depending on the computational domain, this function will 
be rescaled accordingly. 

Equations \eqref{pmlv} and \eqref{pmlu} are discretized by a staggered finite-difference scheme, and the details 
are omitted; see \cite{liutao97}. 

\subsection{$T_0$ and the time reversal operator $A$}
To compute the critical time $T_0$, we first solve the following eikonal equation, 
\begin{eqnarray}
c |\nabla \tilde{T}| &=& 1, \\
T|_{\Gamma}&=&0,
\end{eqnarray}
by using the fast sweeping finite-difference scheme as designed in \cite{zha05,kaooshqia04,qiazhazha07}, 
where $\Gamma$ can be the whole boundary $\partial{\Omega}$ or a part of $\partial{\Omega}$. 
Then $T_0 = \max\{\tilde{T}\}$. 

The Poisson equation defining the harmonic extension is solved by a V-cycle multigrid method \cite{bra77}, 
which converges with the optimal rate independent of the mesh size. 

The time-reversal wave equation \eqref{2} with terminal and boundary conditions is solved 
by a standard second-order finite-difference time-domain scheme.

\section{Numerical results: smooth speed, data on the whole $\bo$}
\label{numerics}

In this section we show two dimensional numerical examples to validate our algorithms.

\begin{figure}[h]
\begin{center}
(a) \includegraphics[width=4.8cm, height=4.59552715654952cm,keepaspectratio]{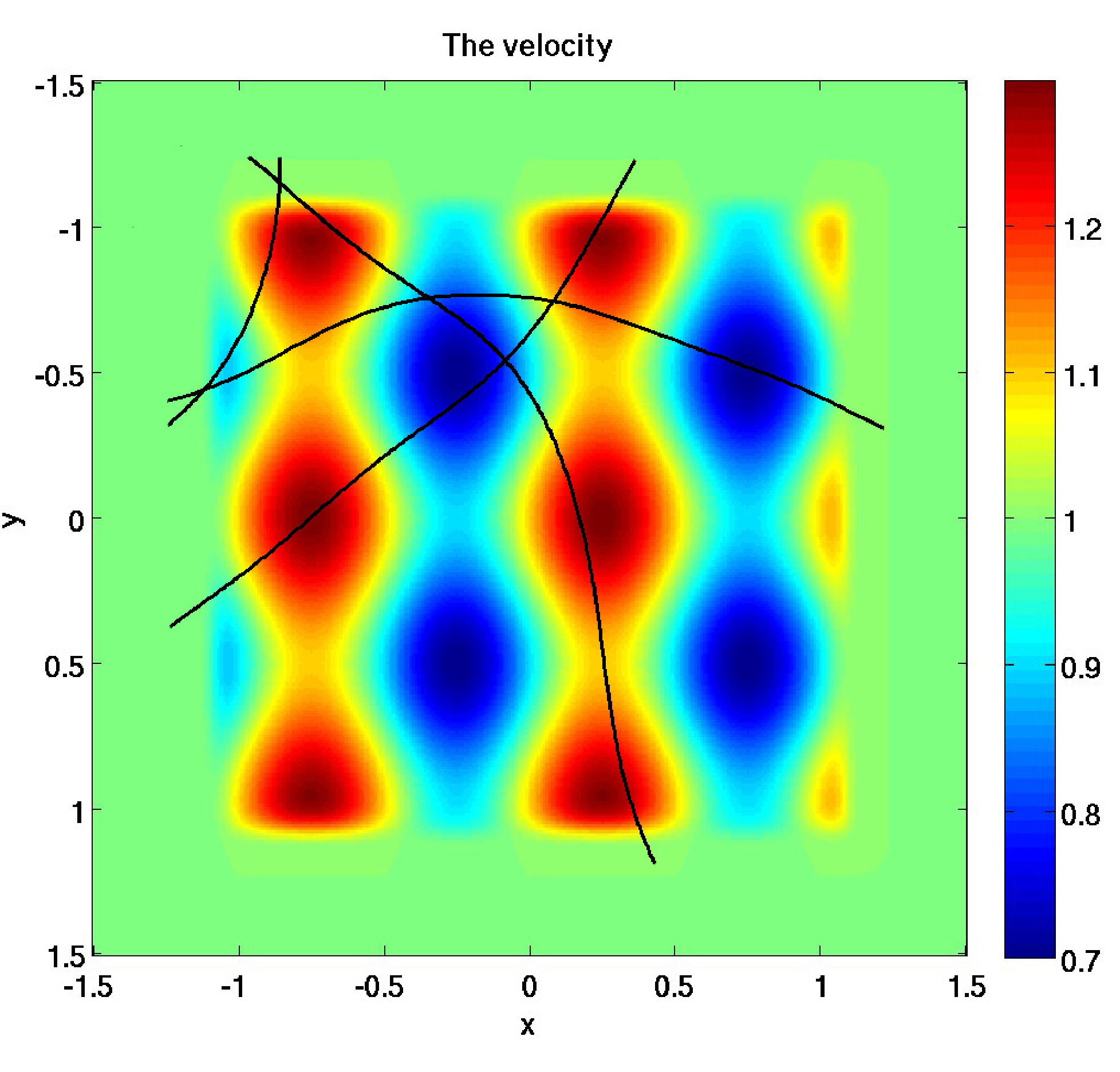}
(b) 
\includegraphics[width=4.8cm,height=4.59552715654952cm,keepaspectratio]{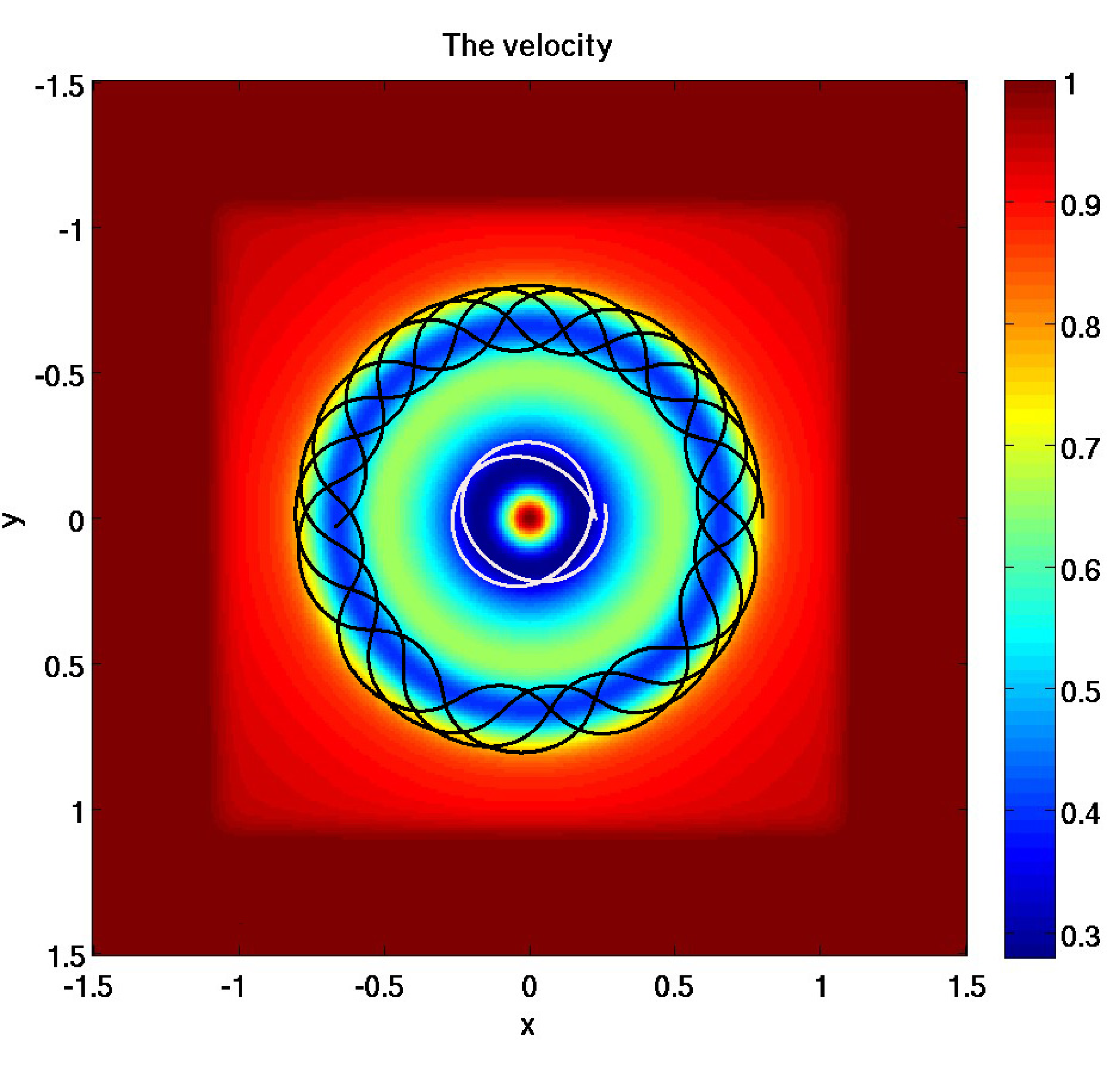}
(c) \includegraphics[width=4.8cm,height=4.59552715654952cm, keepaspectratio]{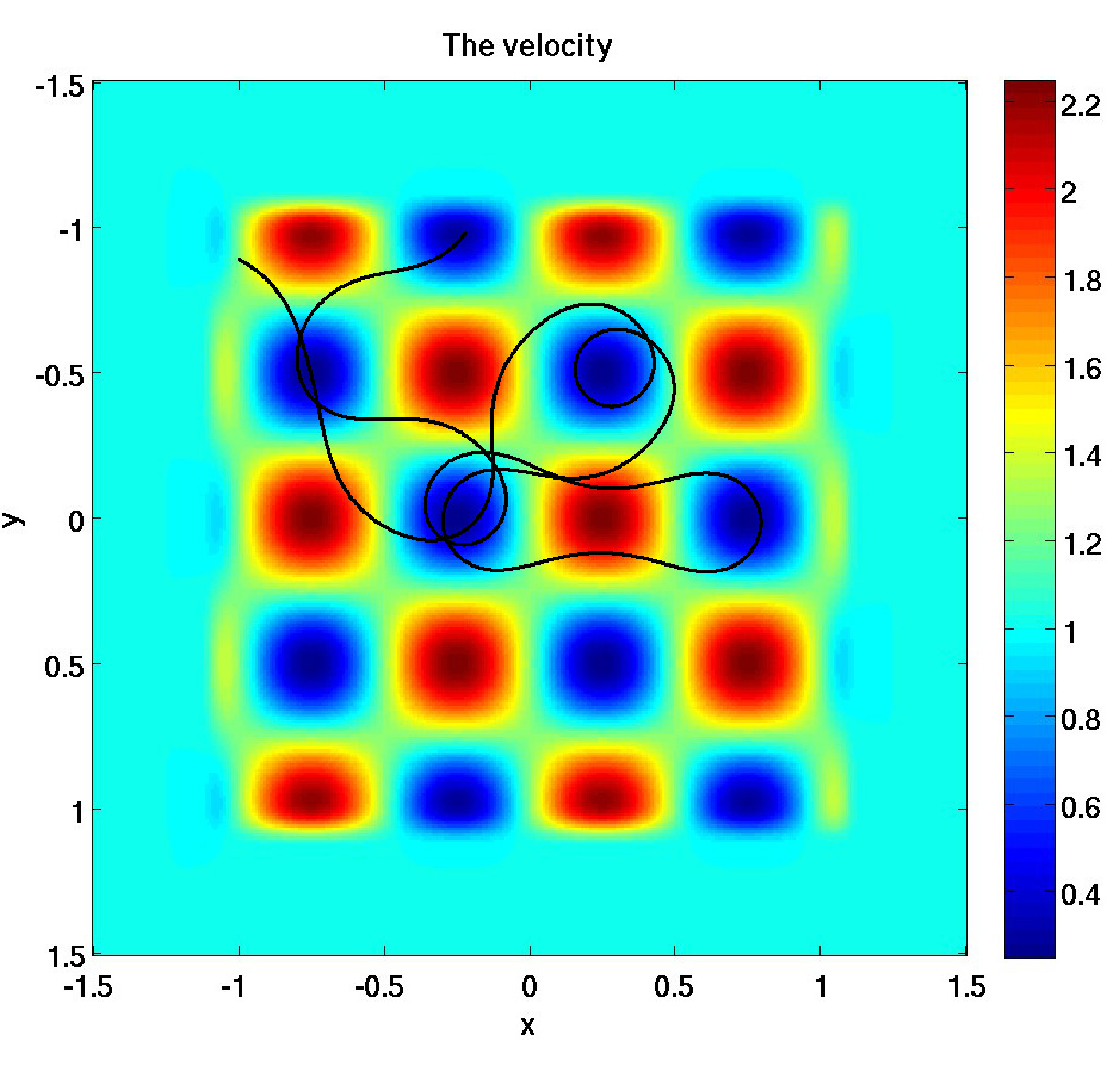} 
\end{center}
\caption{Sound speed models.  (a): the variable  non-trapping sound speed $c_1$. 
(b): the variable radial trapping sound speed $c_2$. (c): the variable  trapping speed $c_3$.}
\label{Fig:Velocities}
\end{figure}

Our computational domain is taken to be $[-1.5,1.5]^2$, and $\Omega=[-1.28,1.28]^2$. We work with three sound speeds:
\begin{eqnarray}
c_1(x,y)&=&1.0+0.2\sin(2\pi x)+0.1\cos(2\pi y), \label{key8} \\
c_2(x,y)&=&\frac{9(x^2+y^2)}{1+9(x^2+y^2)}+\exp{(-90(x^2+y^2))}-0.4\exp{\left(-10(3\sqrt{(x^2+y^2)}-2)^2\right)}, \label{key63} \\
c_3(x,y)&=&1.25 + \sin(2.0\pi x)\cos(2\pi y). \label{key6} 
\end{eqnarray}
See Figure \ref{Fig:Velocities}, with some geodesics plotted.  We also use a smooth cutoff for a smooth transition of those speeds to $1$ near $\bo$ and outside $\Omega$.

In all examples below, we use the abbreviation NS for the Neumann series, and TR for time reversal.  We take a finite 
number, $k+1$, of terms in the Neumann series expansion by stopping when the error gets below $5\%$ in the stable cases, 
and somewhat higher in the non-stable ones. Taking more terms, especially in the cases where we have stability, improves 
the error several times but we do not show those examples. Next, we only compute and comment on the $L^2$ error since 
all test images model non $H^1$ functions. We also show a diagram of the distance function, $\dist(x,\bo)$. 

\subsection{Example 1: The Shepp-Logan phantom}
\begin{figure}
\begin{center}
(a)\epsfig{figure=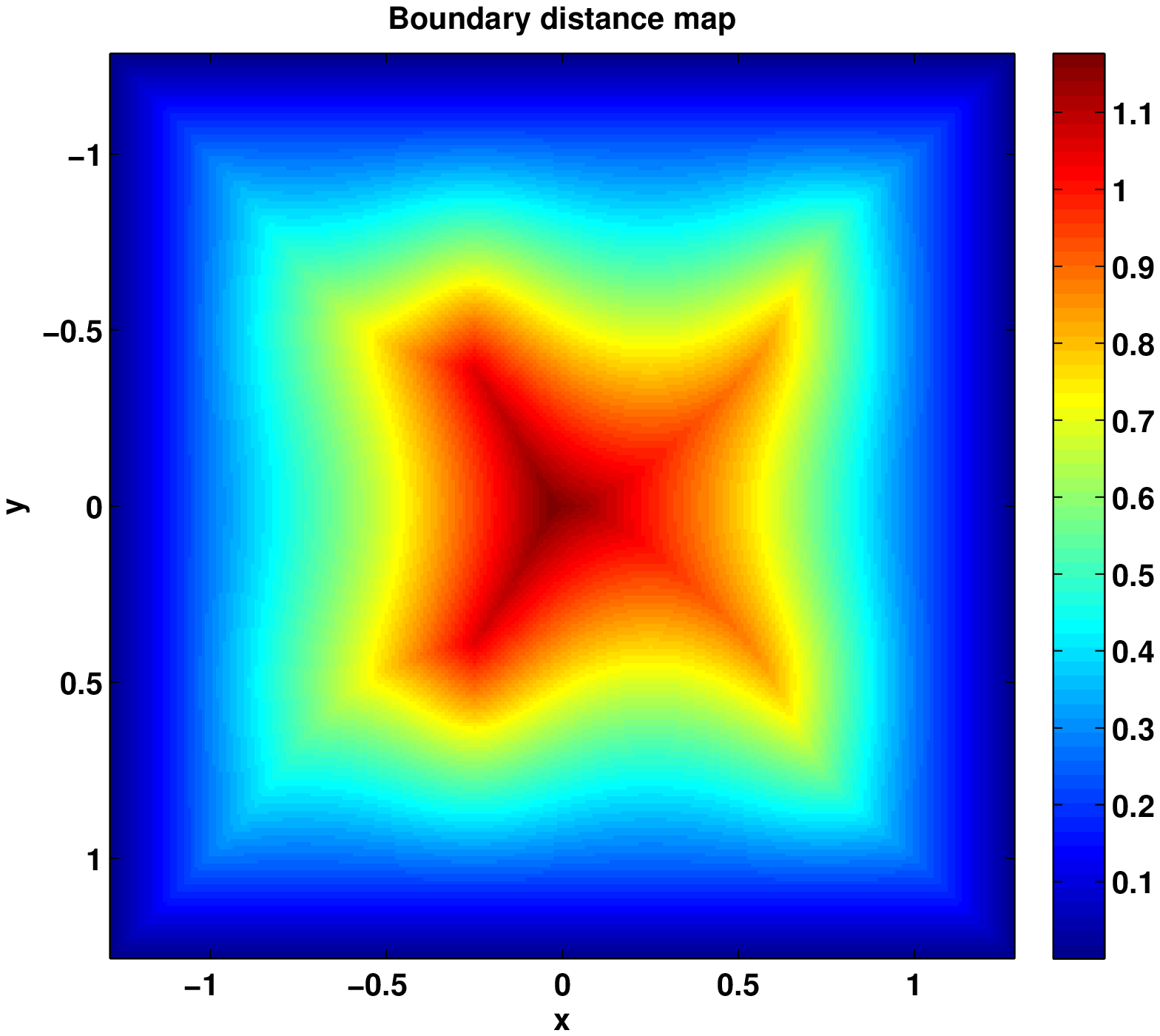,height=5.50cm}
(b)\epsfig{figure=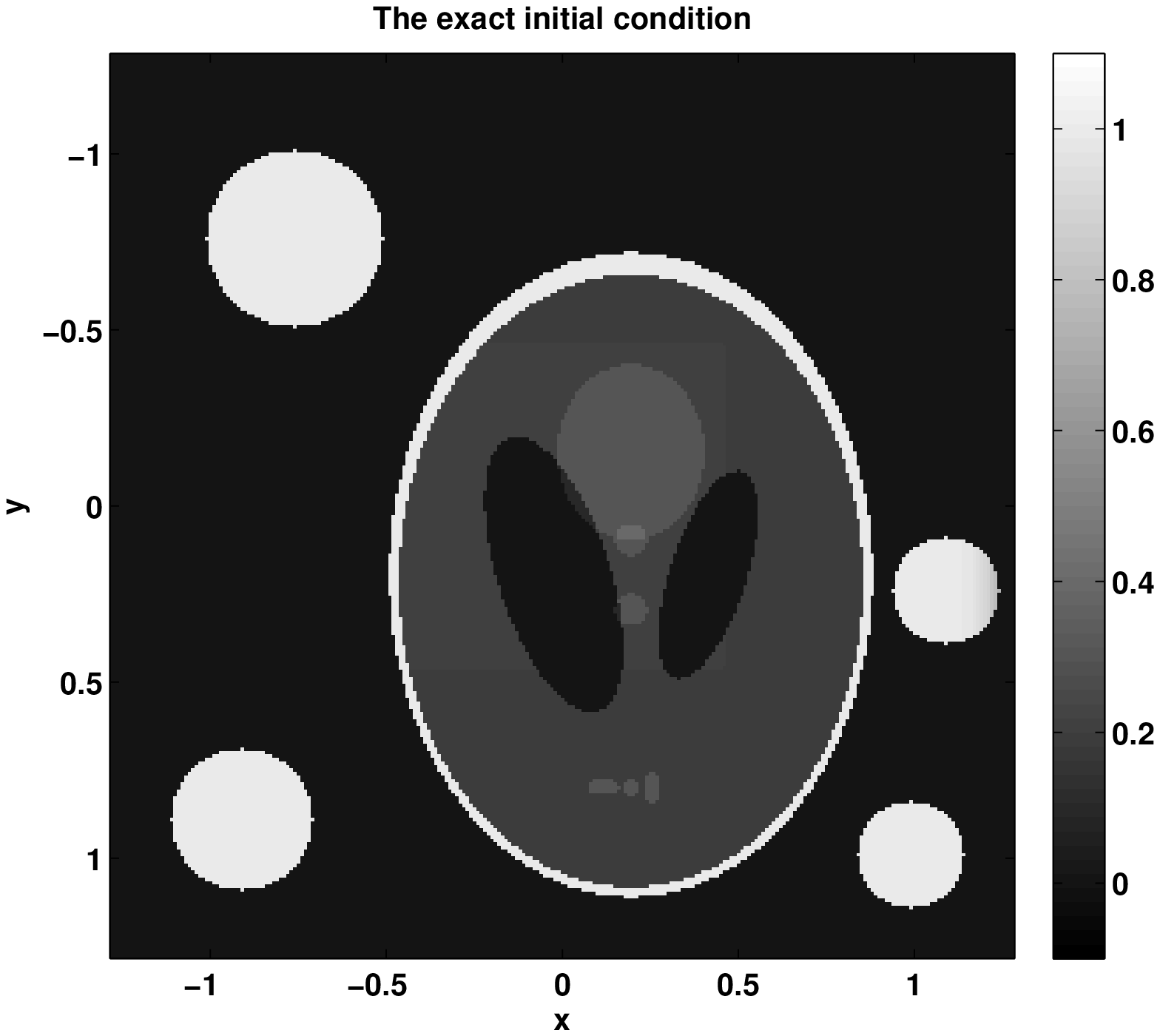,height=5.50cm}\\
(c)\epsfig{figure=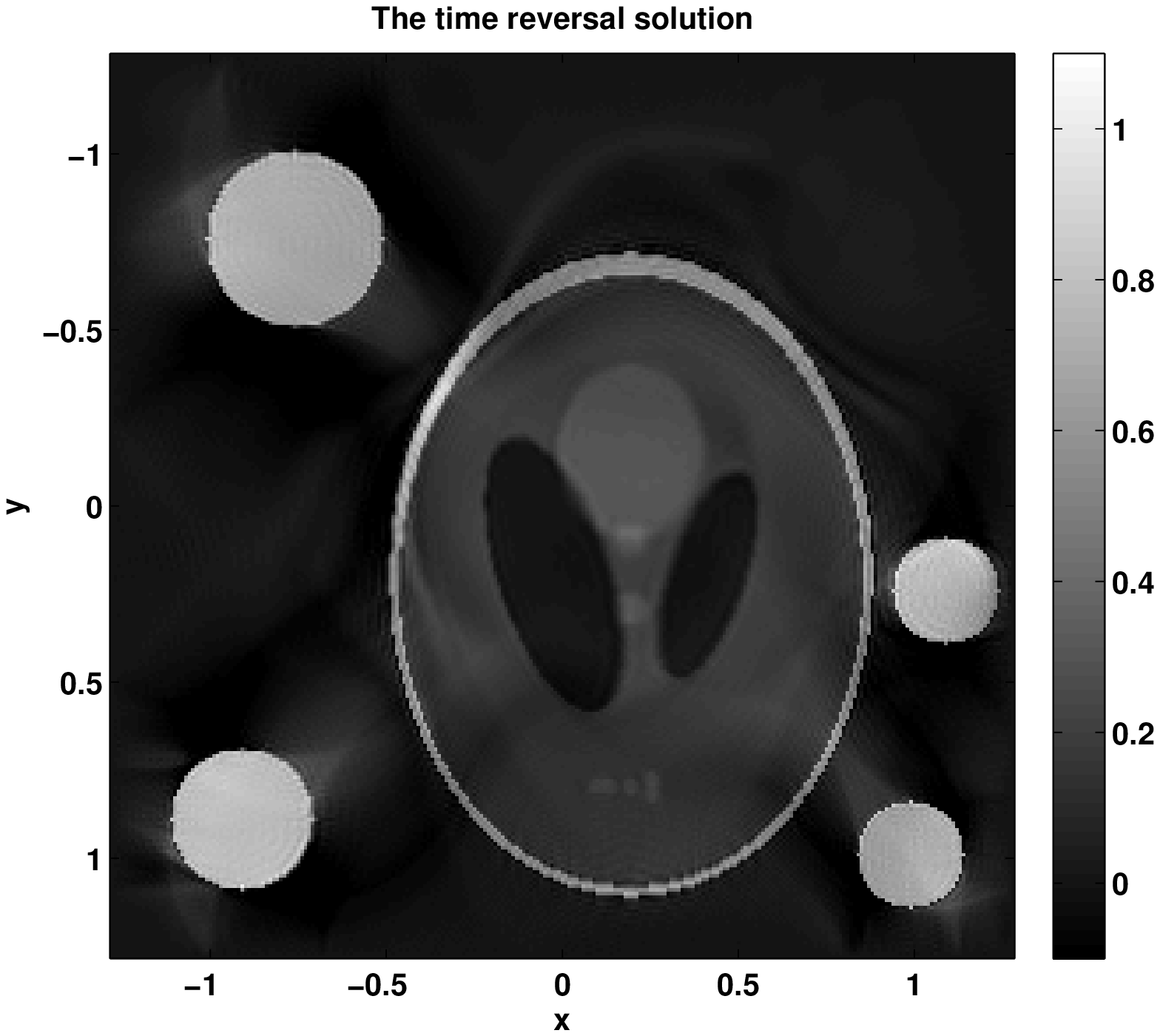,height=5.50cm}
(d)\epsfig{figure=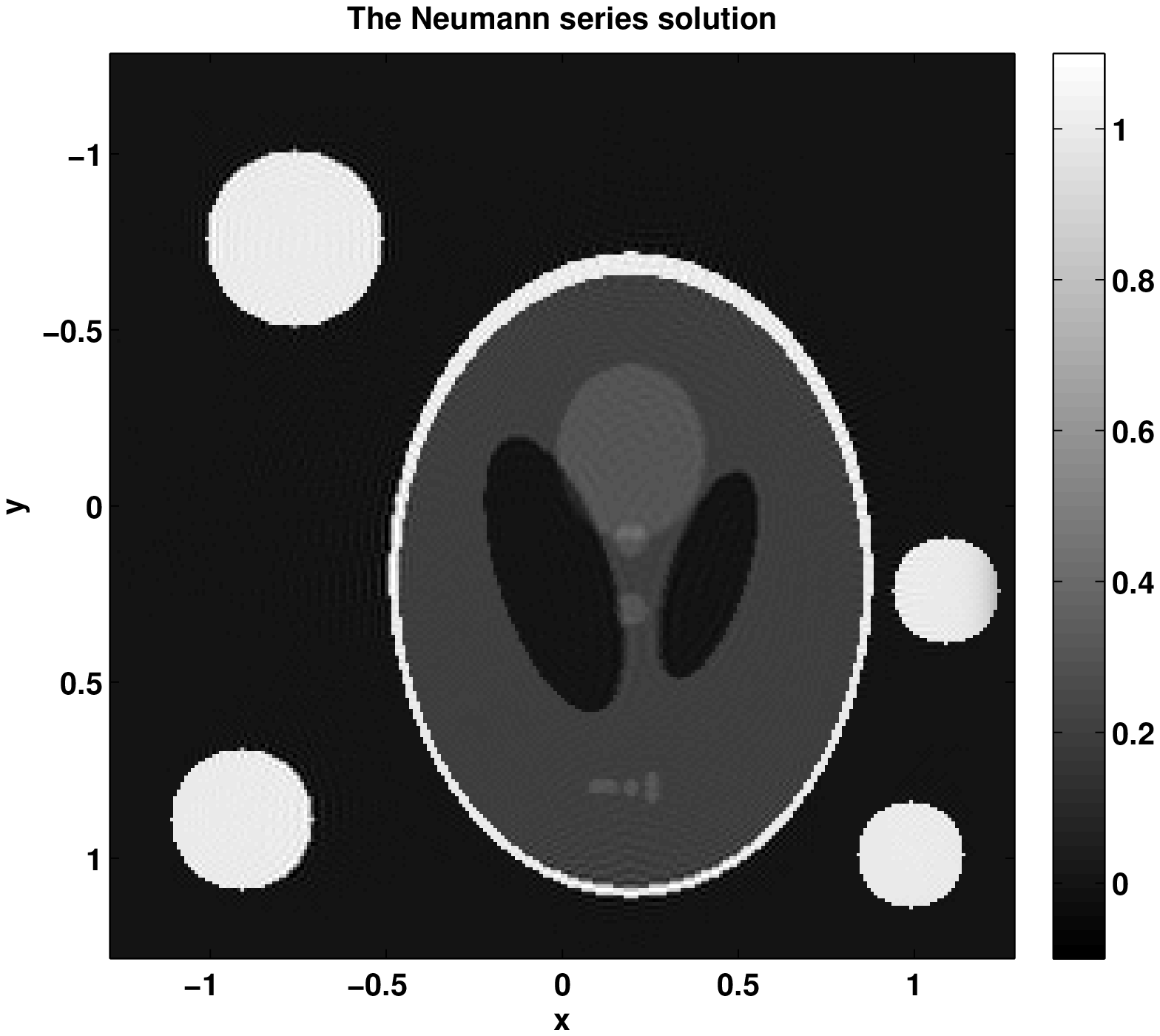,height=5.50cm}\\
(e)\epsfig{figure=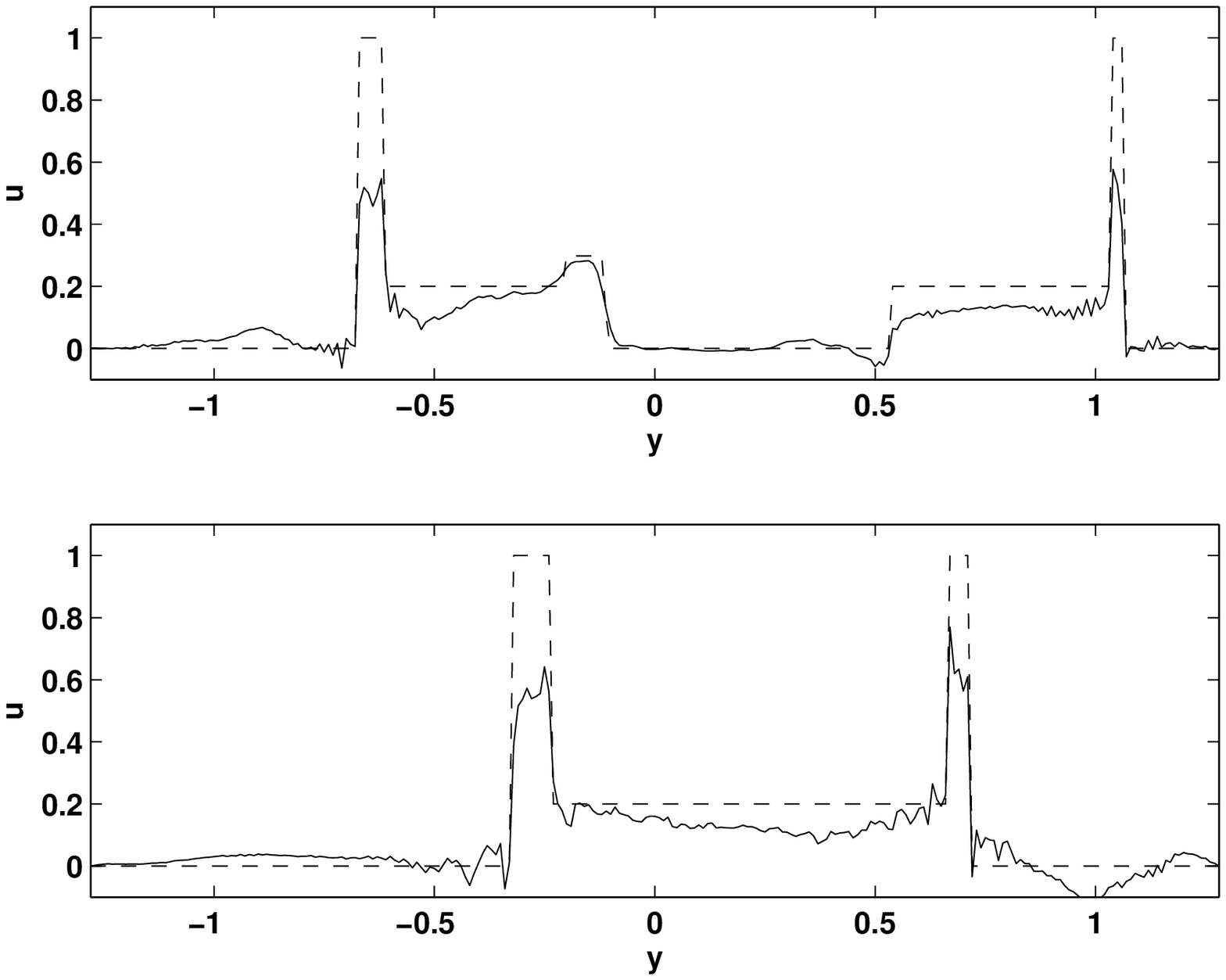,height=4cm,width=6.0cm}
(f)\epsfig{figure=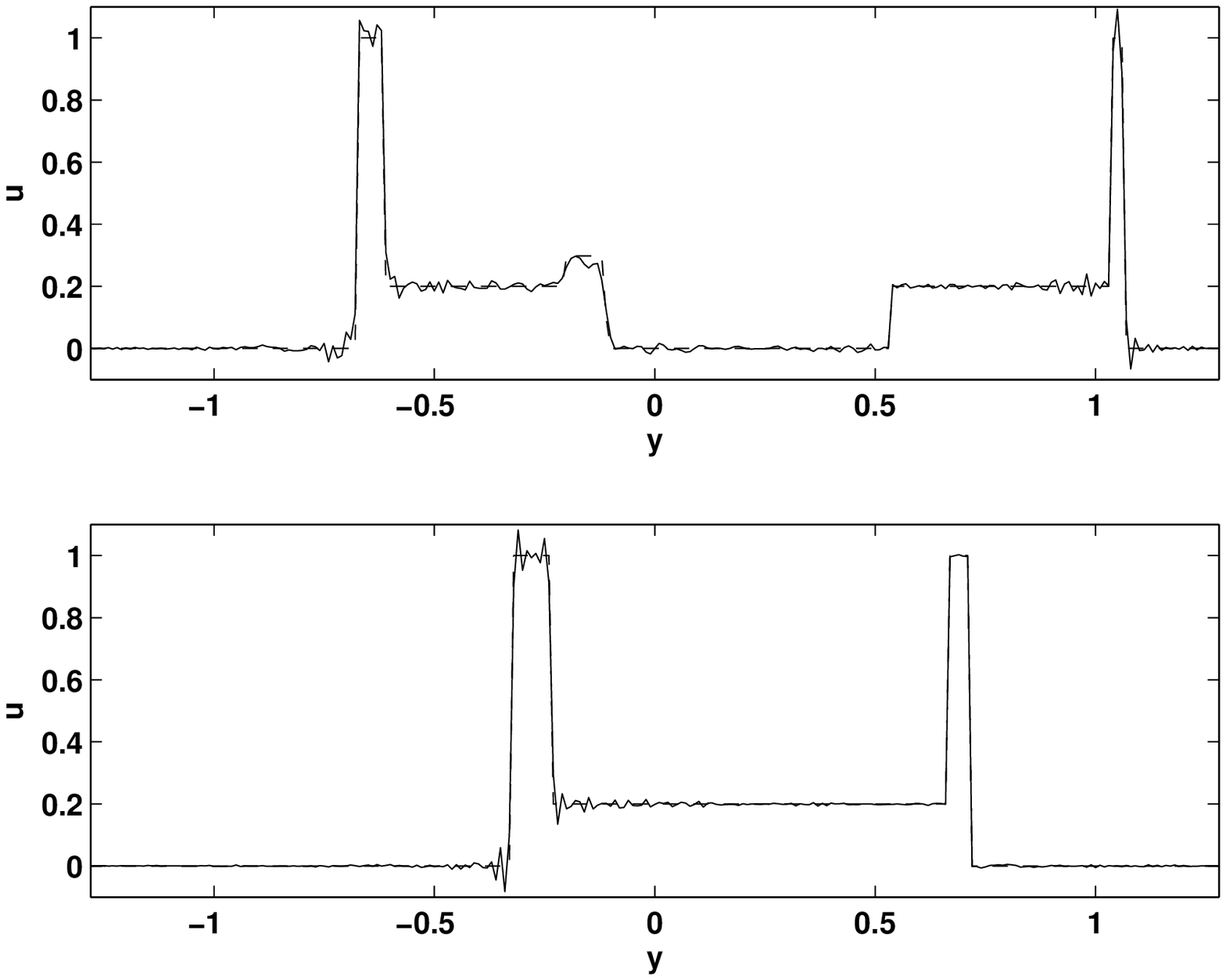,height=4cm,width=6.0cm}\\
(g)\epsfig{figure=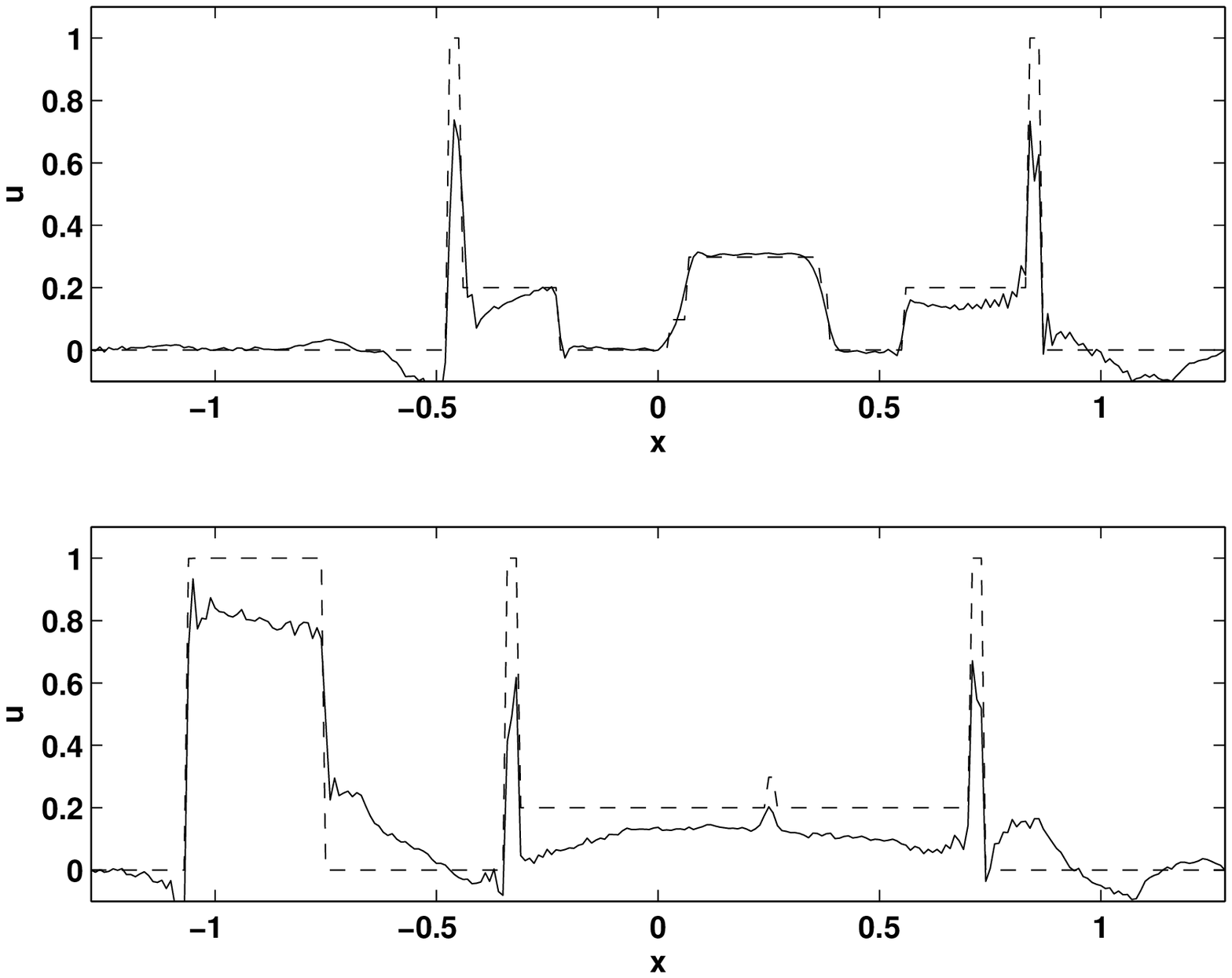,height=4cm,width=6.0cm}
(h)\epsfig{figure=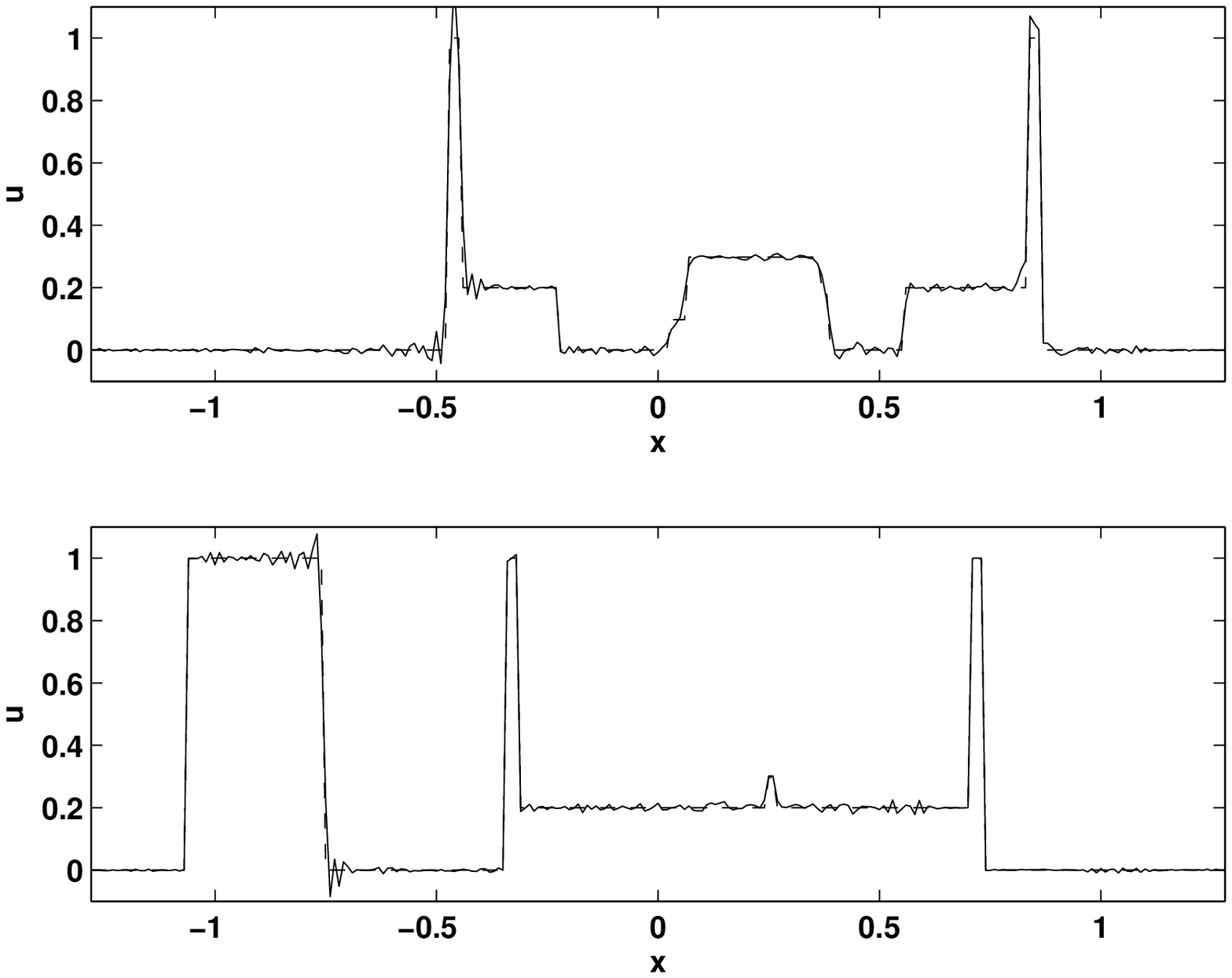,height=4cm,width=6.0cm}
\end{center}
\caption{Example 1 with the non-trapping speed $c_1$. Case 2: $T=4T_0$. 
(a): the boundary distance map. 
(b): the exact initial condition. 
(c): the time reversal solution. 
(d): the Neumann series solution. 
(e): $x$-slices of the time reversal solution (continuous line) and the exact solution (a dashed line). 
(f): $x$-slices of the Neumann series solution (continuous line) and the exact solution (a dashed line). 
(g): $y$-slices of the time reversal solution (continuous line) and the exact solution (a dashed line).
(h): $y$-slices of the Neumann series solution (continuous line) and the exact solution (a dashed line).}
\label{Fig:2dwaveT2pNotrap}
\end{figure}

\begin{figure}
\begin{center}
(a)\epsfig{figure=Fig/Key8Keyini6943Sides100Traveltime.ps,height=5.50cm }
(b)\epsfig{figure=Fig/Key8Keyini6943NX301InitialTarget.ps,height=5.50cm }\\
(c)\epsfig{figure=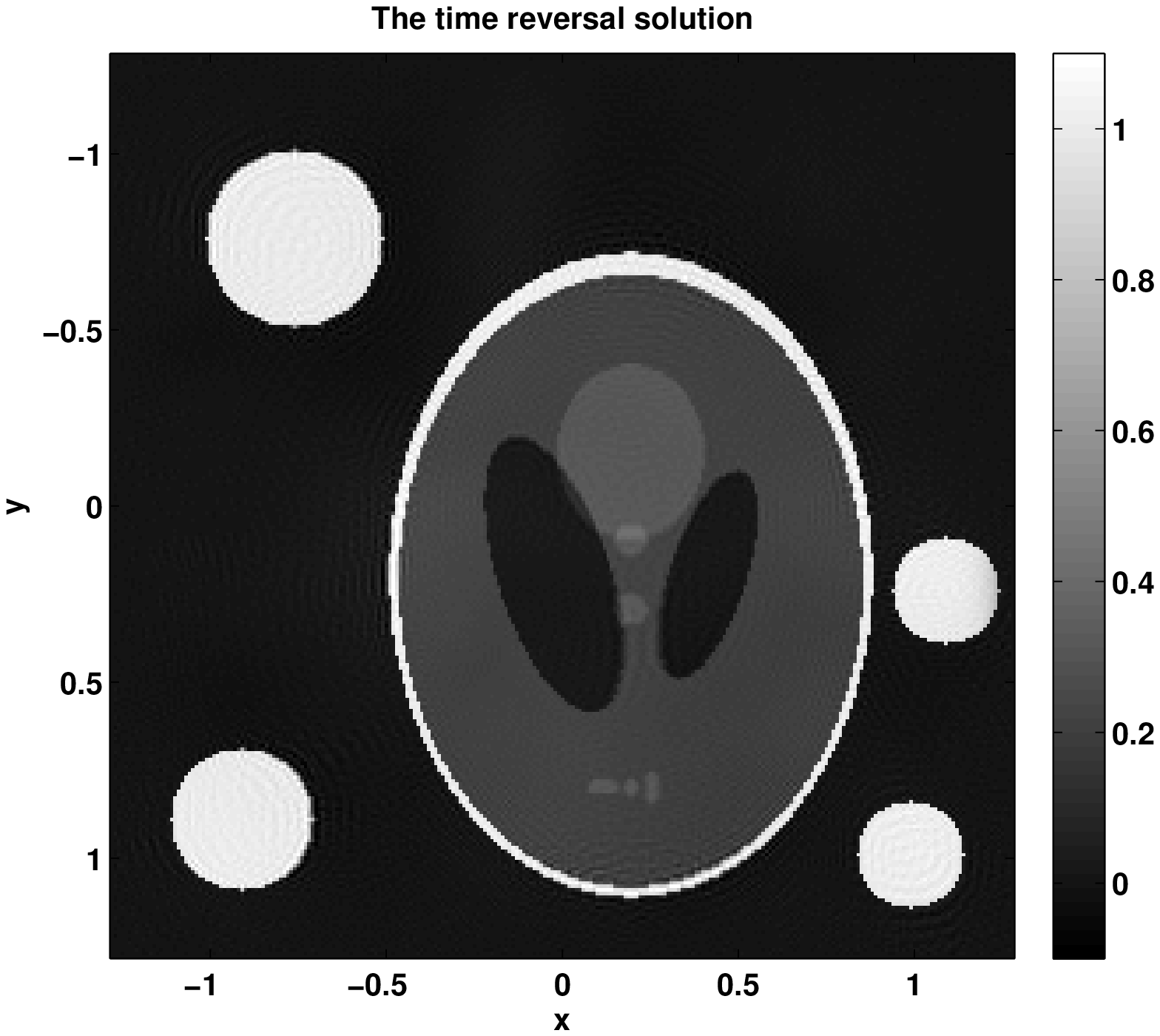,height=5.50cm }
(d)\epsfig{figure=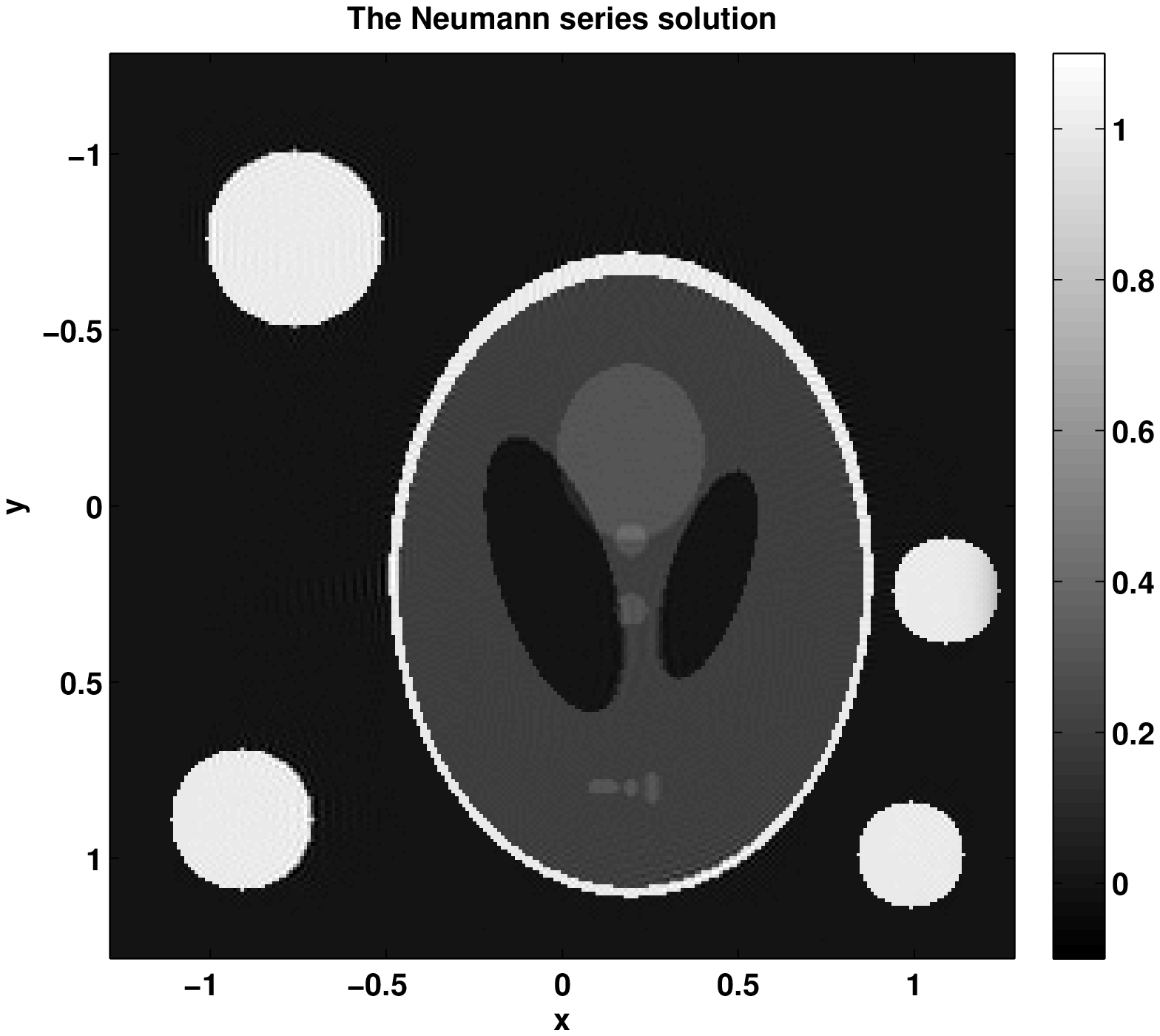,height=5.50cm }\\
(e)\epsfig{figure=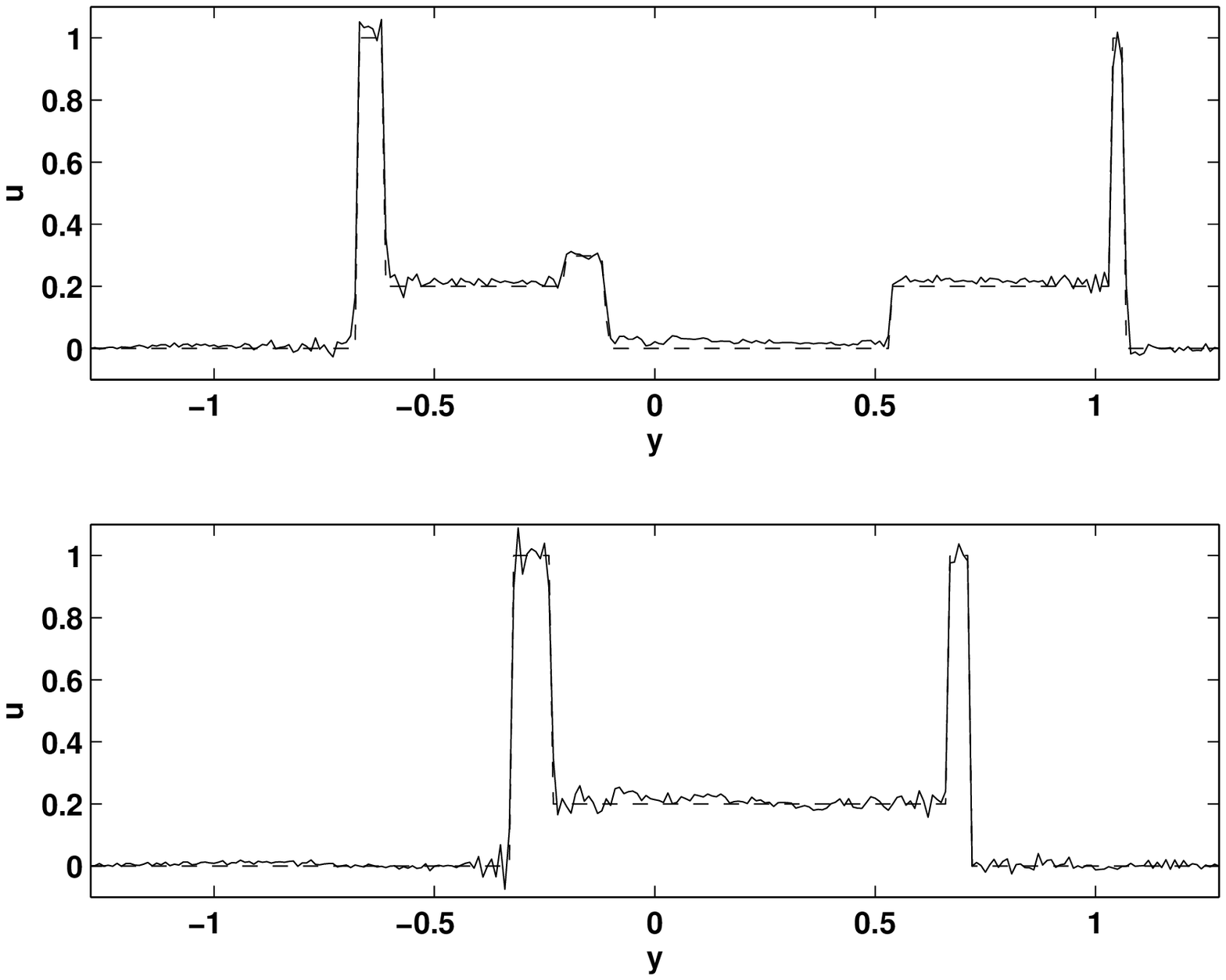,height=4cm,width=6.0cm}
(f)\epsfig{figure=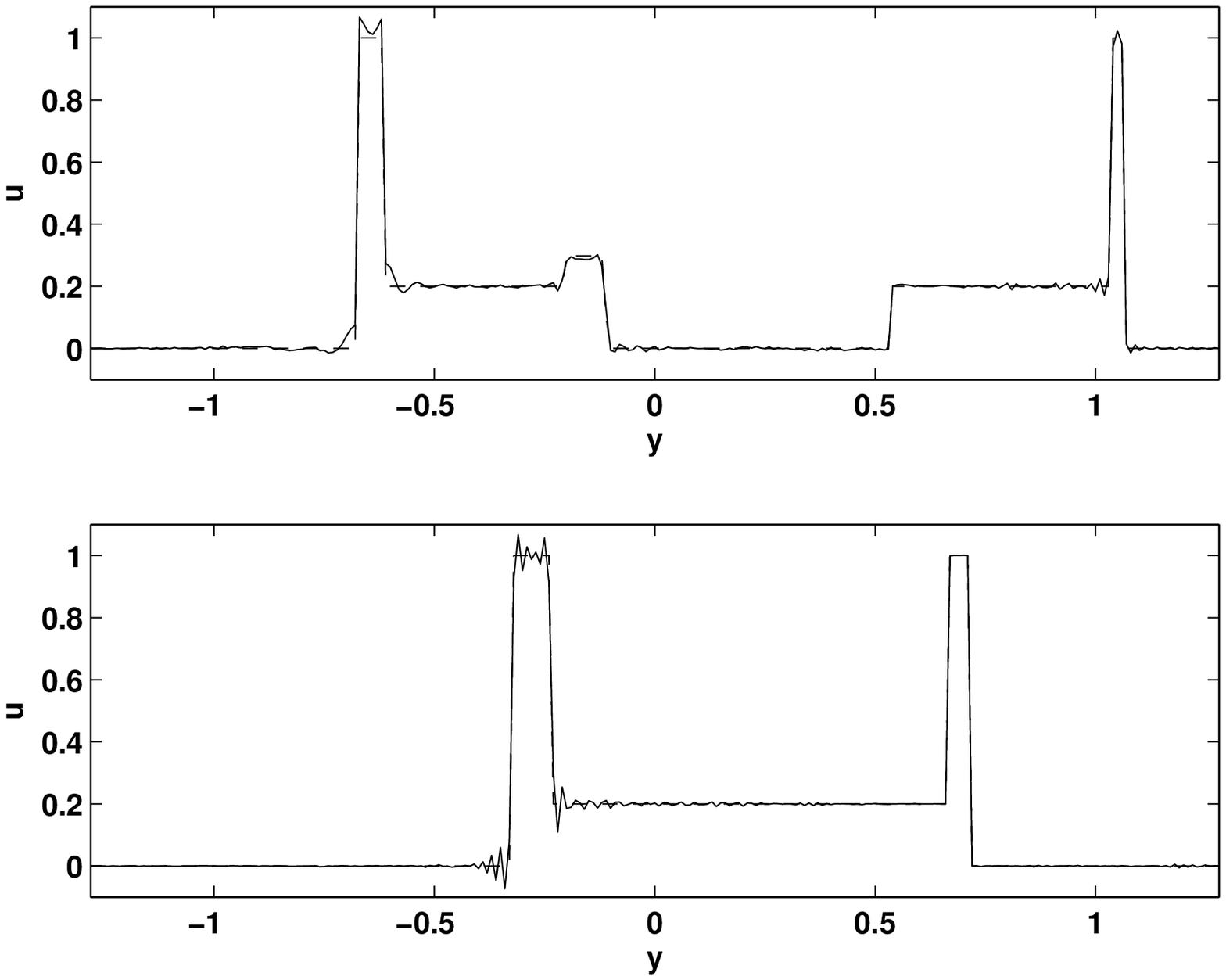,height=4cm,width=6.0cm}\\
(g)\epsfig{figure=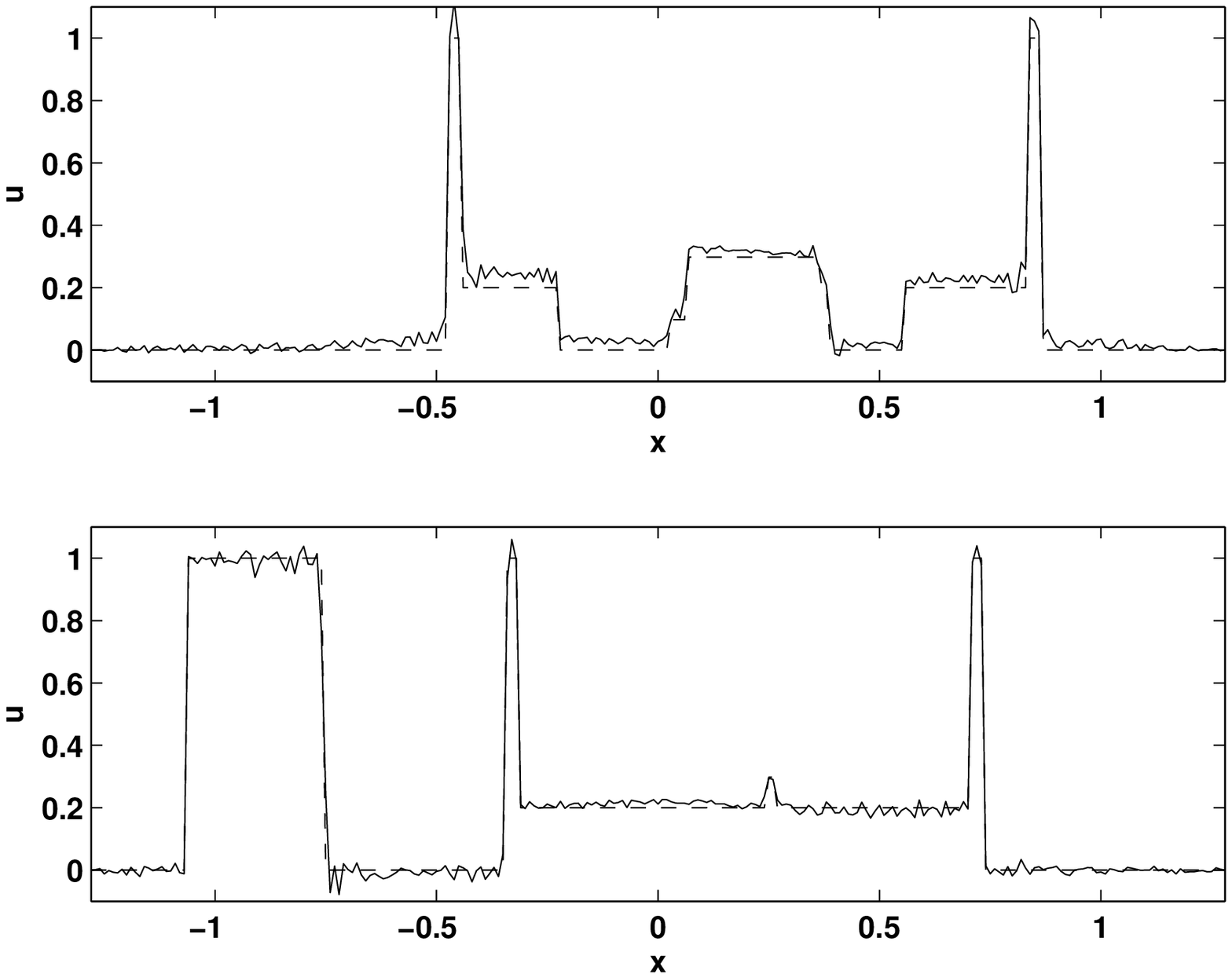,height=4cm,width=6.0cm}
(h)\epsfig{figure=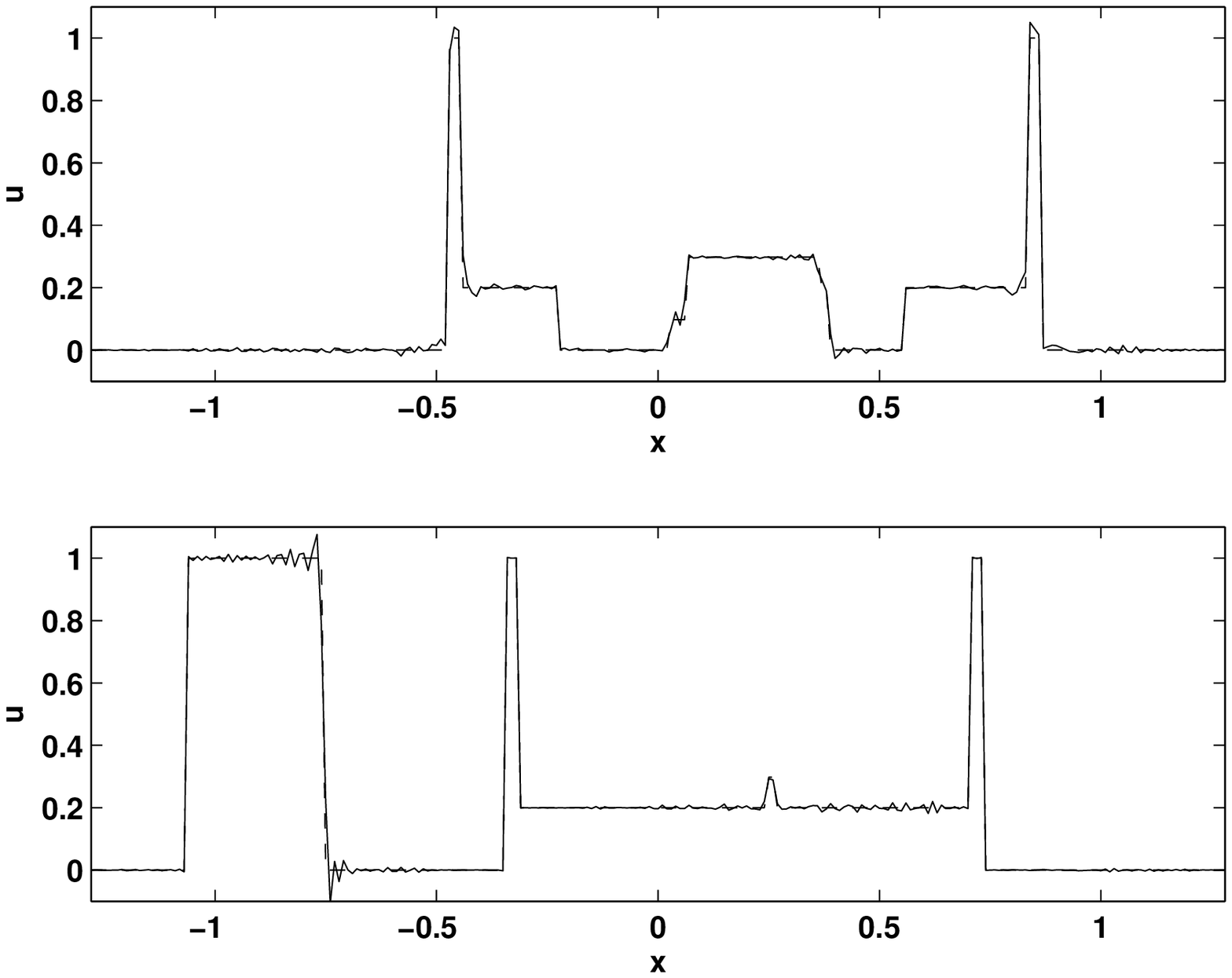,height=4cm,width=6.0cm}
\end{center}
\caption{Example 1 with the non-trapping speed $c_1$. Case 2: $T=4T_0$. 
(a): the boundary distance map. 
(b): the exact initial condition. 
(c): the time reversal solution. 
(d): the Neumann series solution. 
(e): $x$-slices of the time reversal solution (continuous line) and the exact solution (a dashed line). 
(f): $x$-slices of the Neumann series solution (continuous line) and the exact solution (a dashed line). 
(g): $y$-slices of the time reversal solution (continuous line) and the exact solution (a dashed line).
(h): $y$-slices of the Neumann series solution (continuous line) and the exact solution (a dashed line).}
\label{Fig:2dwaveT4pNotrap}
\end{figure}

\begin{figure}
\begin{center}
(a)\epsfig{figure=Fig/Key8Keyini6943Sides100Traveltime.ps,height=5.50cm }
(b)\epsfig{figure=Fig/Key8Keyini6943NX301InitialTarget.ps,height=5.50cm }\\
(c)\epsfig{figure=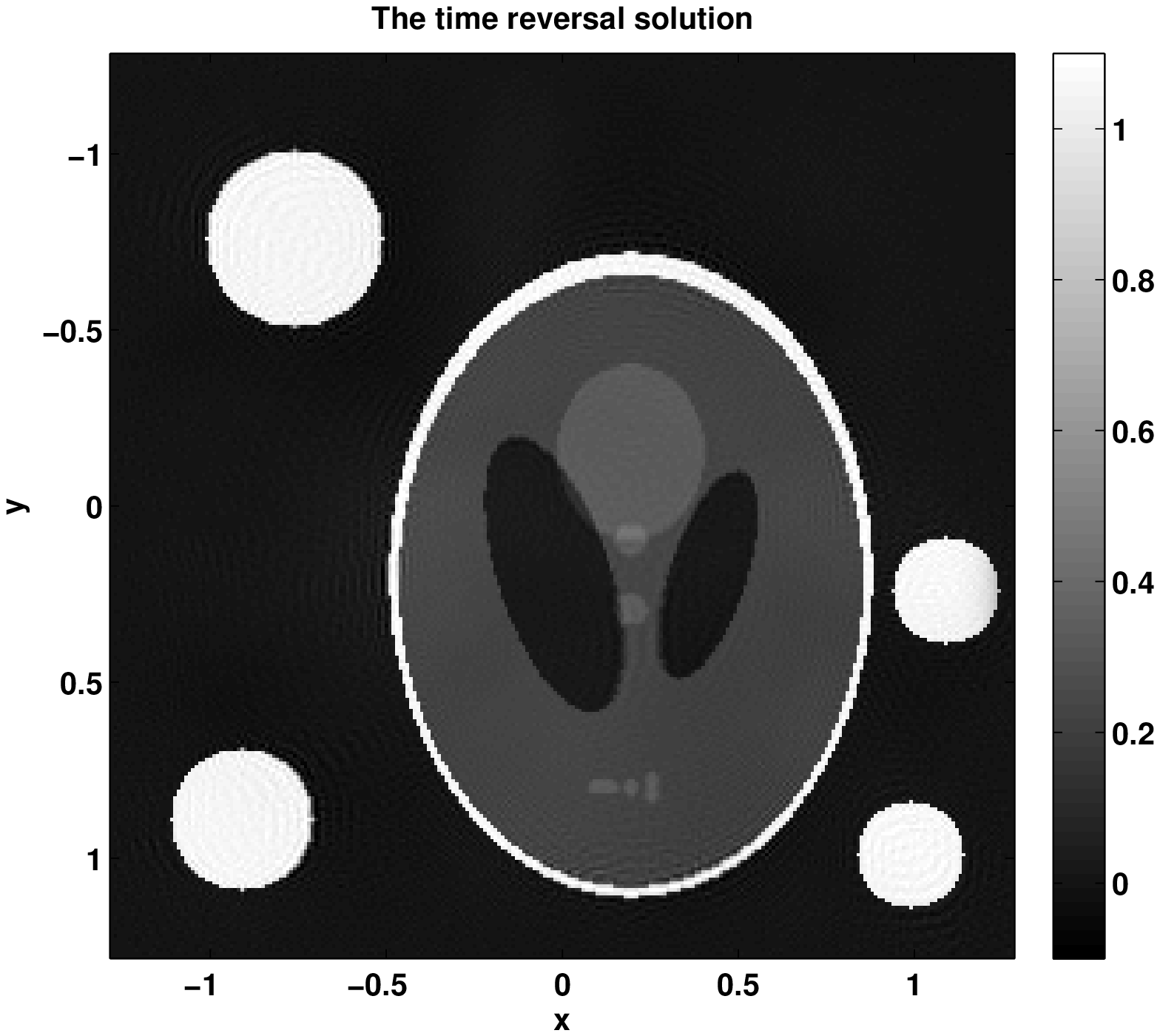,height=5.50cm }
(d)\epsfig{figure=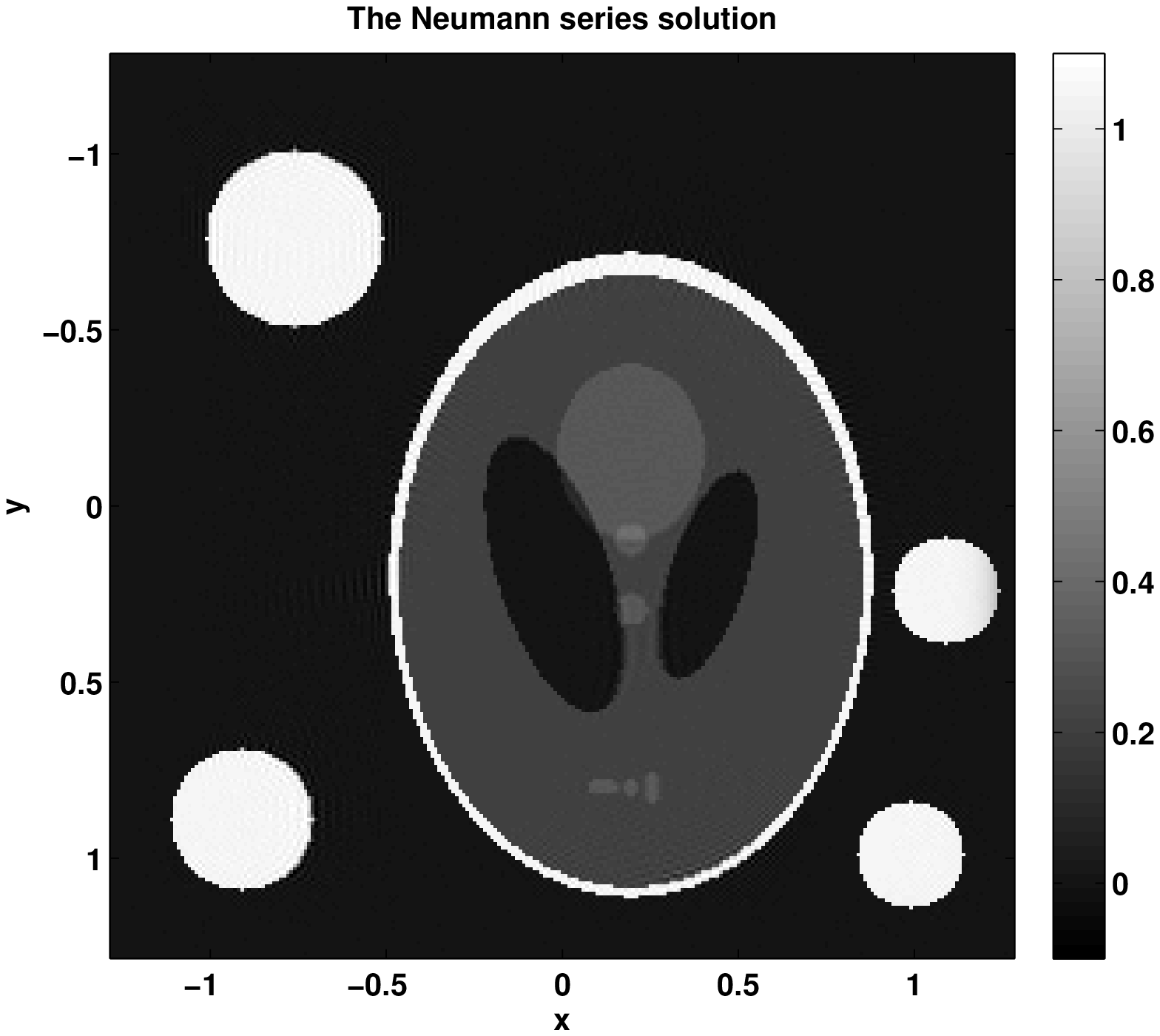,height=5.50cm }\\
(e)\epsfig{figure=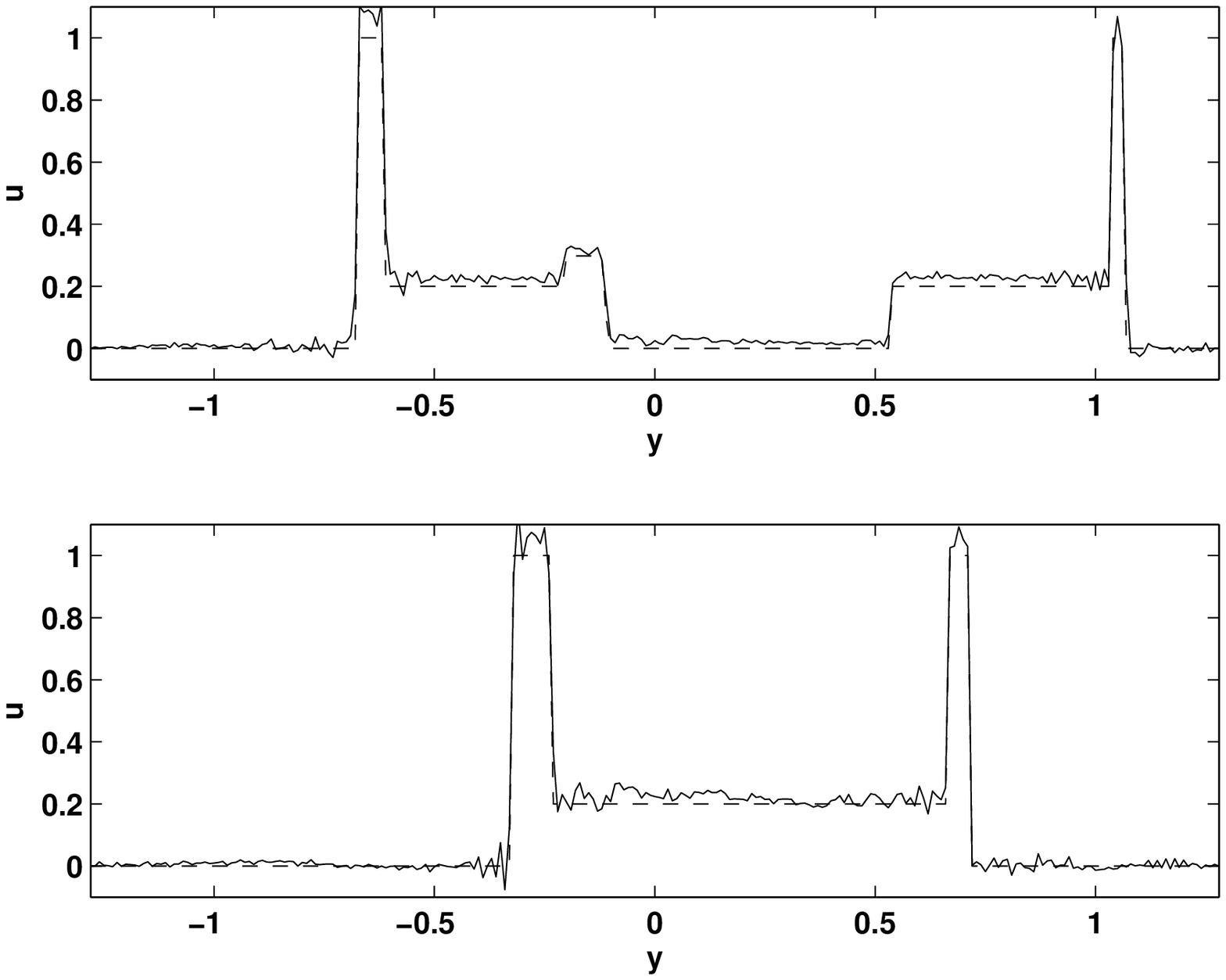,height=4cm,width=6.0cm}
(f)\epsfig{figure=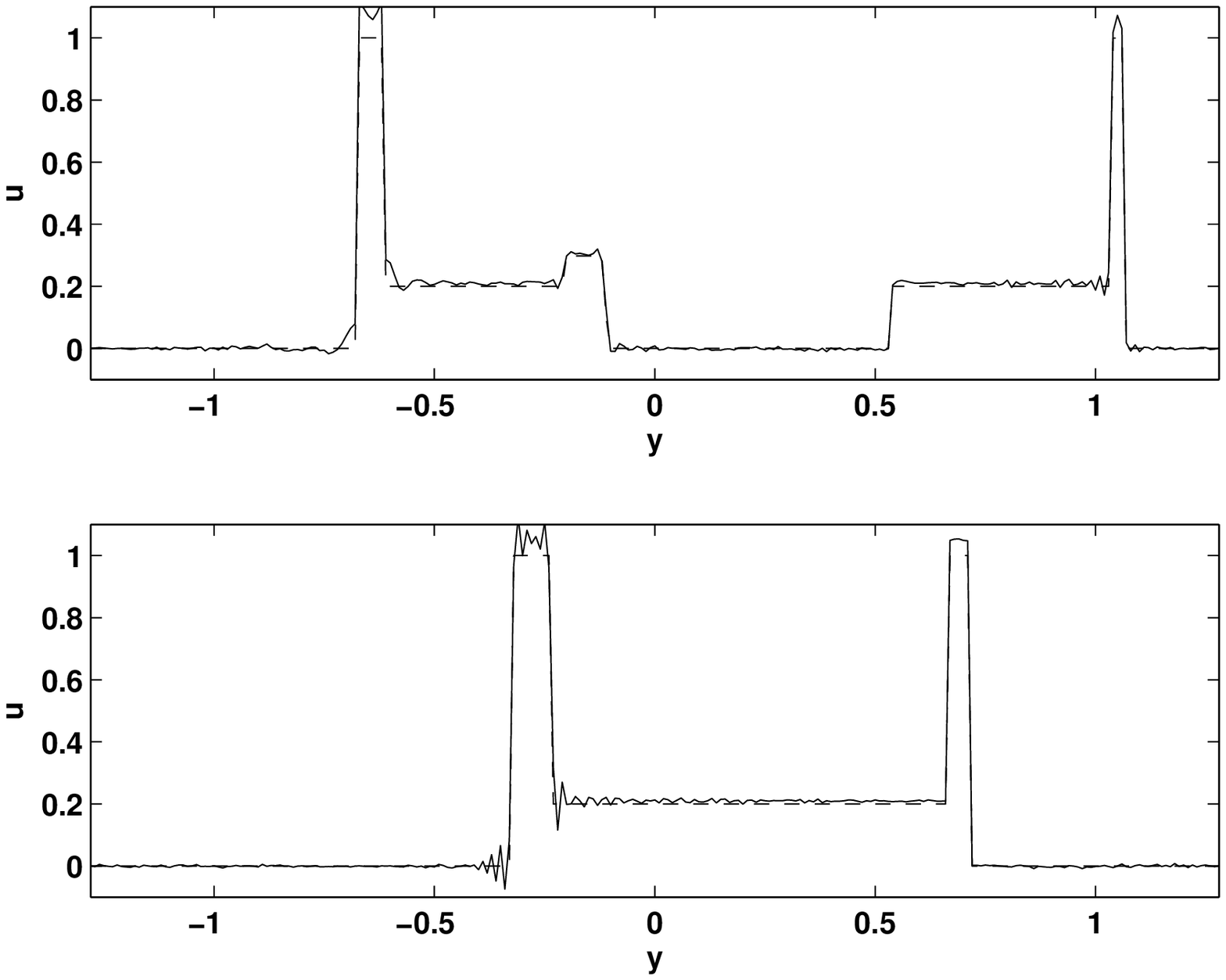,height=4cm,width=6.0cm}\\
(g)\epsfig{figure=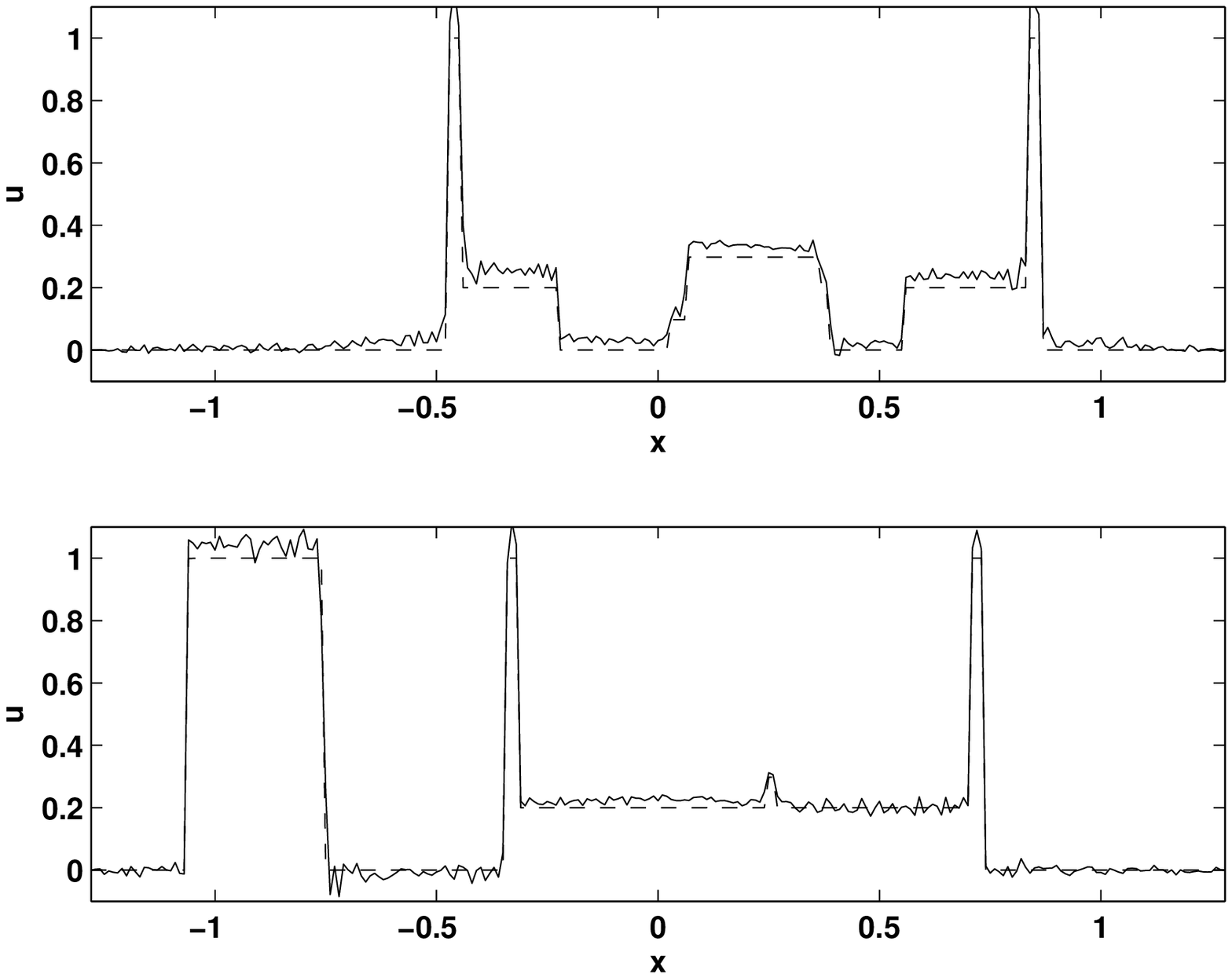,height=4cm,width=6.0cm}
(h)\epsfig{figure=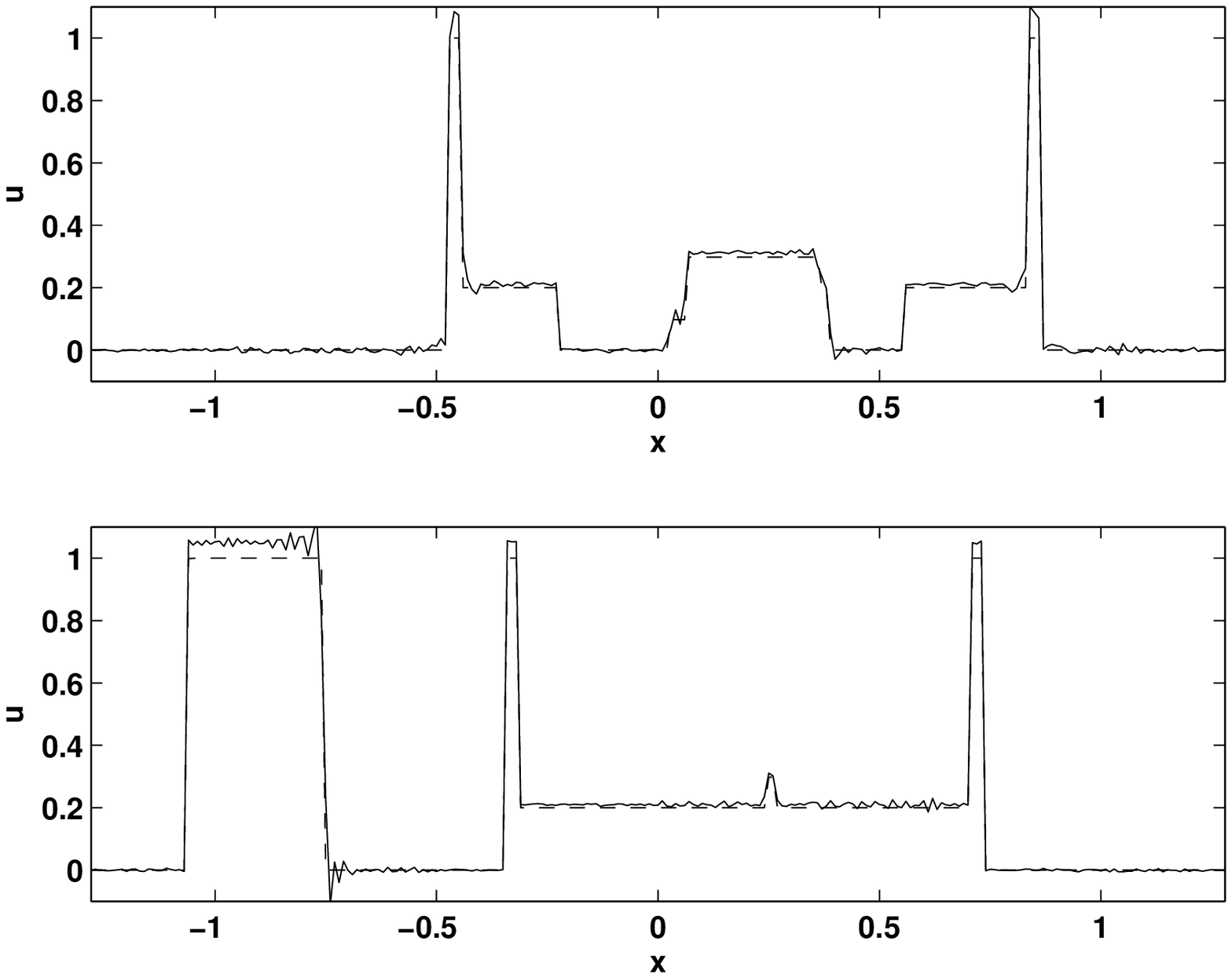,height=4cm,width=6.0cm}
\end{center}
\caption{Example 1 with the non-trapping speed $c_1$. Case 2: $T=4T_0$. 
(a): the boundary distance map. 
(b): the exact initial condition. 
(c): the time reversal solution. 
(d): the Neumann series solution. 
(e): $x$-slices of the time reversal solution (continuous line) and the exact solution (a dashed line). 
(f): $x$-slices of the Neumann series solution (continuous line) and the exact solution (a dashed line). 
(g): $y$-slices of the time reversal solution (continuous line) and the exact solution (a dashed line).
(h): $y$-slices of the Neumann series solution (continuous line) and the exact solution (a dashed line).}
\label{Fig:2dwaveT4pNotrapNoise}
\end{figure}

\subsubsection{Non-trapping speed $c_1$, Figures~\ref{Fig:2dwaveT2pNotrap}--\ref{Fig:2dwaveT4pNotrapNoise}} 

The sound speed is given by equation \eqref{key8}; see  Figure~\ref{Fig:Velocities} (a). The original image as the 
source function $f$ is shown in Figure~\ref{Fig:2dwaveT2pNotrap}(b). Numerically, we estimate $T_0$ to 
be $T_0\approx 1.1767$. A rough estimate of $T_1$ is $3<T_1<4$. We added four bright disks ($f=1$ there) to the classical  Shepp-Logan phantom.

\textbf{Figure~\ref{Fig:2dwaveT2pNotrap}:} $T=2T_0\approx 2.3535$.
The time $T$ is slightly above the stability threshold $T_1/2$ but below $T_1$. The error of the NS solution is $6.63\%$ with $k=8$ ($9$ terms) of the series vs.\ error  $37.76\%$ for the TR one. Since $T<T_1$, the TR solution does not recover the correct size of the jumps --- they are recovered with amplitudes ranging from $1/2$ to $1$, and for many of them, it is just $1/2$; this is clear from the slice diagrams. In contrast, the NS solution has the right amplitudes and would improve with more terms.   

\textbf{Figure~\ref{Fig:2dwaveT4pNotrap}:}
$T=4T_0$. The time $T$ is doubled, and it is greater than $T_1$.  The error of the NS solution is $4.99\%$ with $k=8$ ($9$ terms) of the series vs.\ error  $7.07\%$ for the TR one. The TR reconstruction gets better as expected. While the singularities are recovered correctly, the low frequency part in the TR reconstruction is still not well recovered as evidenced 
in the slice diagrams. For example, in Figure~\ref{Fig:2dwaveT4pNotrap}(e)(g) for the TR reconstruction there are 
oscillations in some flat regions; in contrast, in Figure~\ref{Fig:2dwaveT4pNotrap}(f)(h) for the NS reconstruction those 
oscillations are gone in the same flat regions.

\textbf{Figure~\ref{Fig:2dwaveT4pNotrapNoise}:}
 $T=4T_0$ with $10\%$ noise. The time $T$ is doubled, and it is greater than $T_1$.  The error of the NS solution is $7.07\%$ with $k=8$ ($9$ terms) of the series vs.\ error  $9.71\%$ for the TR one. The TR reconstruction gets better as expected. While the singularities are recovered correctly, the low frequency part in the TR reconstruction is still not well recovered as evidenced 
in the slice diagrams. For example, in Figure~\ref{Fig:2dwaveT4pNotrapNoise}(e)(g) for the TR reconstruction there are 
oscillations in some flat regions; in contrast, in Figure~\ref{Fig:2dwaveT4pNotrapNoise}(f)(h) for the NS reconstruction those 
oscillations are gone in the same flat regions.

\begin{figure}
\begin{center}
(a)\epsfig{figure=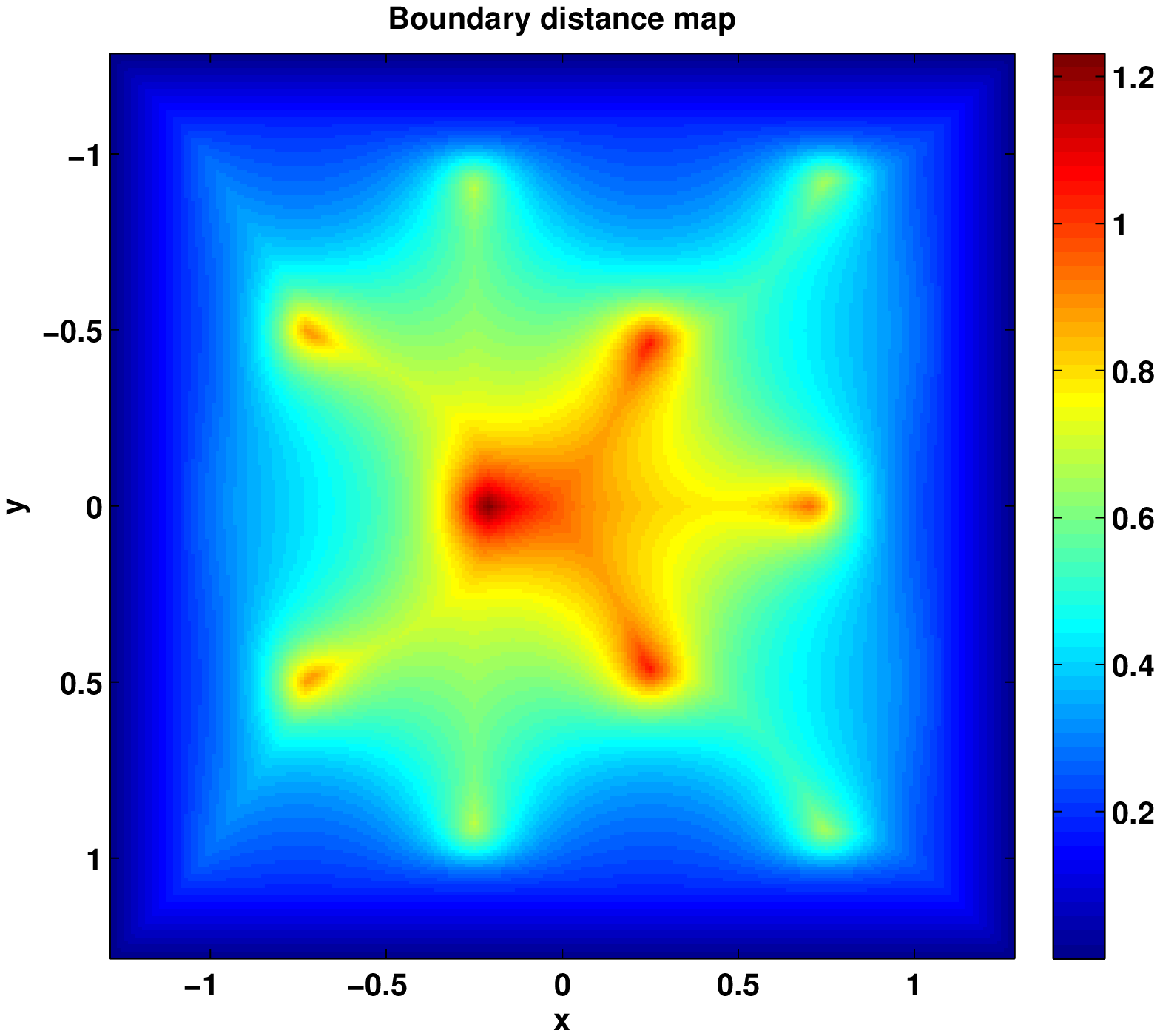,height=5.5cm }
(b)\epsfig{figure=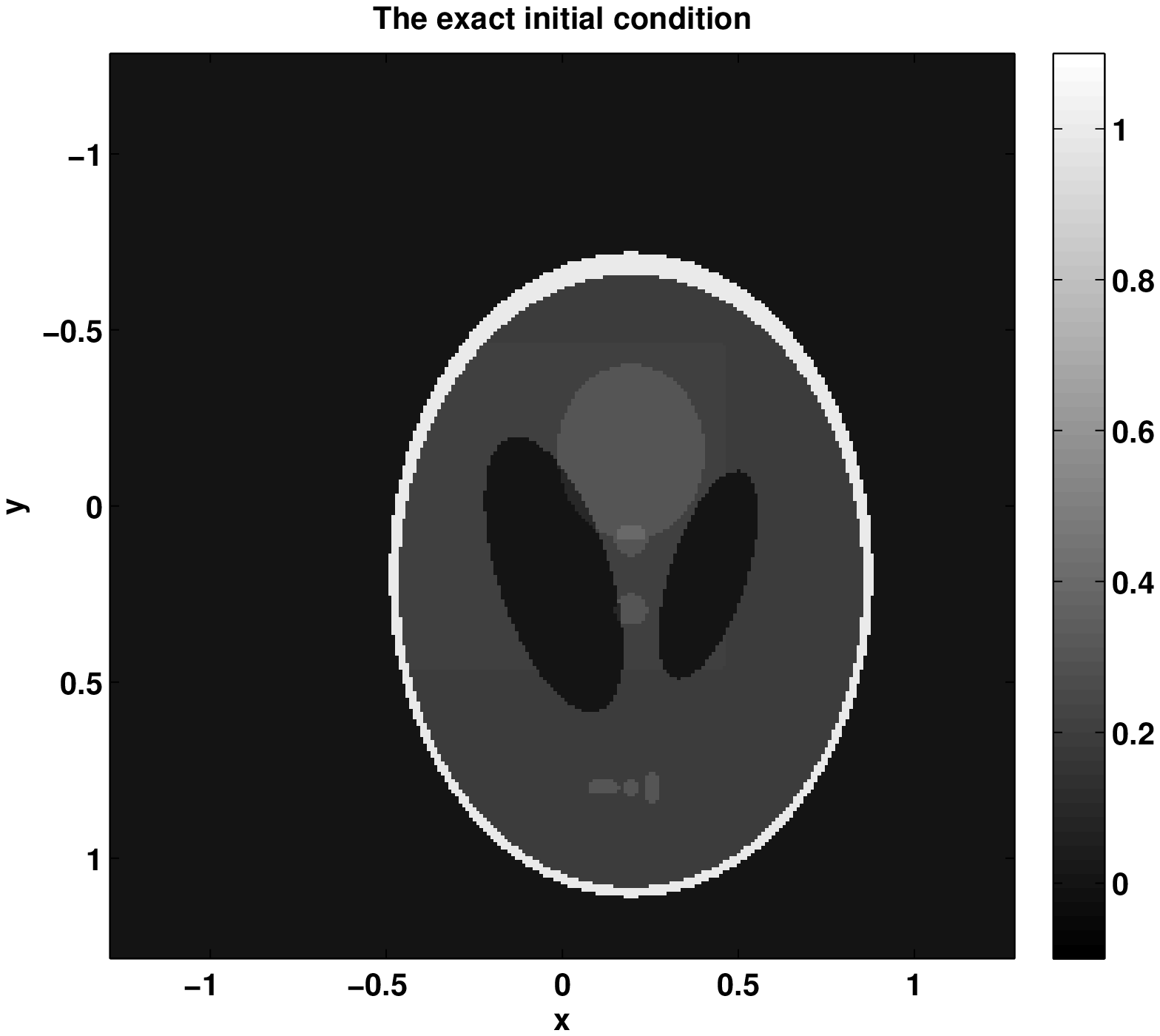,height=5.5cm }\\
(c)\epsfig{figure=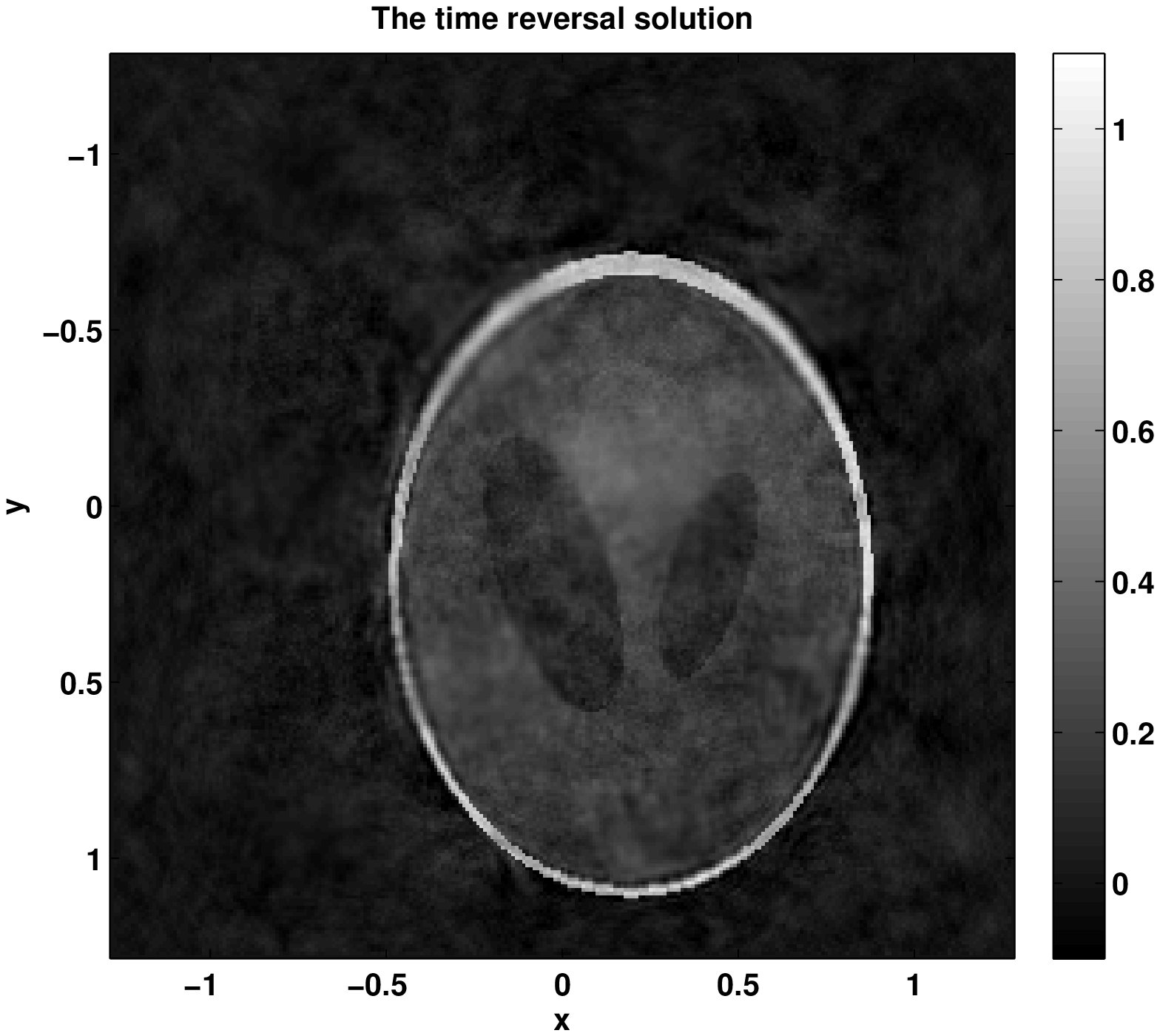,height=5.5cm }
(d)\epsfig{figure=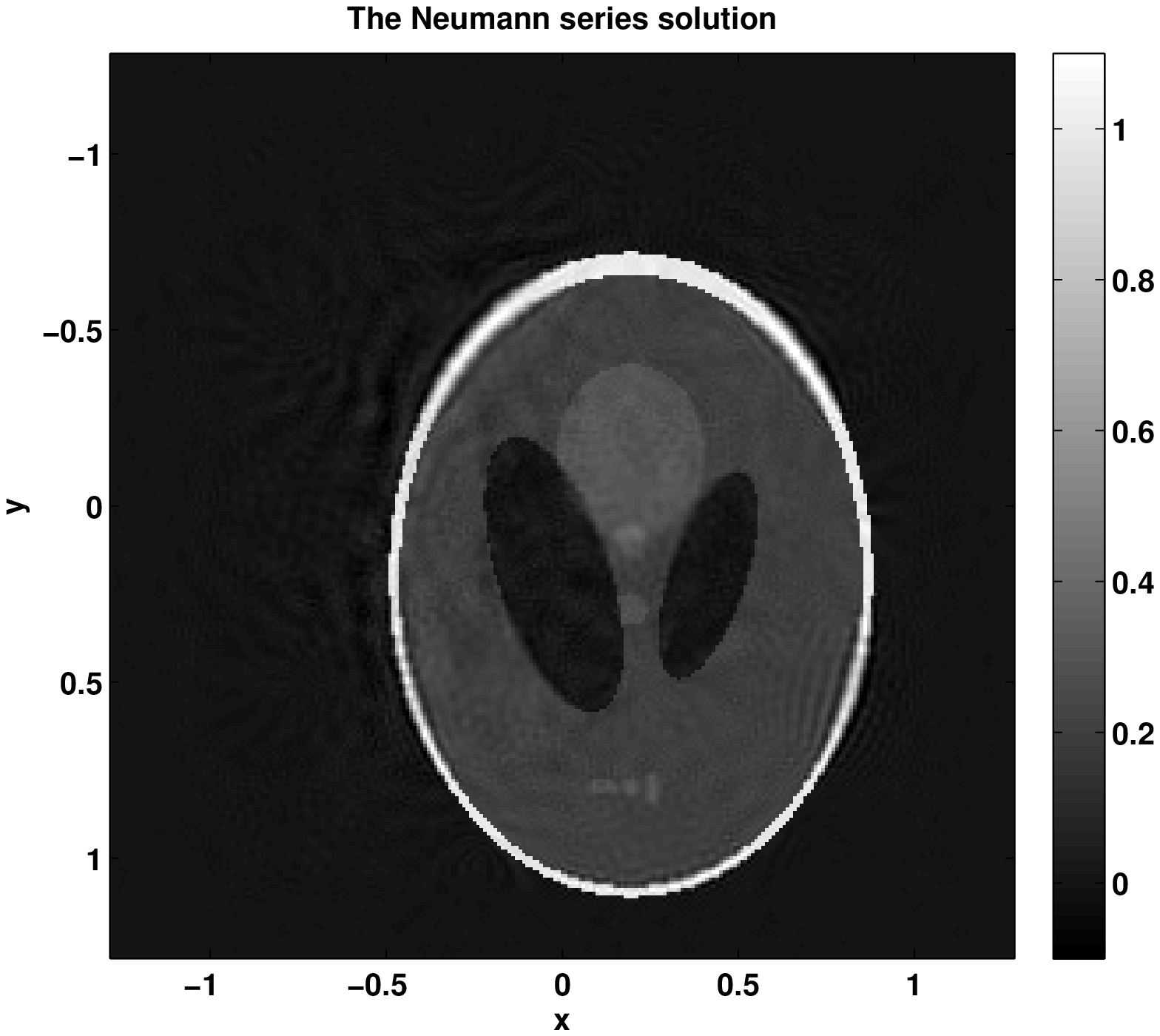,height=5.5cm }\\
(e)\epsfig{figure=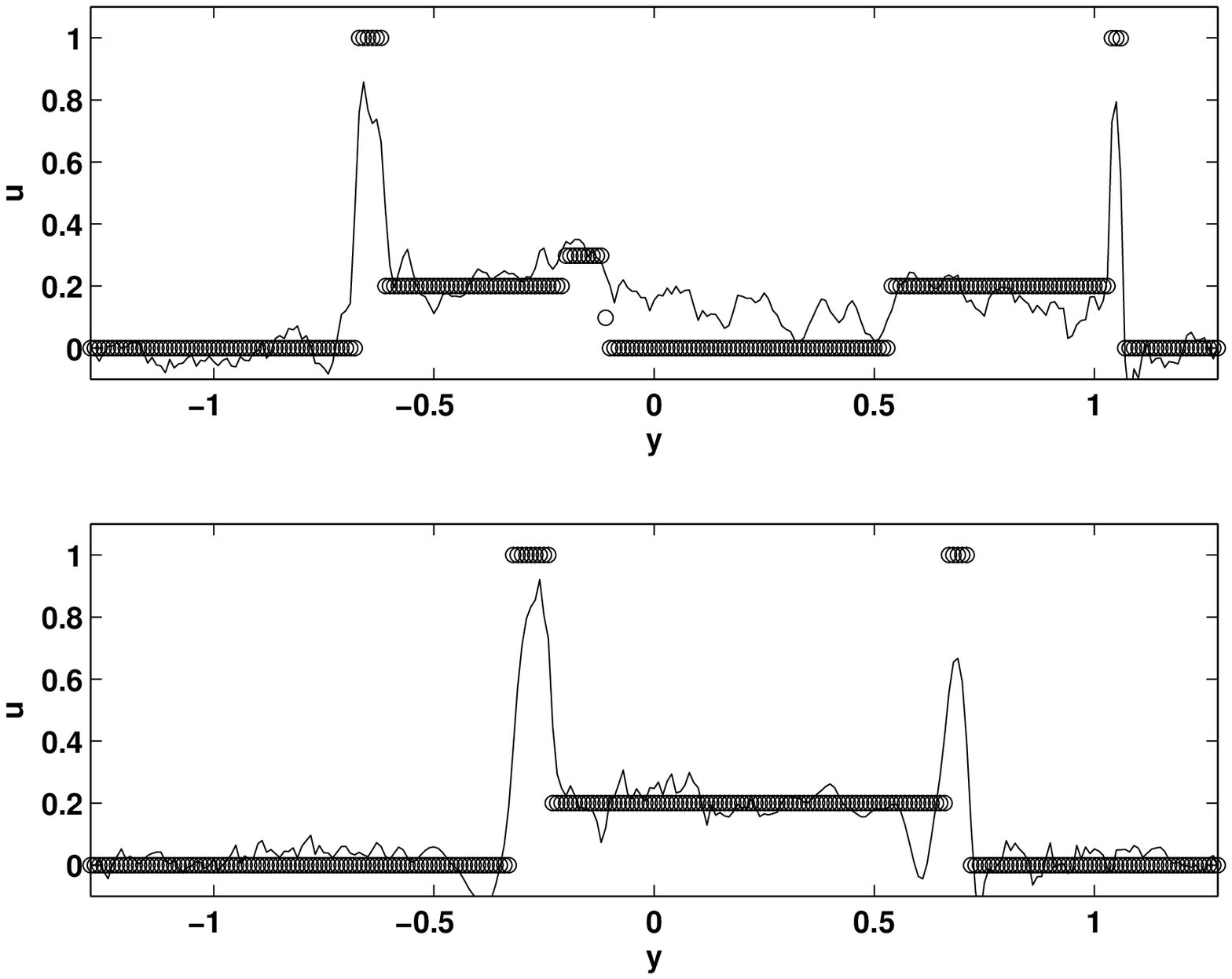,height=4cm,width=6.0cm}
(f)\epsfig{figure=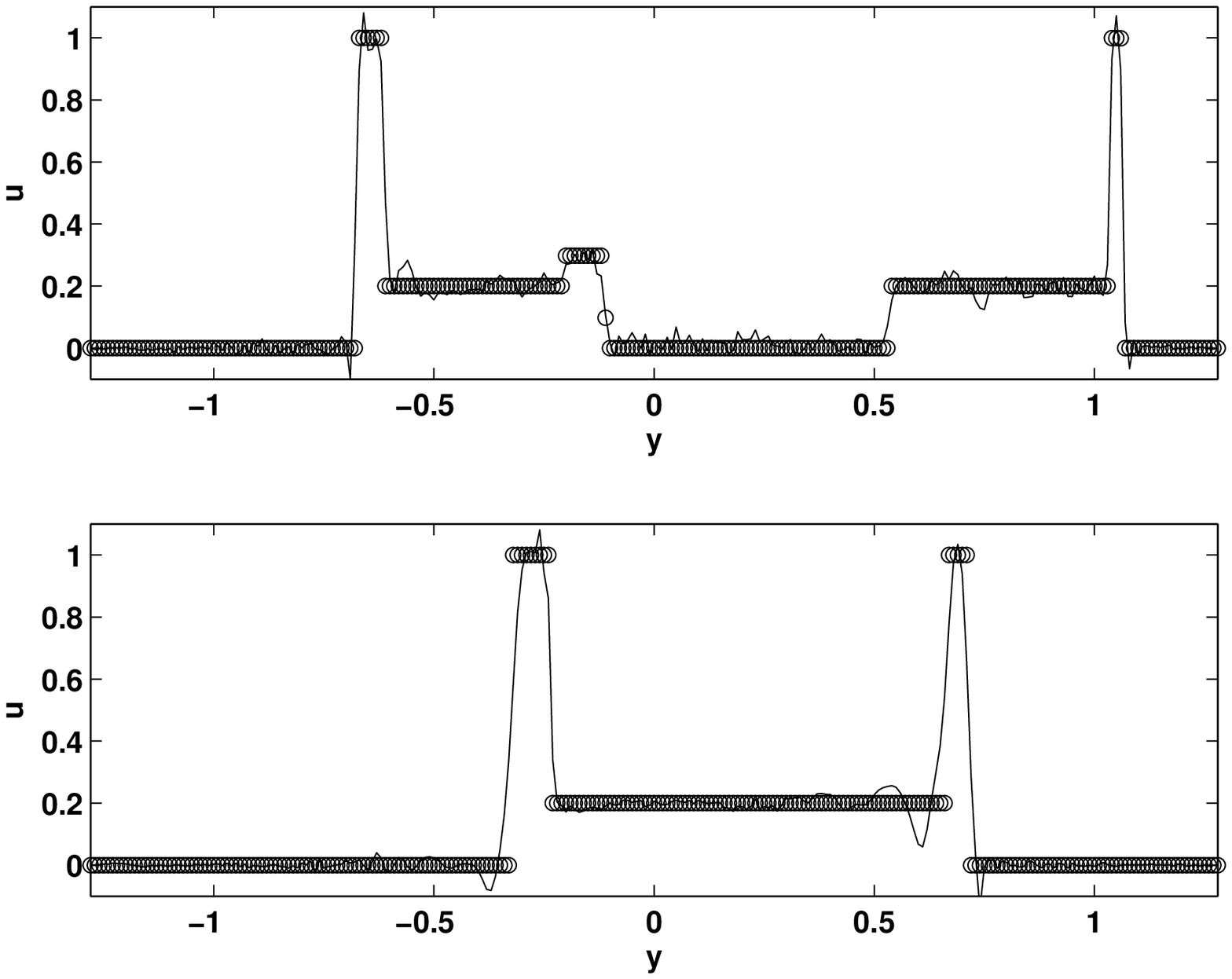,height=4cm,width=6.0cm}\\
(g)\epsfig{figure=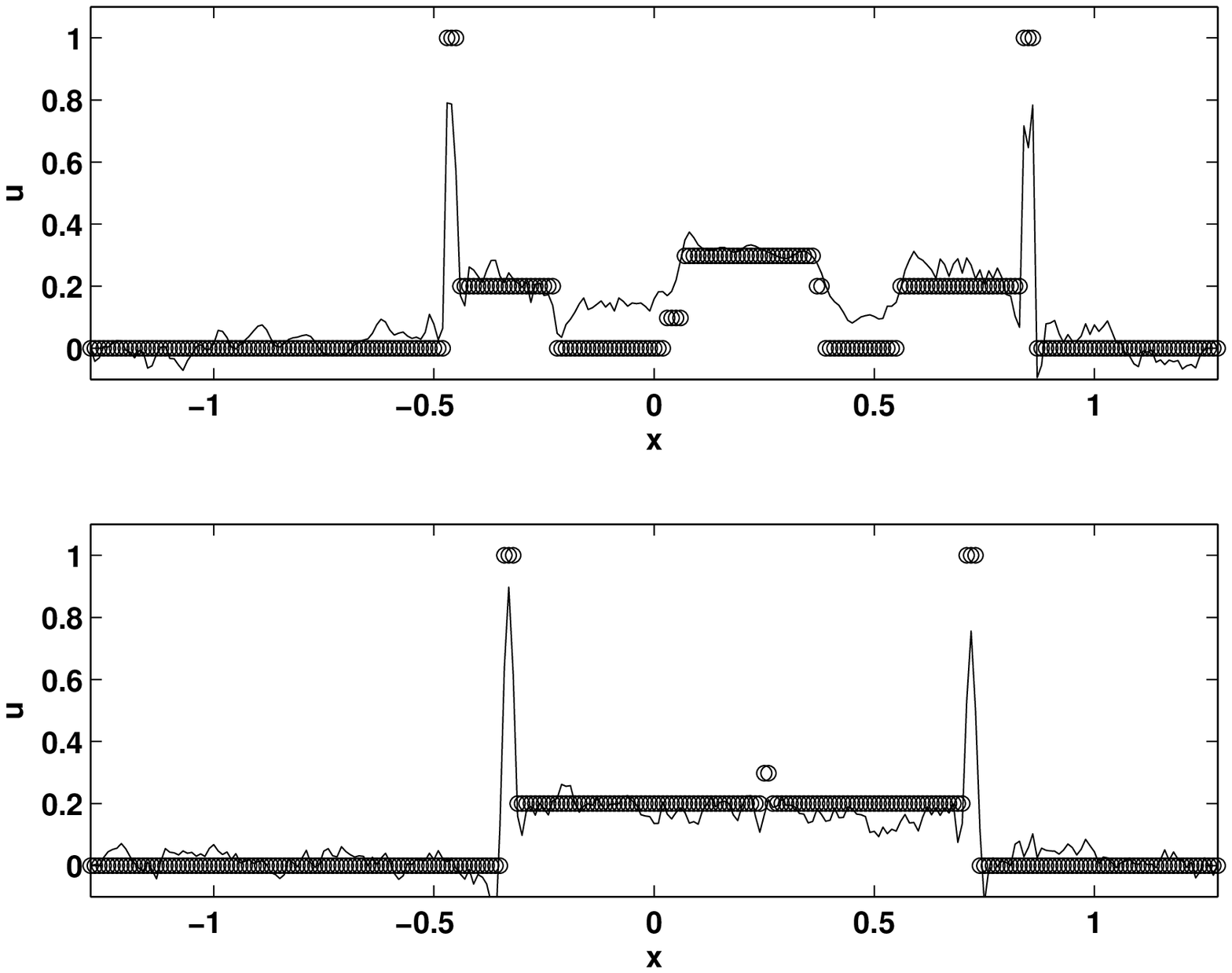,height=4cm,width=6.0cm}
(h)\epsfig{figure=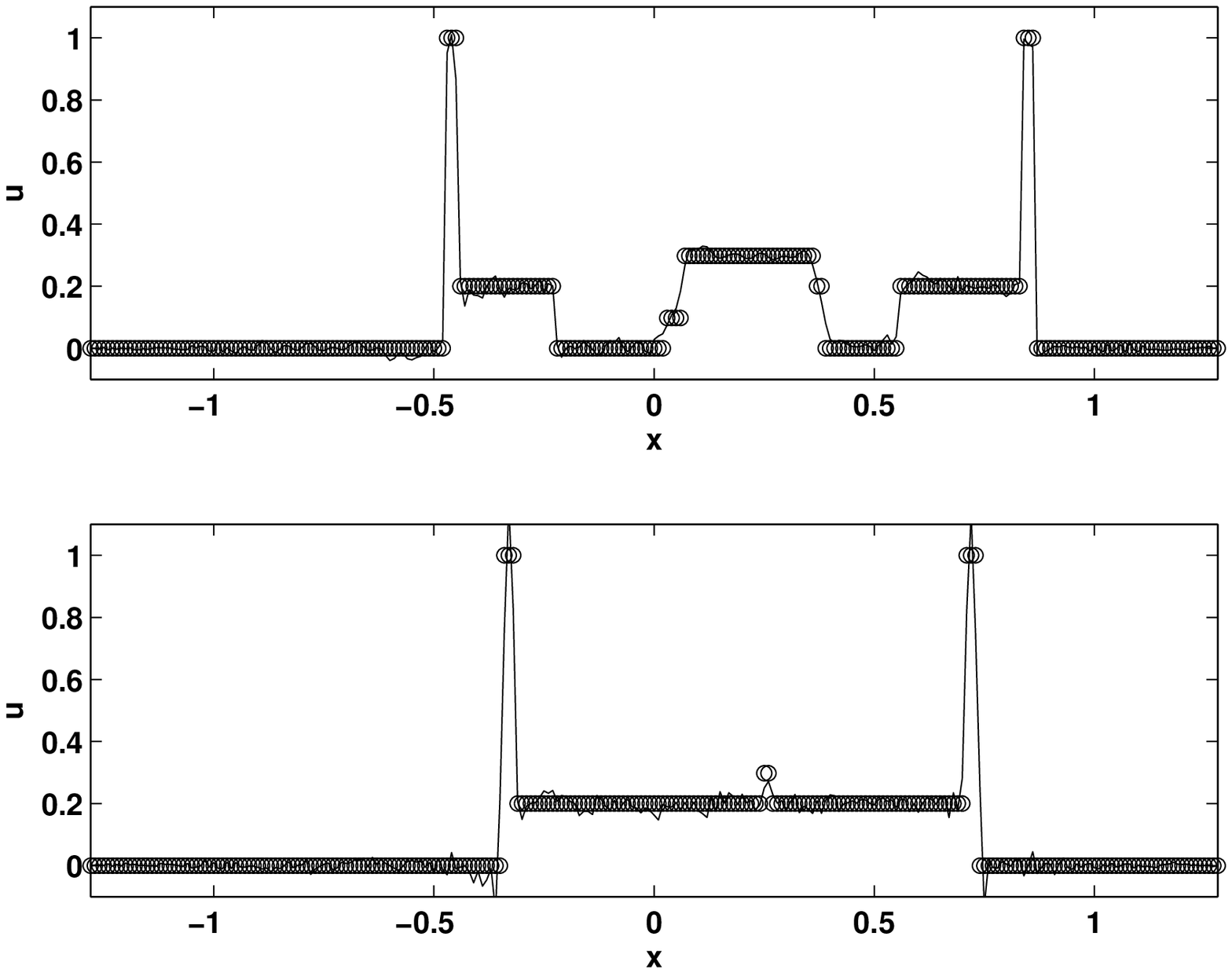,height=4cm,width=6.0cm}
\end{center}
\caption{Example 1 with the trapping speed $c_3$. $T=4T_0$. 
(a): the boundary distance map. 
(b): the exact initial condition. 
(c): the time reversal solution. 
(d): the Neumann series solution. 
(e): $x$-slices of the time reversal solution (``-'') and the exact solution (``o''). 
(f): $x$-slices of the Neumann series solution (``-'') and the exact solution (``o''). 
(g): $y$-slices of the time reversal solution (``-'') and the exact solution (``o'').
(h): $y$-slices of the Neumann series solution (``-'') and the exact solution (``o'').}
\label{Fig:2dsquare4T0Trap}
\end{figure}

\subsubsection{Trapping sound speed $c_3$, Figure \ref{Fig:2dsquare4T0Trap}
}
The sound speed is $c_3$ given by equation \eqref{key6} and  Figure~\ref{Fig:Velocities} (c). Although it is hard to show that the speed is actually trapping, one can easily construct numerically geodesics of length at least $20$; therefore, based on this numerical evidence, $T_1$ is at least $20$ and possibly $\infty$.  The times $T$ that we choose are much smaller than $20/2$ so there is a considerable set of singularities that are not detected. 

We have $T_0\approx 1.25$. In  Figure \ref{Fig:2dsquare4T0Trap}, we present reconstructions with  $T= 4T_0$. The phantom in the middle of the NS image (error $18.11\%$, $k=8$)  can be separated from the background, while in the TR one (error $39.56\%$) --- much harder. The invisible singularities in this case are distributed in a chaotic way and what appears as noise in both images is due to the fact that there are singularities there that cannot be resolved --- even if the actual image is not singular there!

\subsection{Example 2: Zebras, Figures~\ref{Fig:2dZebra4T0Notrap}--\ref{Fig:2dZebra4T0Speed63}}
\subsubsection{Non-trapping speed $c_1$}

The sound speed is given by equation \eqref{key8} and it is visualized in Figure~\ref{Fig:Velocities}(c). Here $f$ now is represented by the zebras image that has more complex structure. 

\begin{figure}
\begin{center}
(a)\epsfig{figure=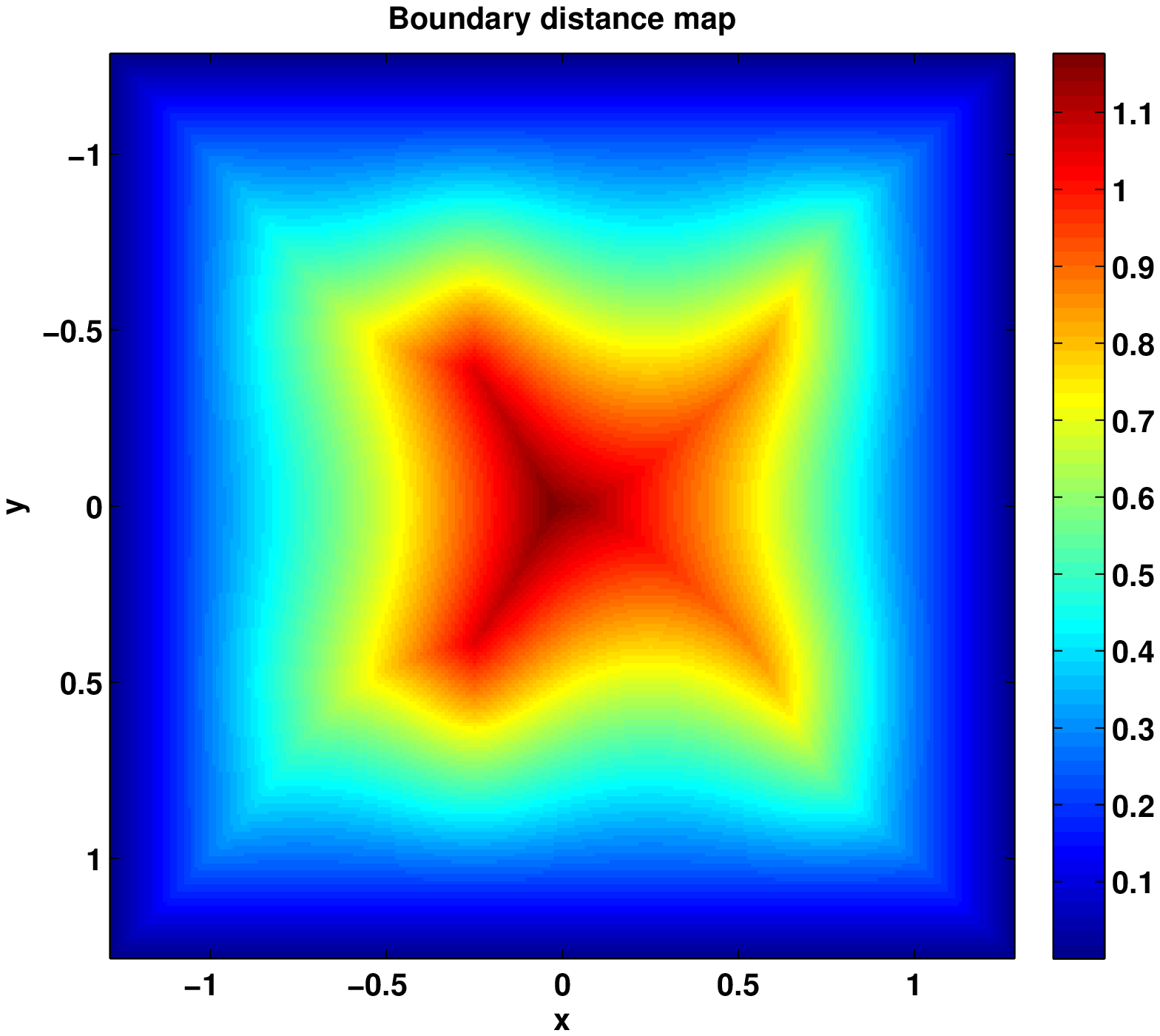,height=5.5cm }
(b)\epsfig{figure=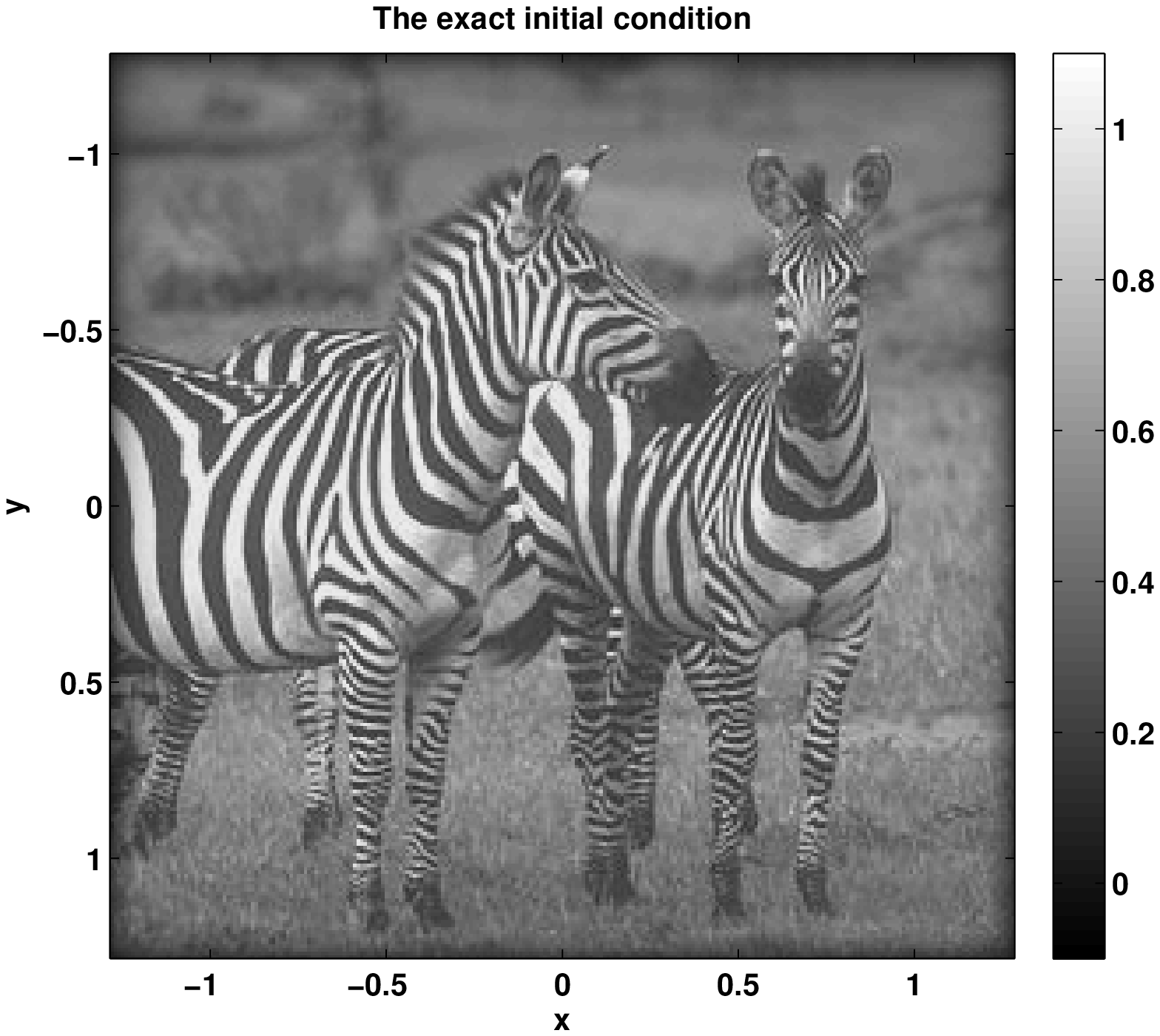,height=5.5cm }\\
(c)\epsfig{figure=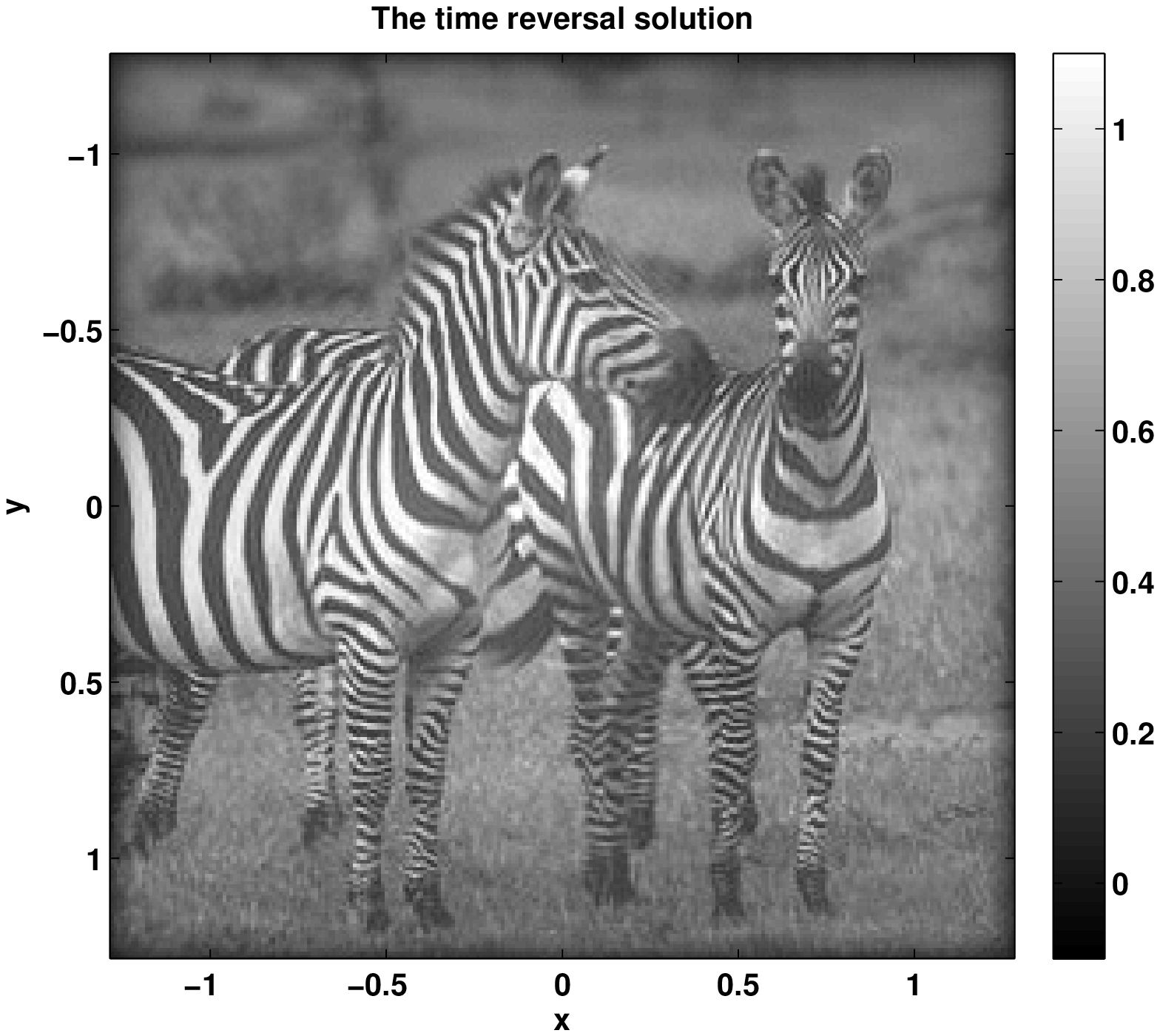,height=5.5cm }
(d)\epsfig{figure=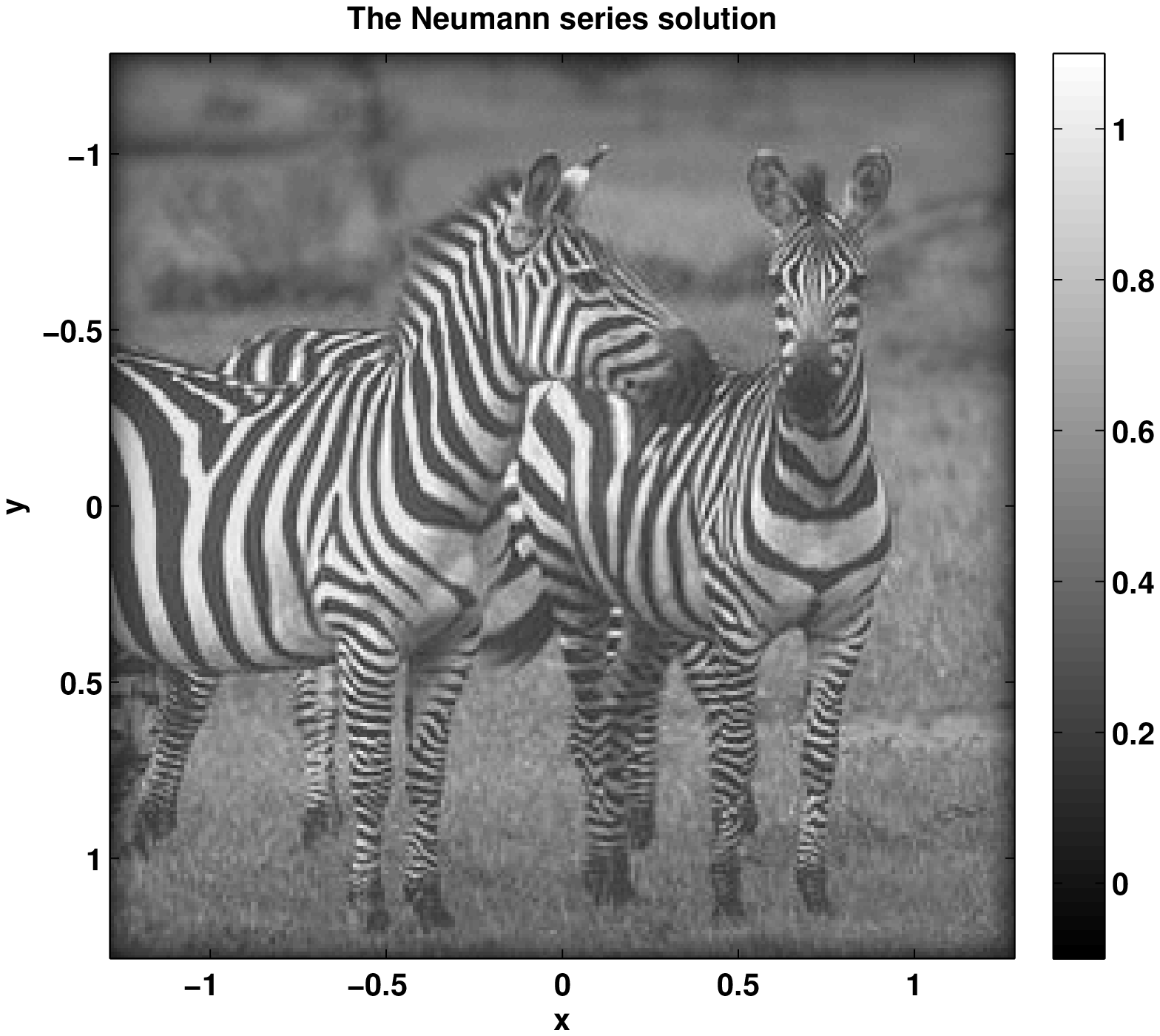,height=5.5cm }\\
(e)\epsfig{figure=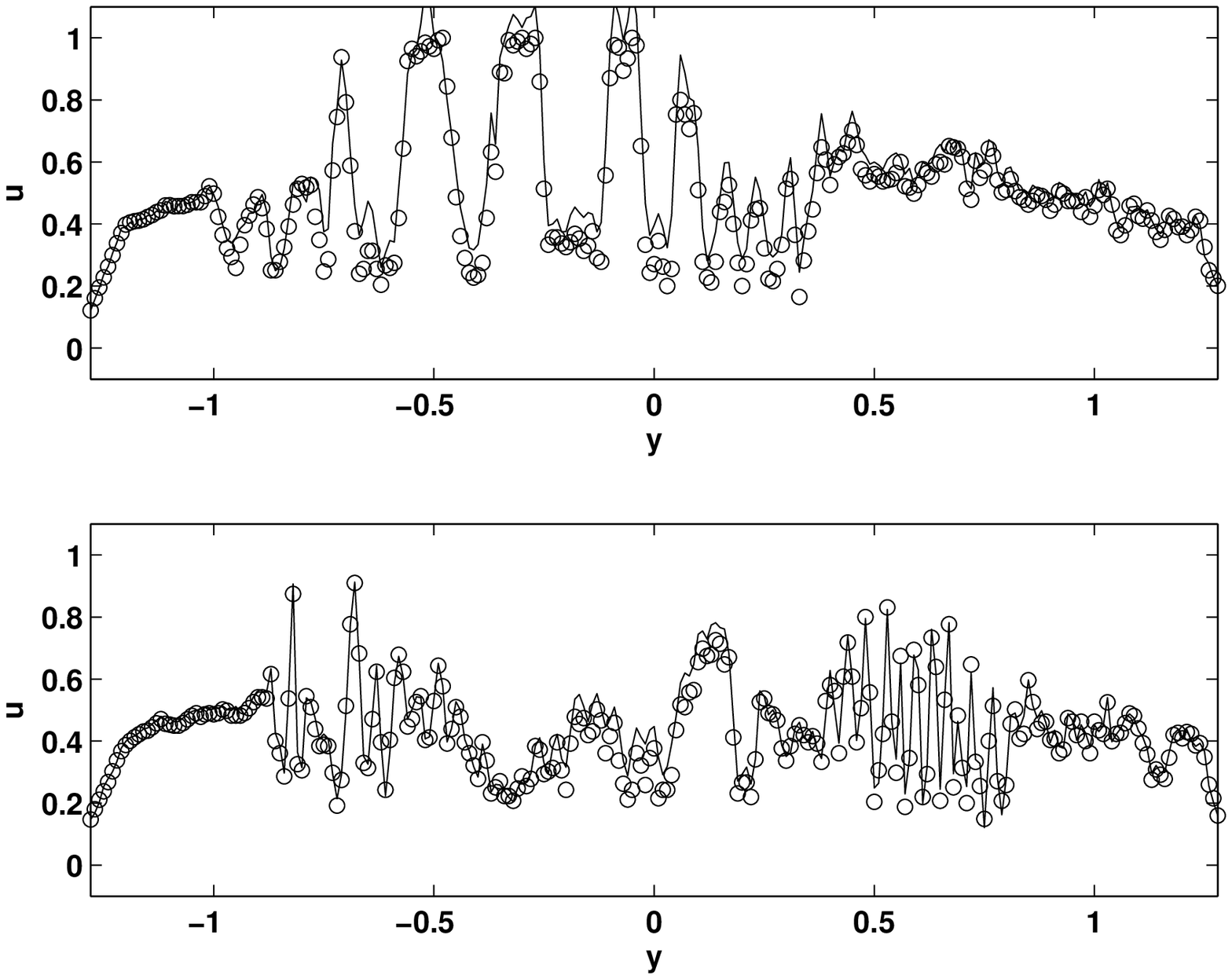,height=4cm,width=6.0cm}
(f)\epsfig{figure=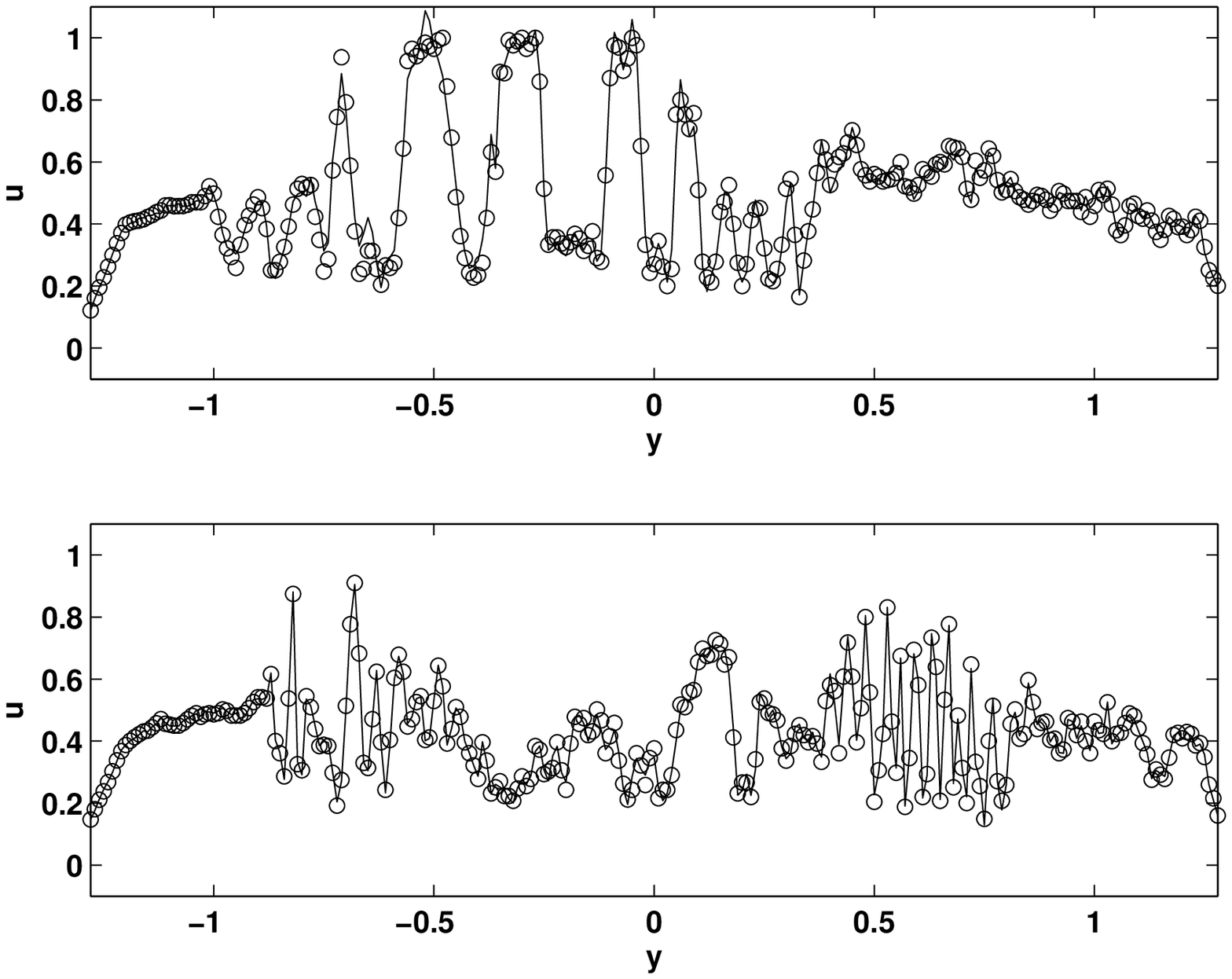,height=4cm,width=6.0cm}\\
(g)\epsfig{figure=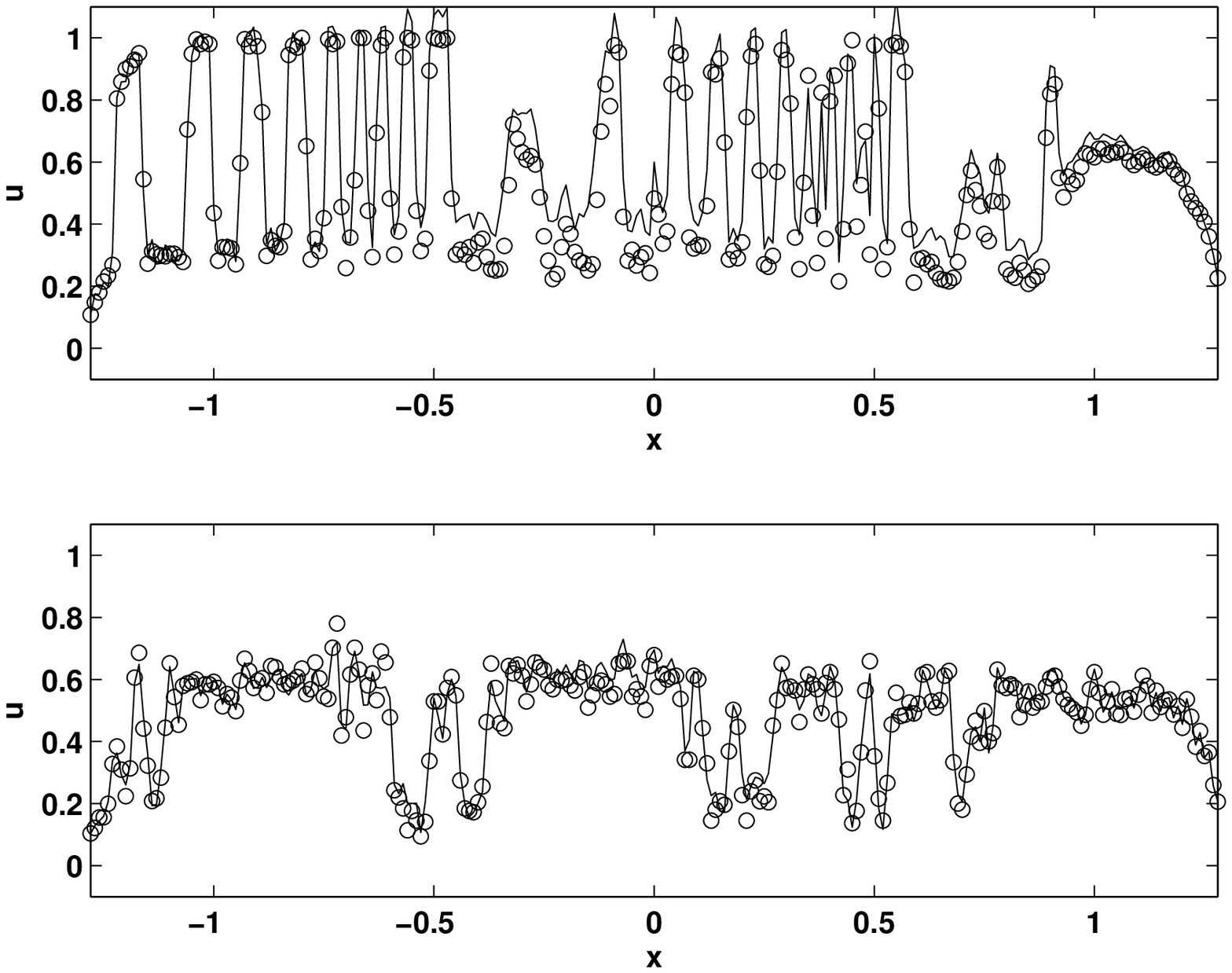,height=4cm,width=6.0cm}
(h)\epsfig{figure=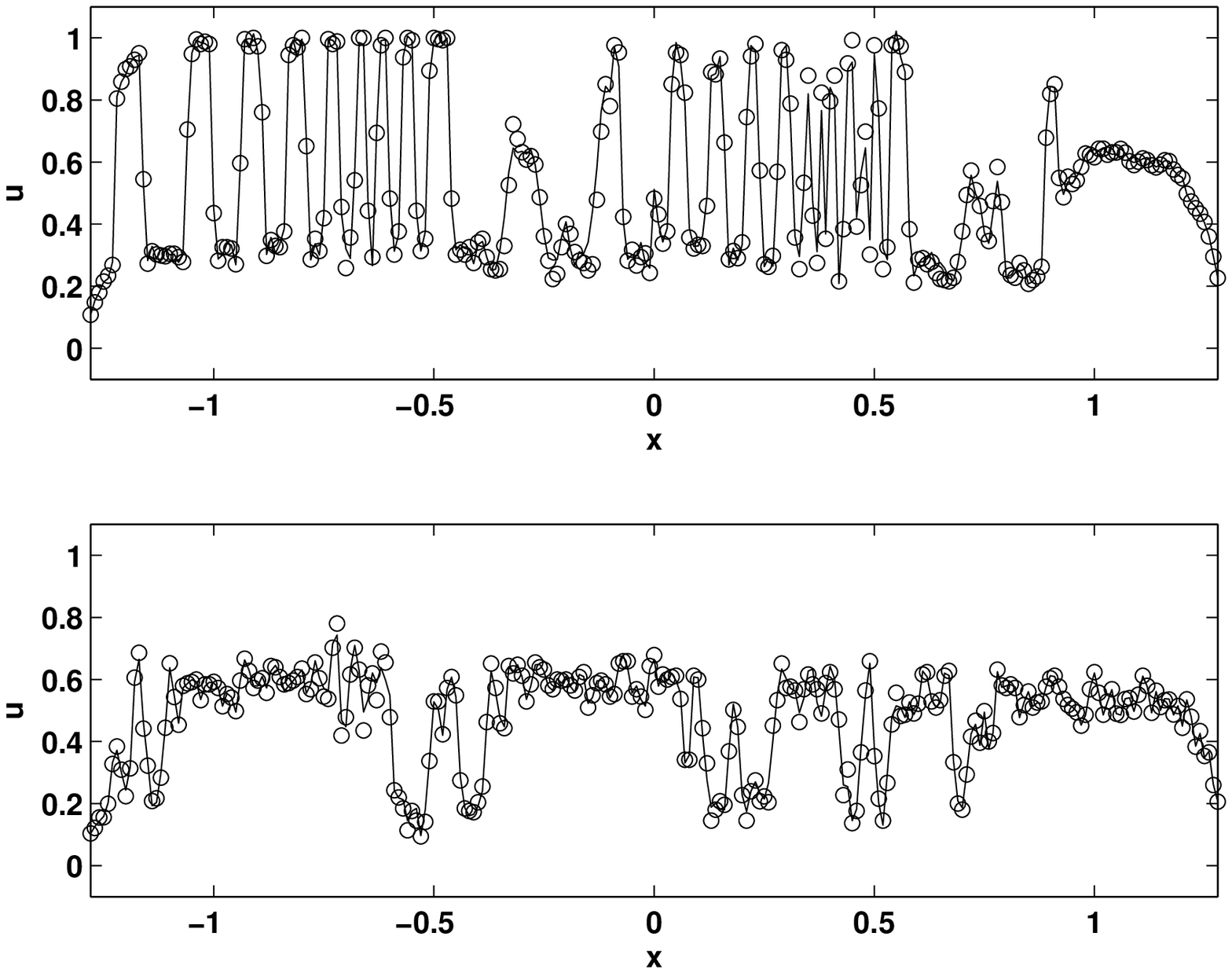,height=4cm,width=6.0cm}
\end{center}
\caption{Example 2 with the non-trapping speed $c_1$. Case 1: $T=4T_0$. 
(a): the boundary distance map. 
(b): the exact initial condition. 
(c): the time reversal solution. 
(d): the Neumann series solution. 
(e): $x$-slices of the time reversal solution (``-'') and the exact solution (``o''). 
(f): $x$-slices of the Neumann series solution (``-'') and the exact solution (``o''). 
(g): $y$-slices of the time reversal solution (``-'') and the exact solution (``o'').
(h): $y$-slices of the Neumann series solution (``-'') and the exact solution (``o'').}
\label{Fig:2dZebra4T0Notrap}
\end{figure}

\textbf{Figure~\ref{Fig:2dZebra4T0Notrap}:} $T=4T_0$.
 Both methods give a good reconstruction, and with $k=2$ only, the NS one has error $4.6\%$ that is about $1/2$ the  of the TR one. The error can be improved significantly with more terms. 

\begin{figure}
\begin{center}
(a)\epsfig{figure=Fig/Key8Keyini6341Traveltime.ps,height=5.5cm }
(b)\epsfig{figure=Fig/Key8Keyini6341NX301InitialTarget.ps,height=5.5cm }\\
(c)\epsfig{figure=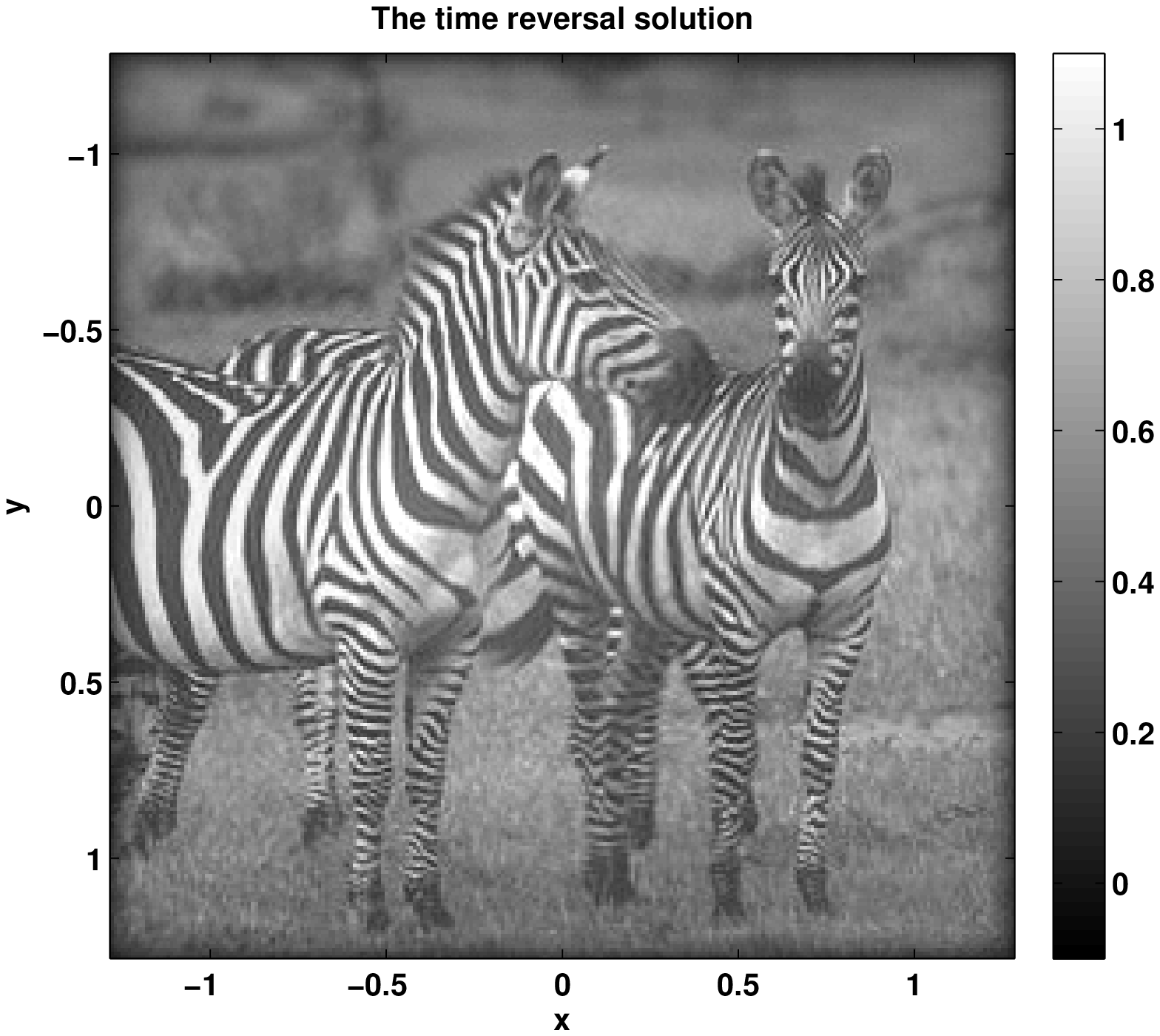,height=5.5cm }
(d)\epsfig{figure=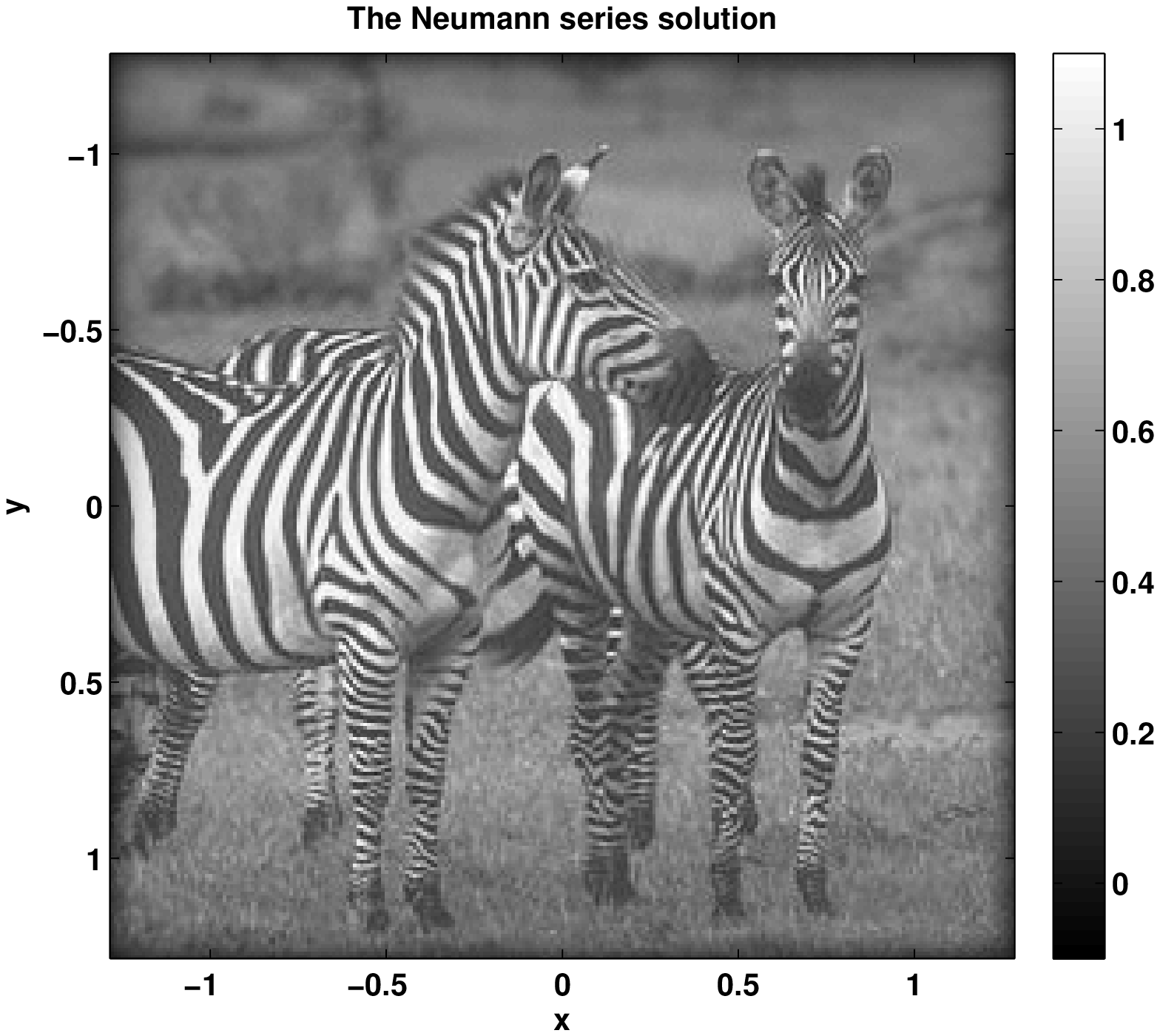,height=5.5cm }\\
(e)\epsfig{figure=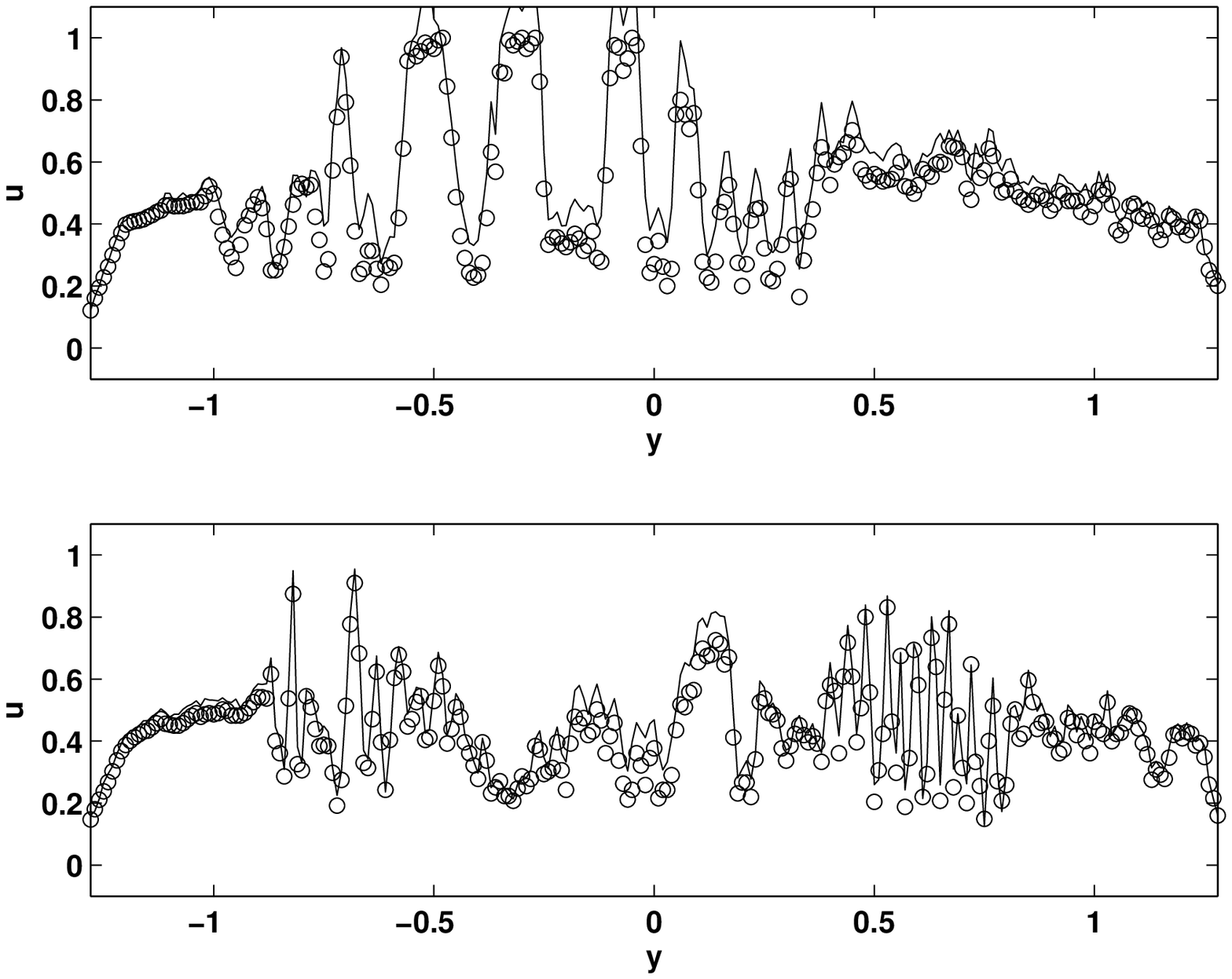,height=4cm,width=6.0cm}
(f)\epsfig{figure=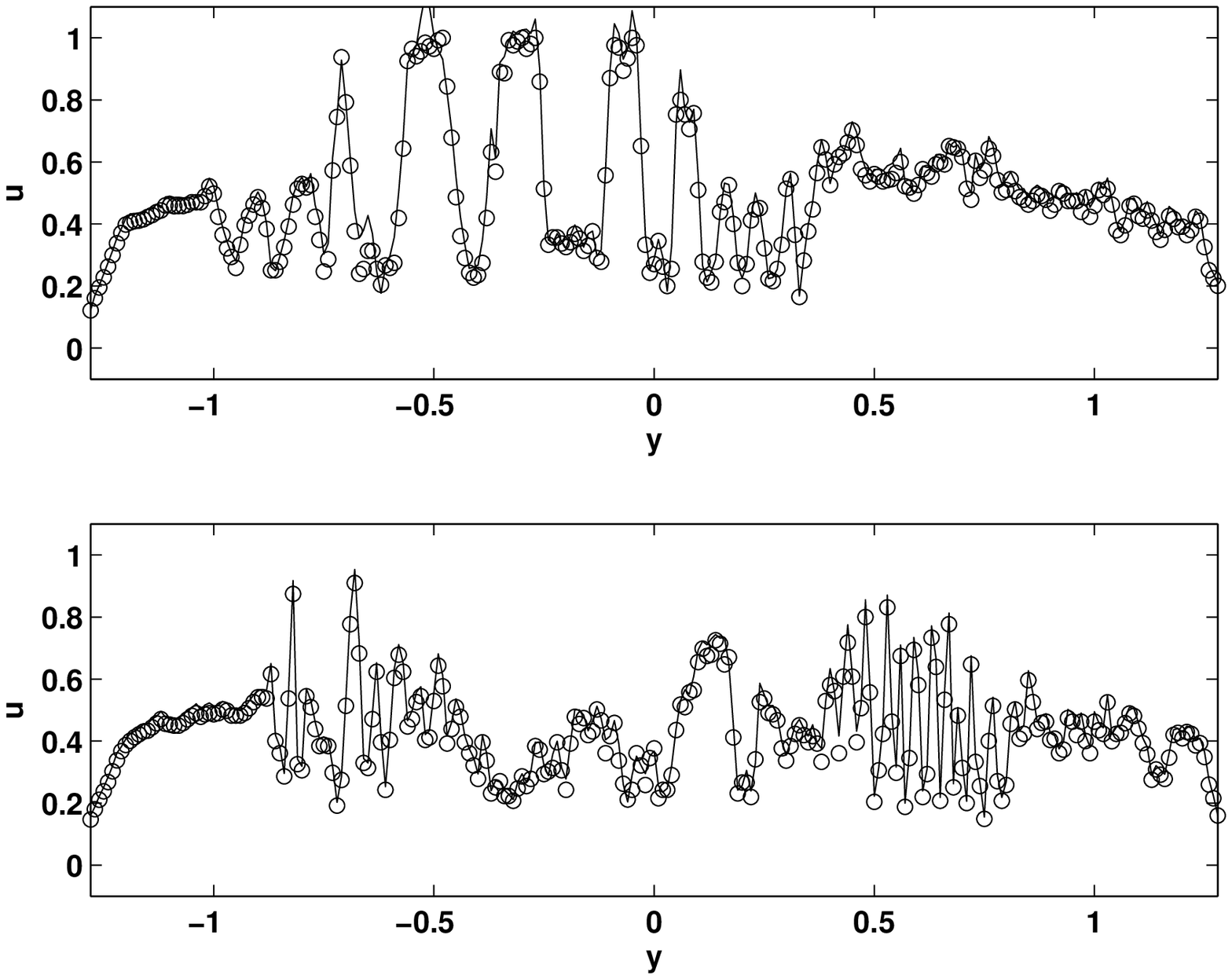,height=4cm,width=6.0cm}\\
(g)\epsfig{figure=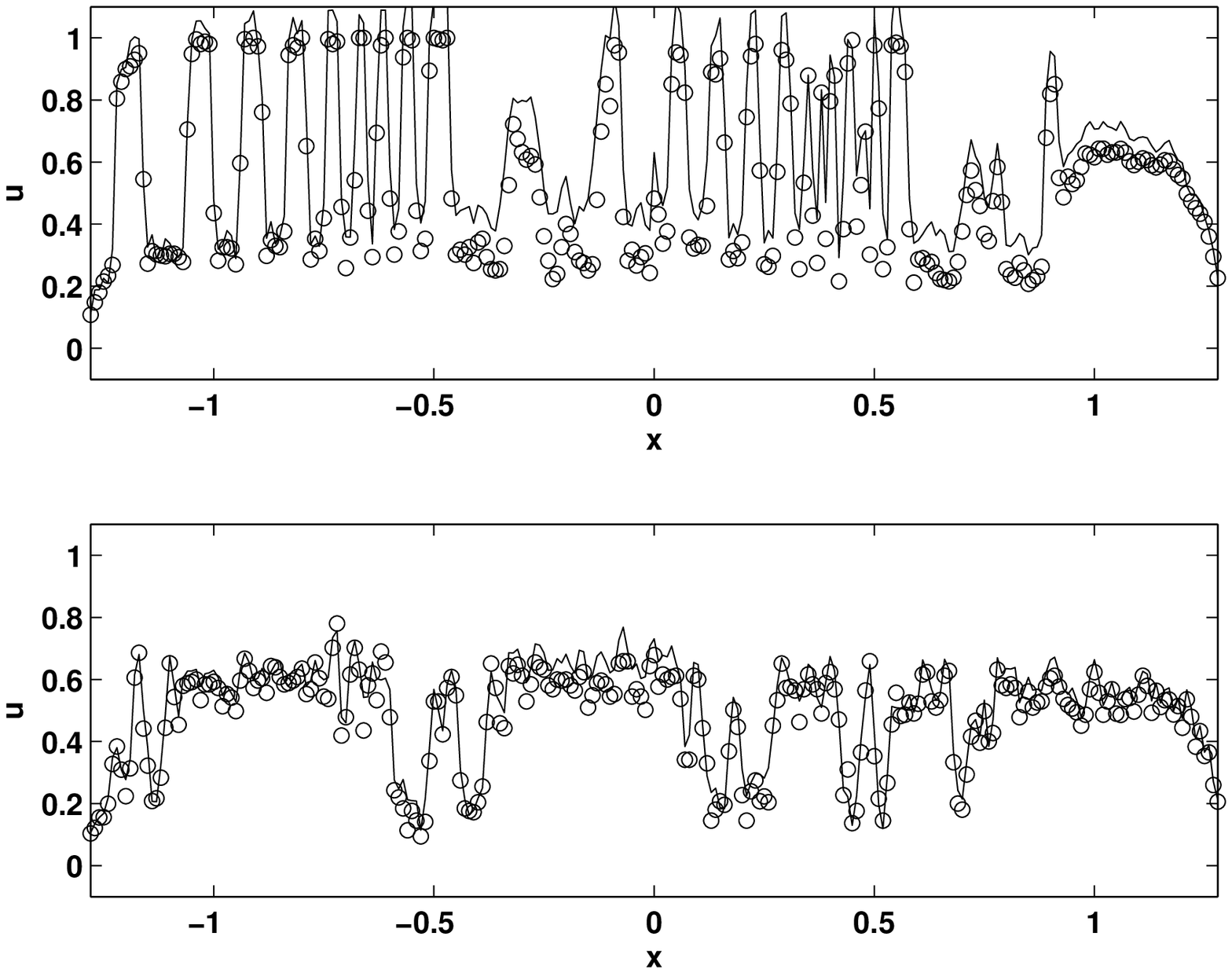,height=4cm,width=6.0cm}
(h)\epsfig{figure=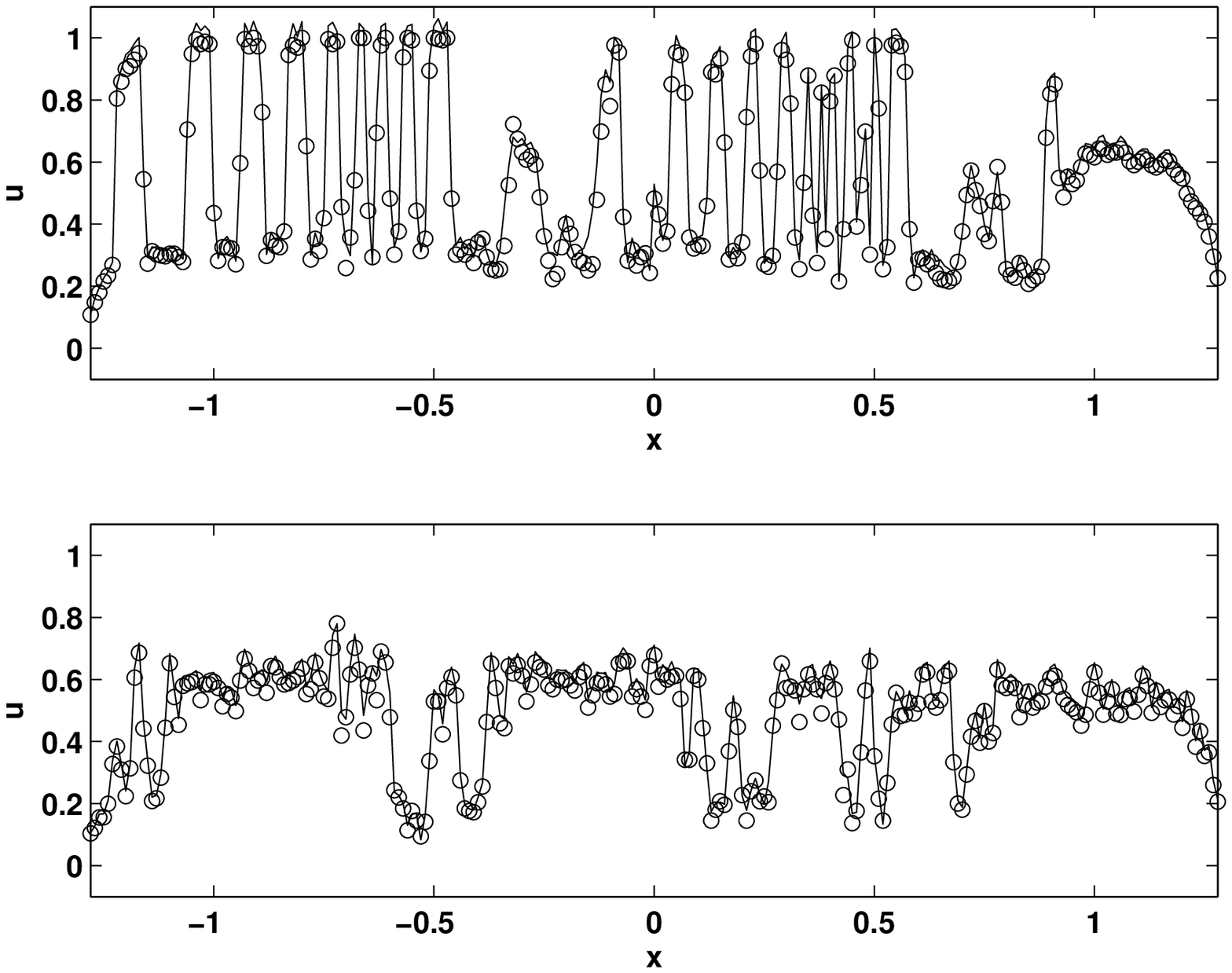,height=4cm,width=6.0cm}
\end{center}
\caption{Example 2 with the non-trapping speed $c_1$. Case 2: $T=4T_0$ with $10\%$ noise. 
(a): the boundary distance map. 
(b): the exact initial condition. 
(c): the time reversal solution. 
(d): the Neumann series solution. 
(e): $x$-slices of the time reversal solution (``-'') and the exact solution (``o''). 
(f): $x$-slices of the Neumann series solution (``-'') and the exact solution (``o''). 
(g): $y$-slices of the time reversal solution (``-'') and the exact solution (``o'').
(h): $y$-slices of the Neumann series solution (``-'') and the exact solution (``o'').}
\label{Fig:2dZebra4T0NotrapNoise}
\end{figure}

\textbf{Figure~\ref{Fig:2dZebra4T0NotrapNoise}:}  
$T=4T_0$ with $10\%$ noise. Here, $k=8$, with an error $6.25\%$, about  $1/2$ the  of the TR one.

\subsubsection{Trapping sound speed $c_3$}
The sound speed is given by equation \eqref{key6}. As we indicated above, $T_0\approx  1.23$, $T_1>30$, and the geodesic 
flow is somewhat chaotic with invisible singularities that have also chaotic distribution.

\begin{figure} 
\begin{center}
(a)\epsfig{figure=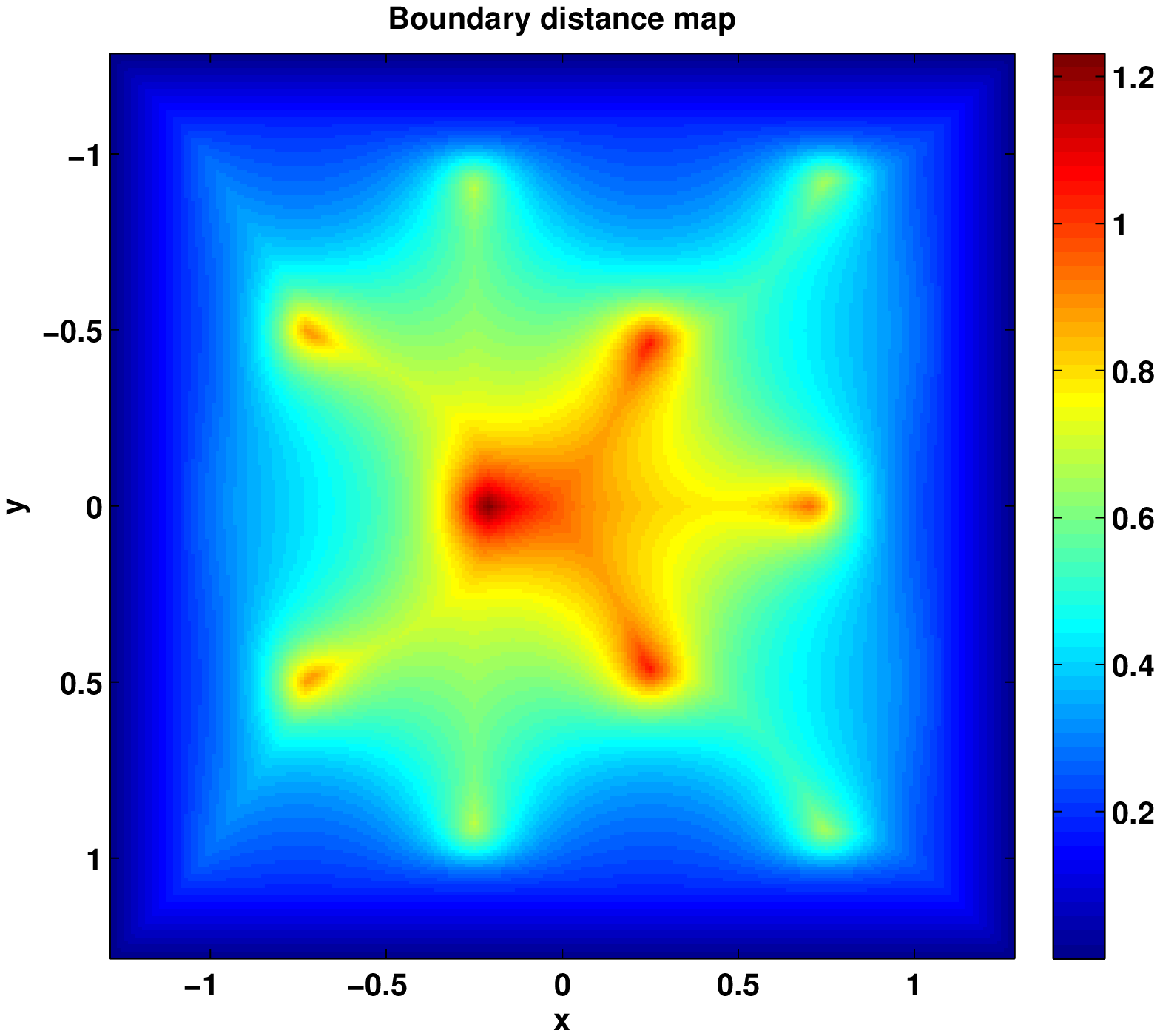,height=5.5cm }
(b)\epsfig{figure=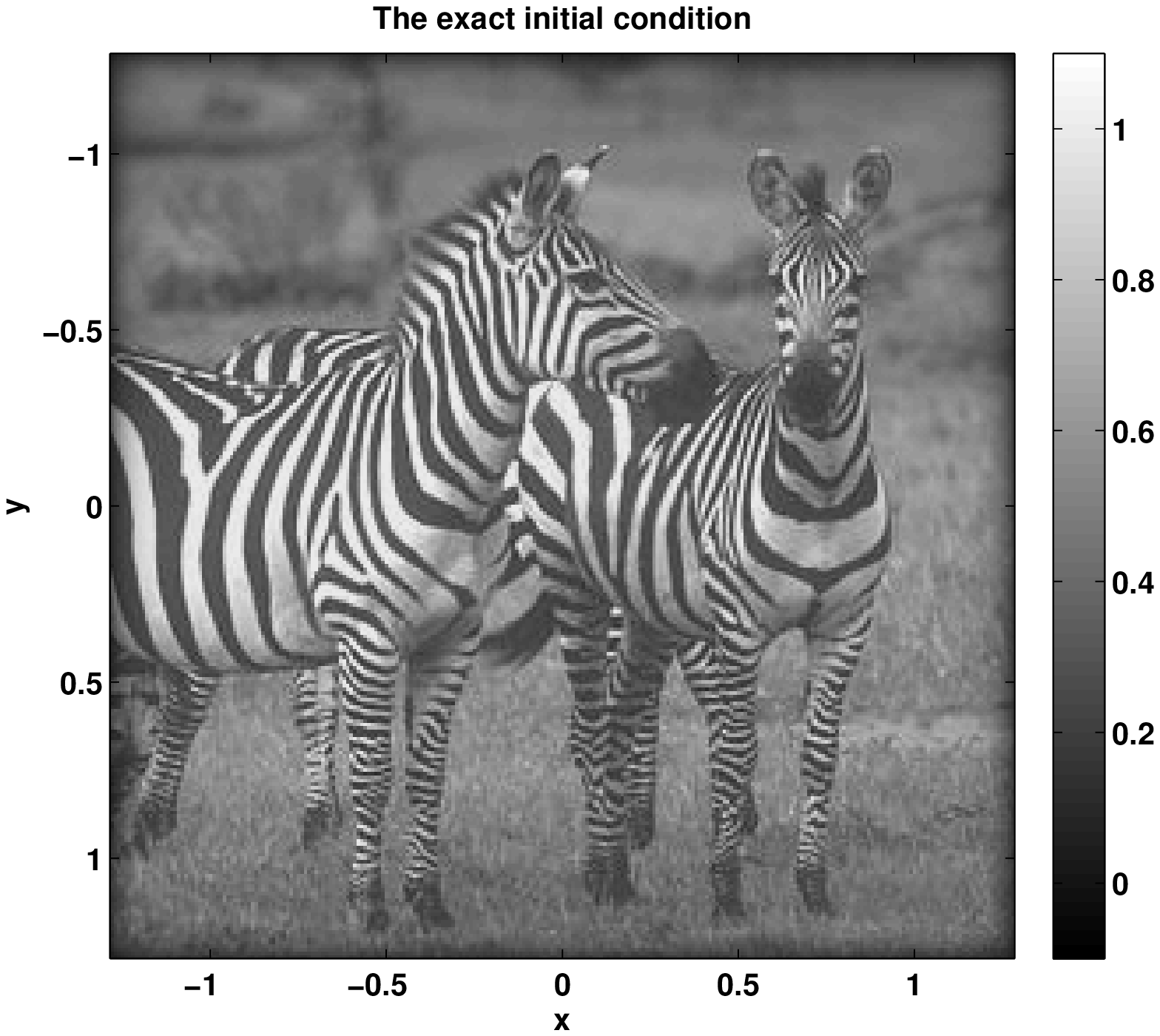,height=5.5cm }\\
(c)\epsfig{figure=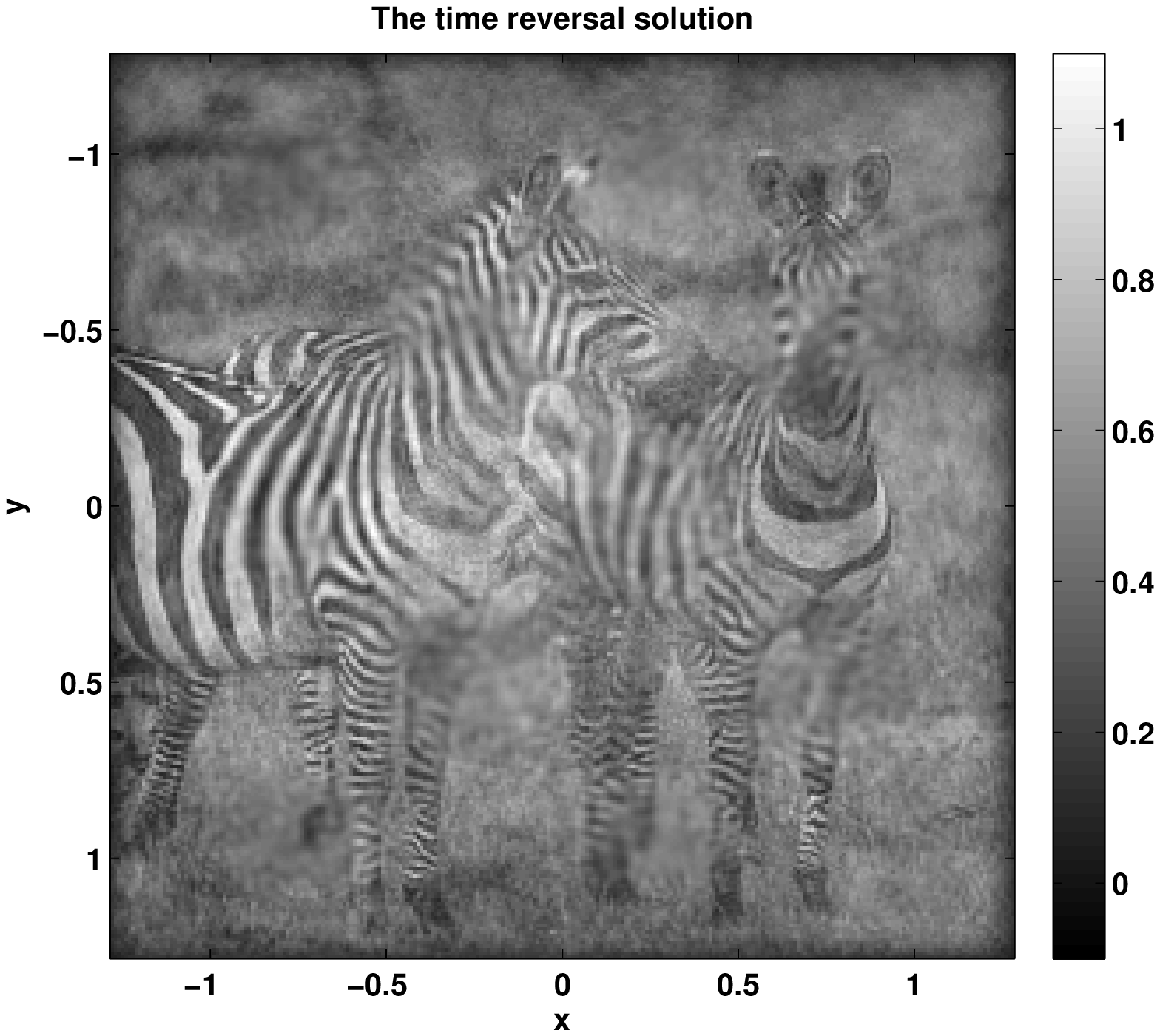,height=5.5cm }
(d)\epsfig{figure=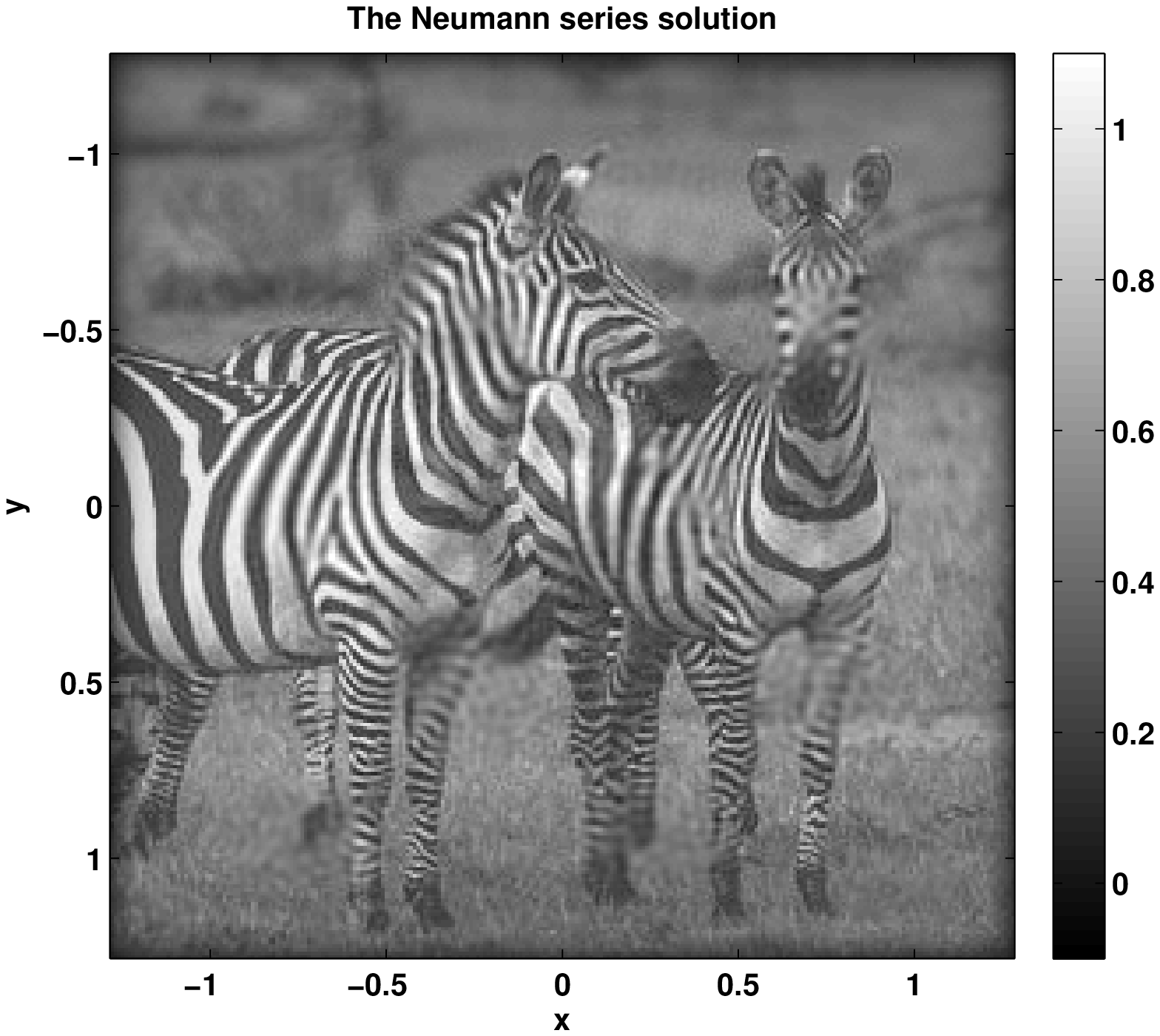,height=5.5cm }
\end{center}
\caption{Example 2 with the trapping speed $c_3$. Case 1: $T=4T_0$. 
(a): the boundary distance map. 
(b): the exact initial condition. 
(c): the time reversal solution. 
(d): the Neumann series solution. 
}
\label{Fig:2dZebra4T0trap}
\end{figure}

\textbf{Figure \ref{Fig:2dZebra4T0trap}:}  $T= 4T_0$. The NS reconstruction is much cleaner, error $9.64\%$ with $k=20$ that is a bit less than $1/2$ of the TR error. 
 
\begin{figure} 
\begin{center}
(a)\epsfig{figure=Fig/Key6Keyini6341Traveltime.ps,height=5.5cm }
(b)\epsfig{figure=Fig/Key6Keyini6341NX301InitialTarget.ps,height=5.5cm }\\
(c)\epsfig{figure=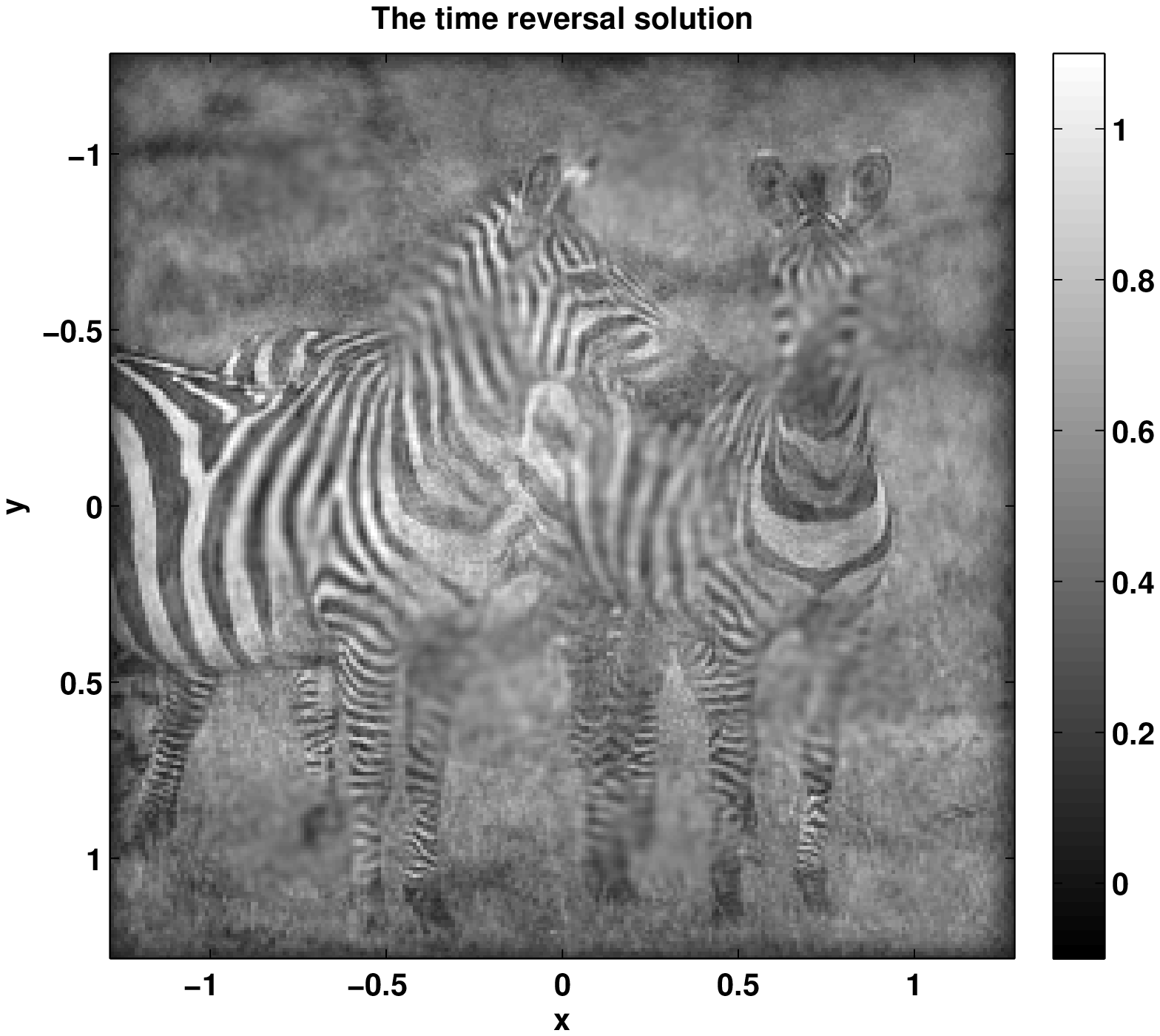,height=5.5cm }
(d)\epsfig{figure=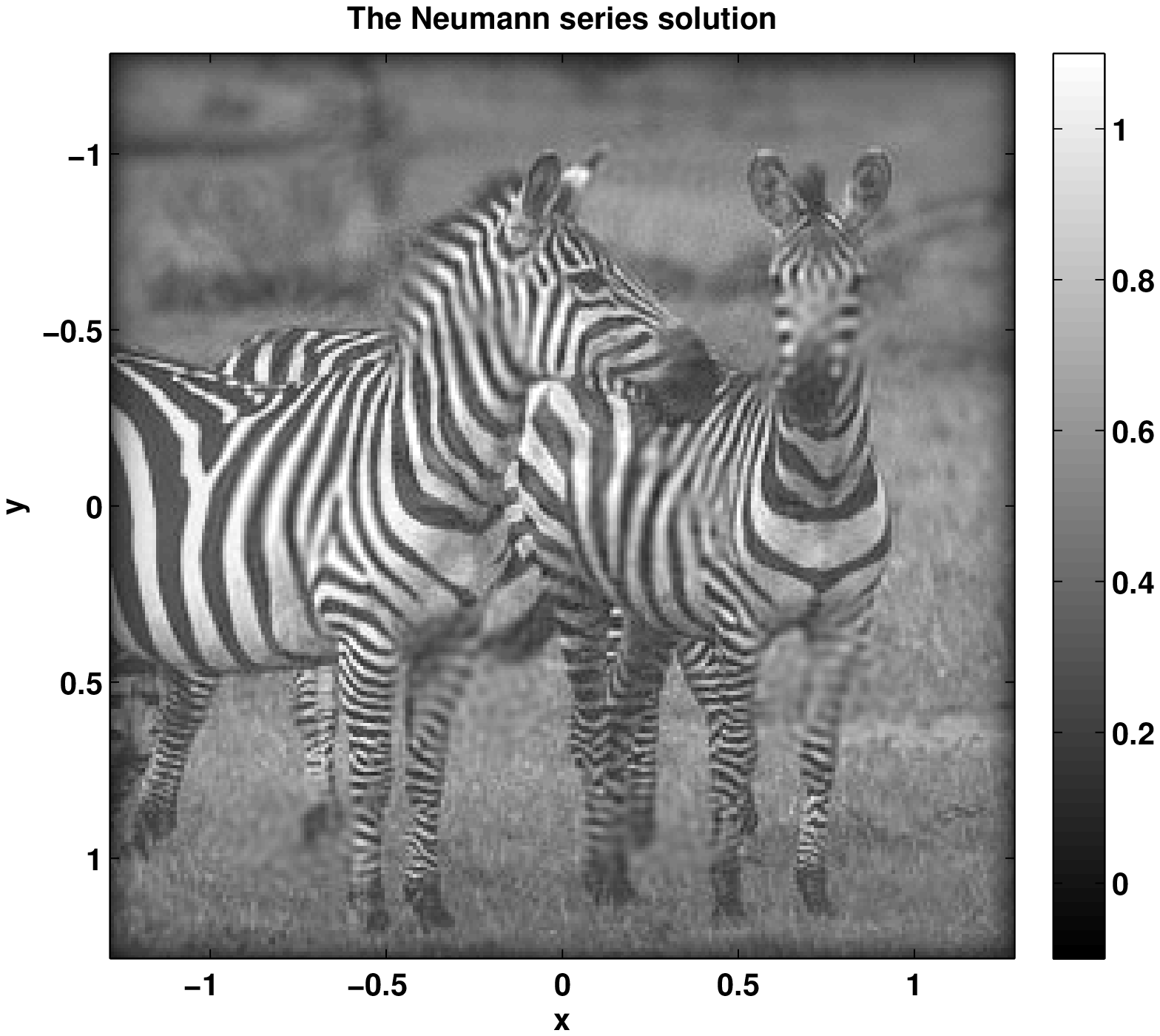,height=5.5cm }\\
\end{center}
\caption{Example 2 with the trapping speed $c_3$. Case 2: $T=4T_0$ with $10\%$ random noise.
(a): the boundary distance map. 
(b): the exact initial condition. 
(c): the time reversal solution. 
(d): the Neumann series solution. 
}
\label{Fig:2dZebra4T0trapNoise}
\end{figure}

\textbf{Figure \ref{Fig:2dZebra4T0trapNoise}:} $T= 4T_0$ with noise. The noise increases slightly the error in both cases.

\subsubsection{Trapping speed $c_2$}

\begin{figure}%
\begin{center}
(a)\epsfig{figure=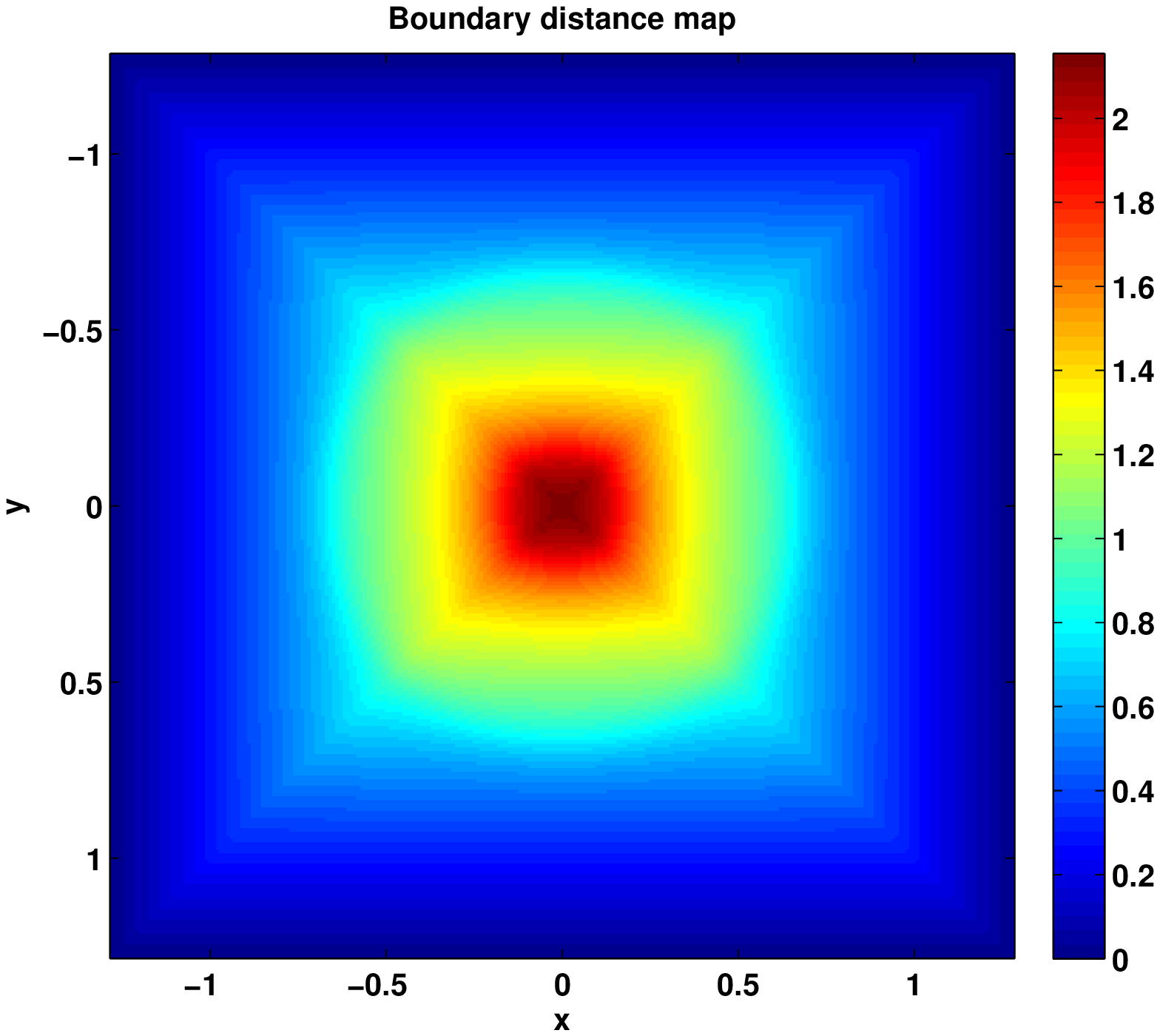,height=5.5cm }
(b)\epsfig{figure=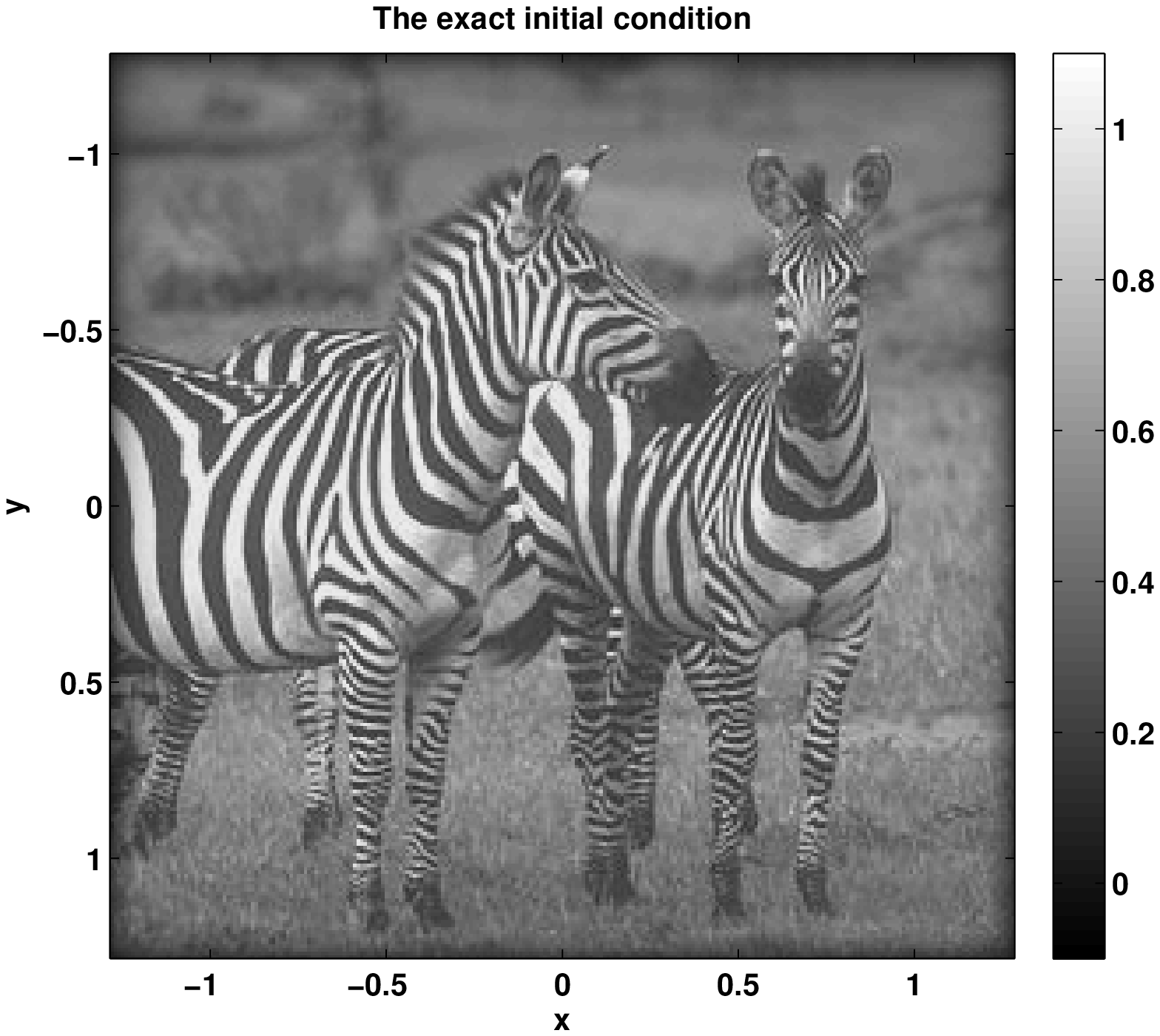,height=5.5cm }\\
(c)\epsfig{figure=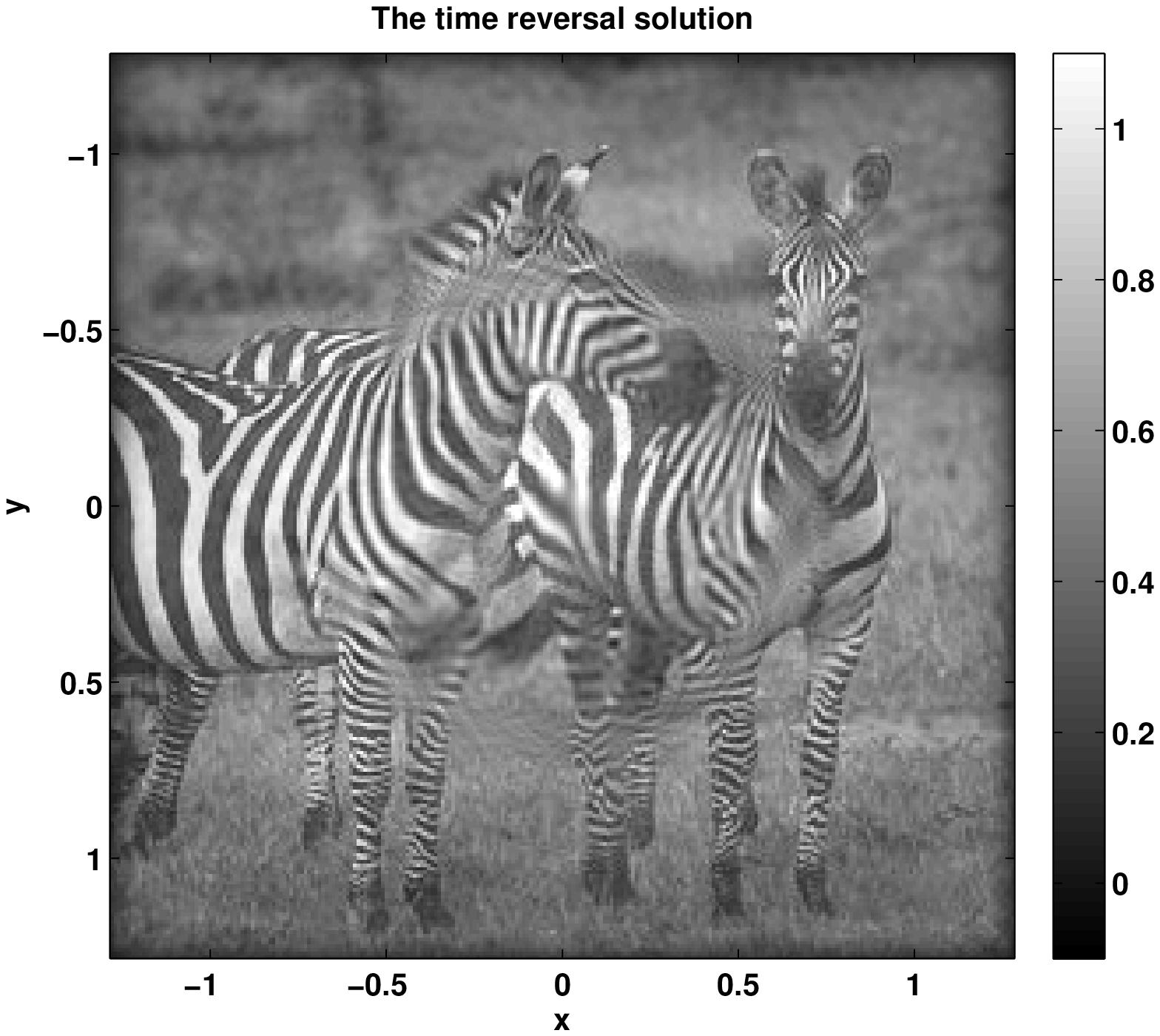,height=5.5cm }
(d)\epsfig{figure=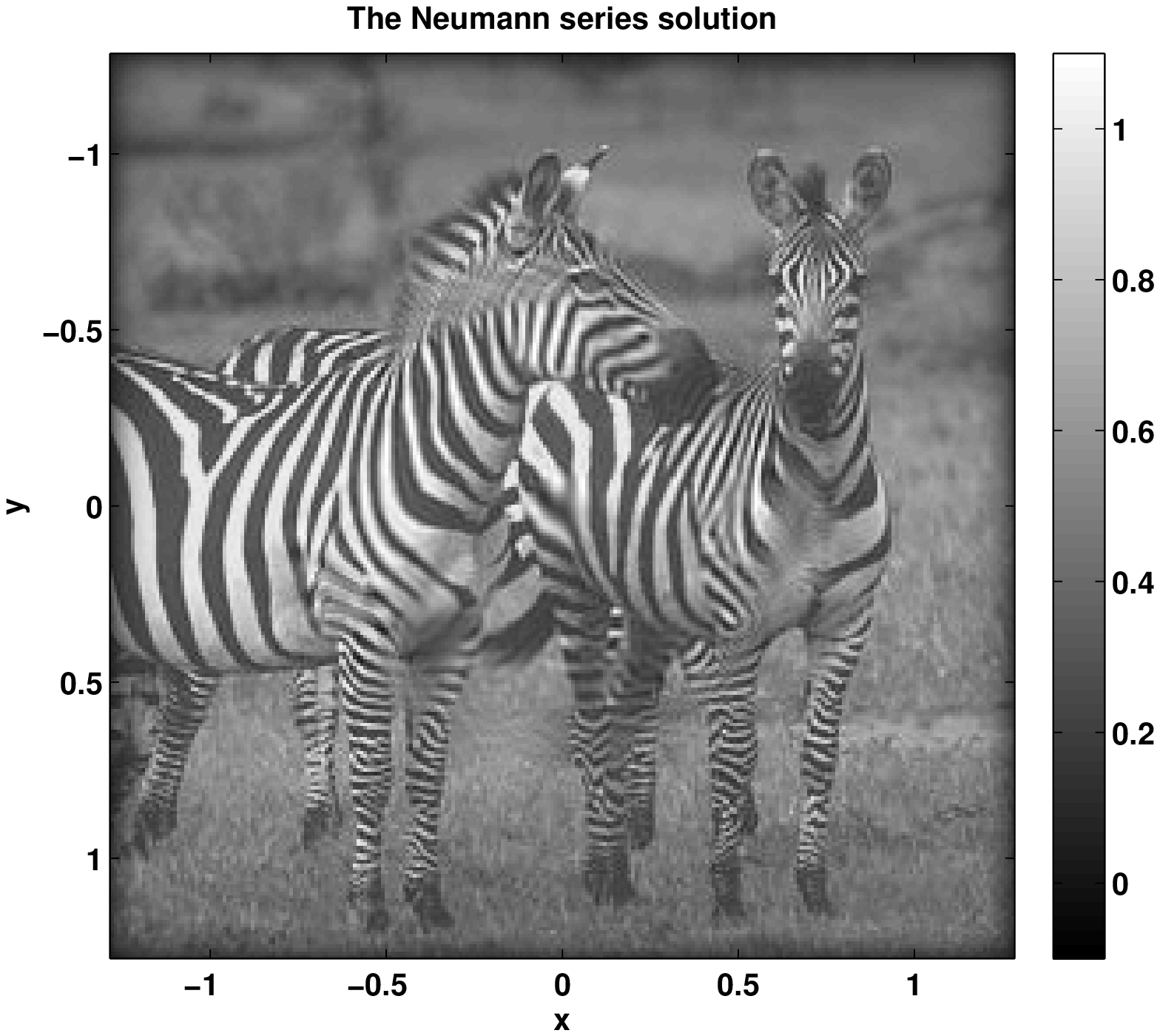,height=5.5cm }\\
\end{center}
\caption{Example 2 with the trapping speed $c_2$. $T=4T_0$. 
(a): the boundary distance map. 
(b): the exact initial condition. 
(c): the time reversal solution. 
(d): the Neumann series solution. 
}
\label{Fig:2dZebra4T0Speed63}
\end{figure}

\textbf{Figure \ref{Fig:2dZebra4T0Speed63}:} The sound speed is given by equation \eqref{key63} and we estimate $T_0$ to be $T_0\approx 2.1547$. It is trapping with the circles of radii approximately $0.23$ and $0.67$ being stable trapped rays. The invisible singularities are much more structured now and are in some neighborhoods of those two circles. As a result, jumps across radial and close to radial lines near those circles are affected more. Notice that the NS solution is much cleaner in the smooth parts close to the boundary, and the TR one is brighter than the original close to the center. This shows that the NS solutions reconstructs the low frequency modes better. 


\section{Numerical results: discontinuous sound speed}
\label{Sec:NumericalDis}



\begin{figure}
\begin{center}
(a)\epsfig{figure=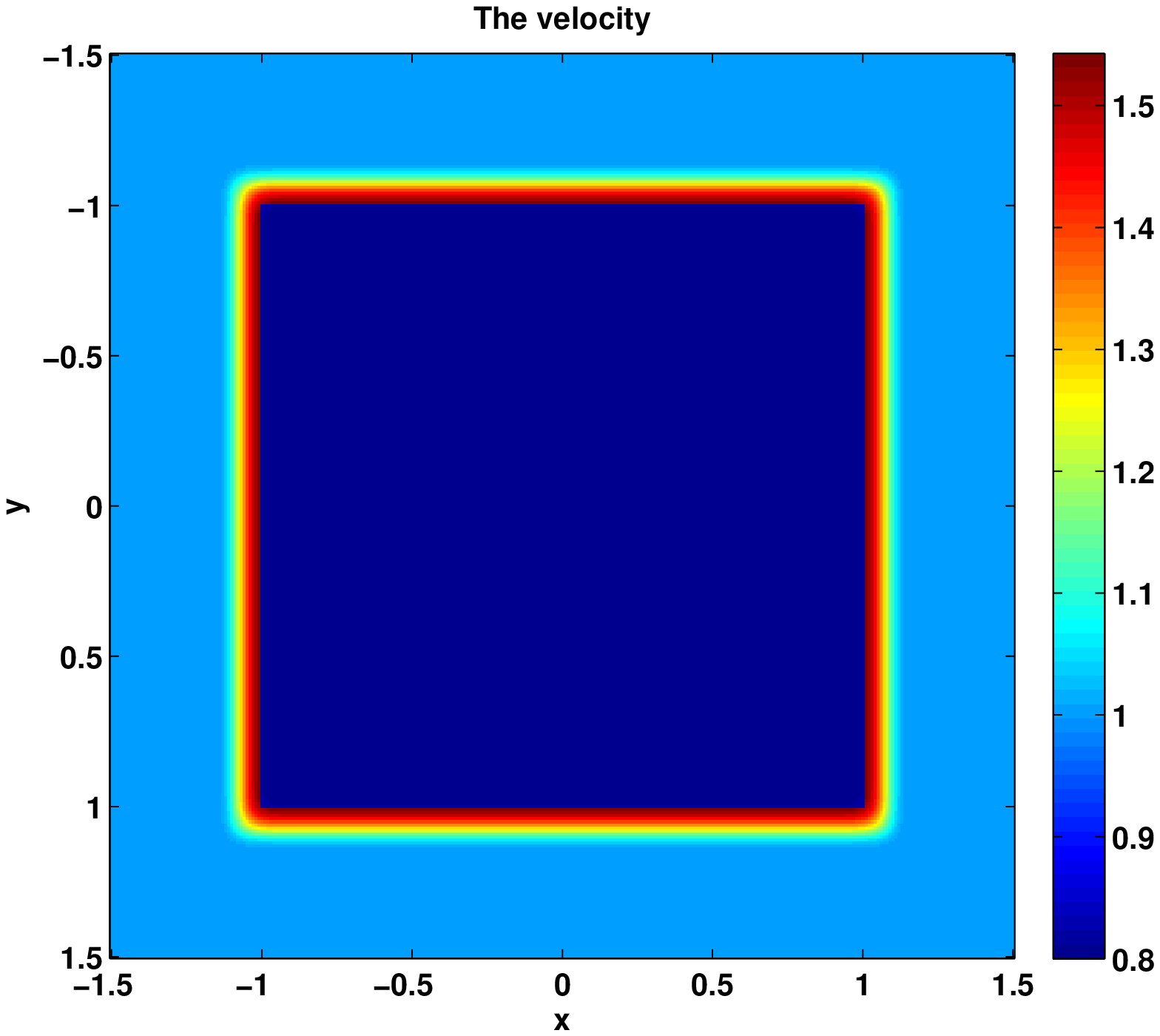,height=5.0cm}
(b)\epsfig{figure=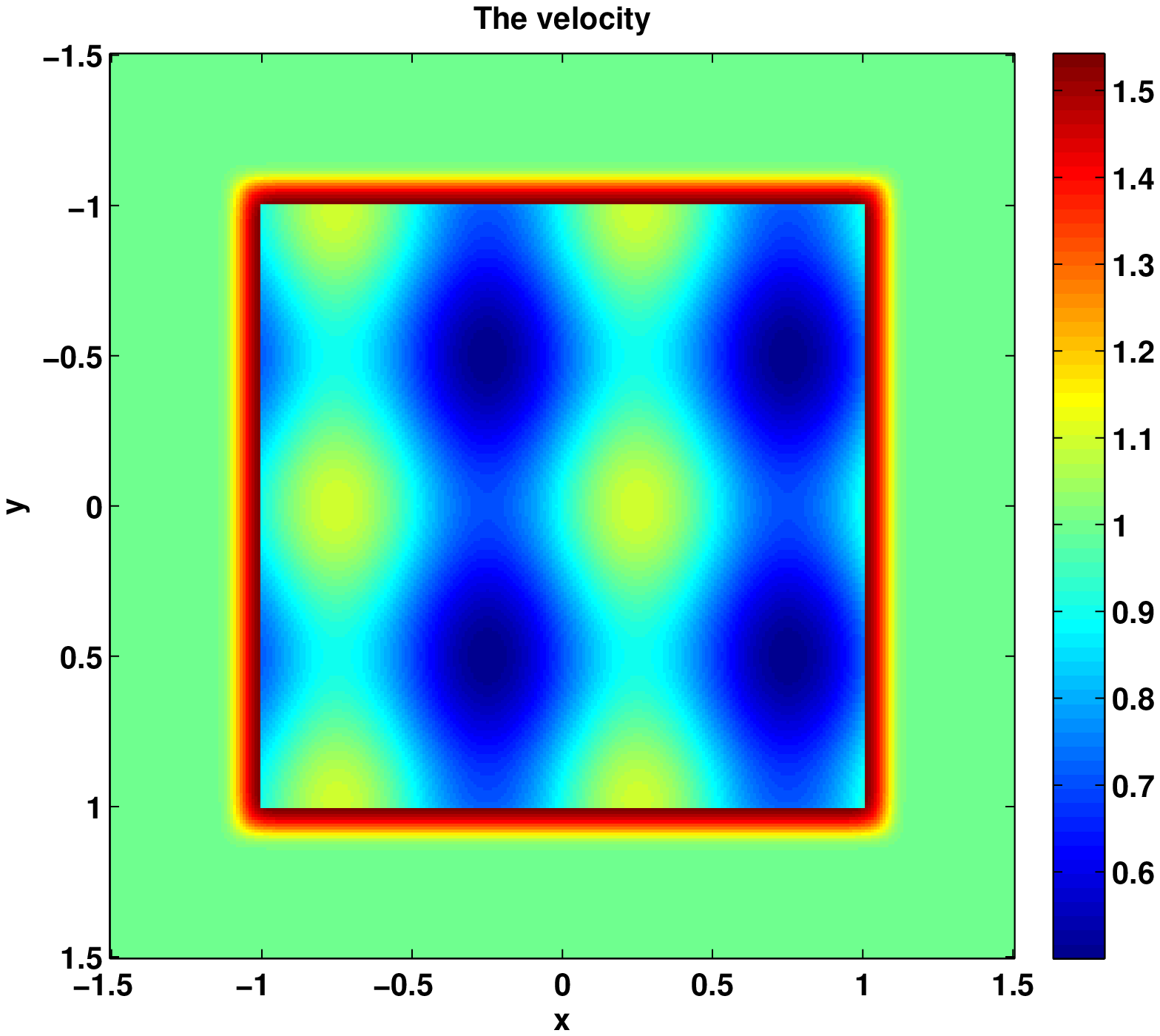,height=5.0cm}
\end{center}
\caption{Sound speed models.  (a): a discontinuous piecewise sound speed $c_4$. (b): a non- piecewise constant discontinuous sound speed $c_5$.}
\label{Fig:VelocitiesDis}
\end{figure}

We present here numerical examples with two discontinuous sound speeds illustrated in Figure \ref{Fig:VelocitiesDis}. The first one, that we denote by $c_4$, is equal to $0.8$ in the square $[-1,1]^2$, then jumps by about a factor of $2$ from the interior to the exterior, and then jumps to $1$. The second one, $c_5$, has similar jumps but inside the square $[-1,1]^2$ is variable, equal to the speed \r{key8}. 

If the speed jumps by a factor of $2$ when going out of a square, all rays hitting the boundary at an angle less than $60$ degrees are completely reflected. Those rays that hit the boundary at an angle greater than $30$ degrees generate a reflected ray hitting the boundary again without a transmitted component, etc. Therefore, all rays hitting the boundary of the smaller square at angles between $30$ and $60$ degrees are completely trapped for all times. 

\subsection{Example 3: The Shepp-Logan phantom, Figures~\ref{Fig:2dZebraDis1}--\ref{Fig:2dZebraDis2} }
\subsubsection{Piecewise constant discontinuous speed $c_4$}
\begin{figure}
\begin{center}
(a)\epsfig{figure=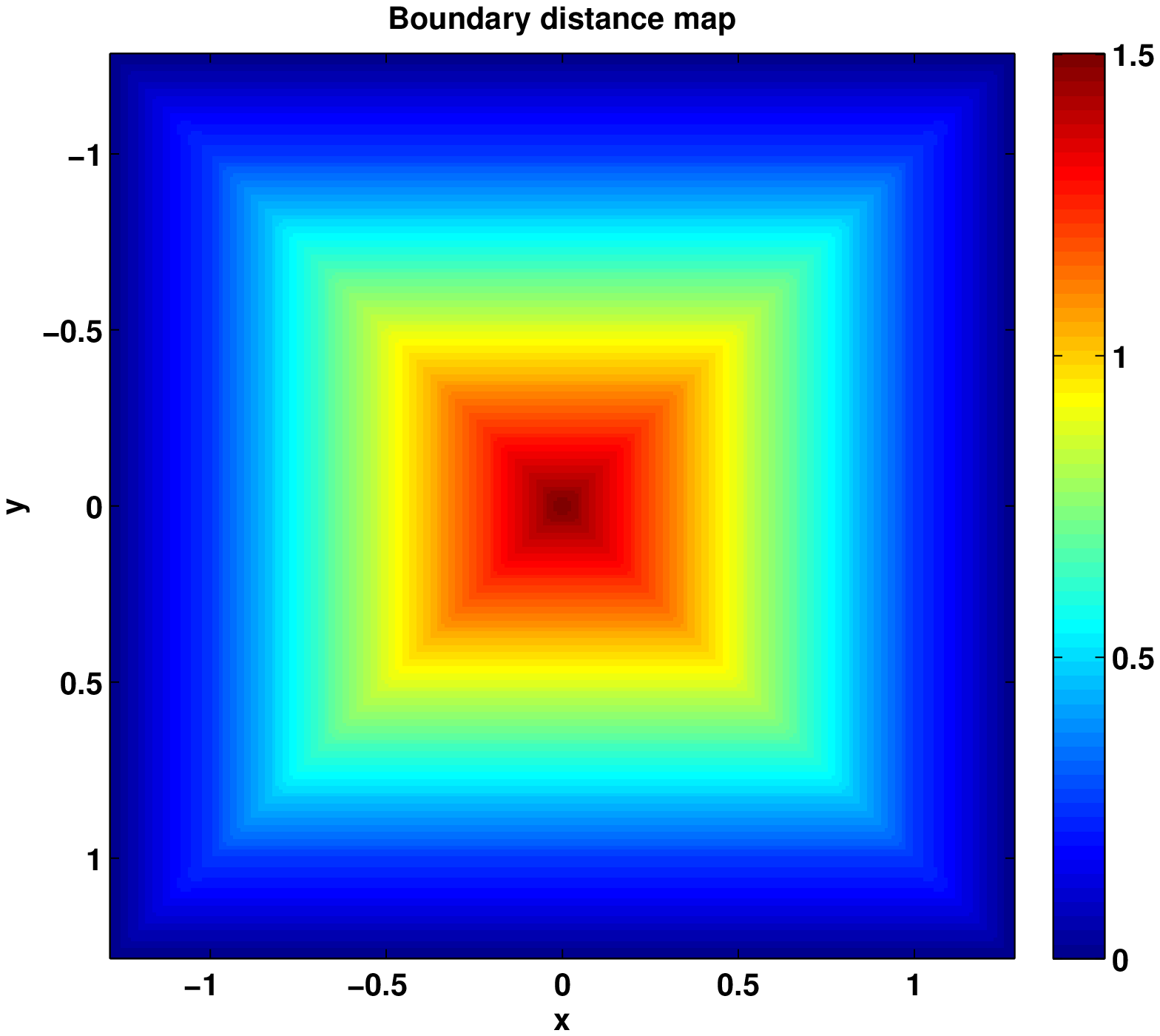,height=5.5cm }
(b)\epsfig{figure=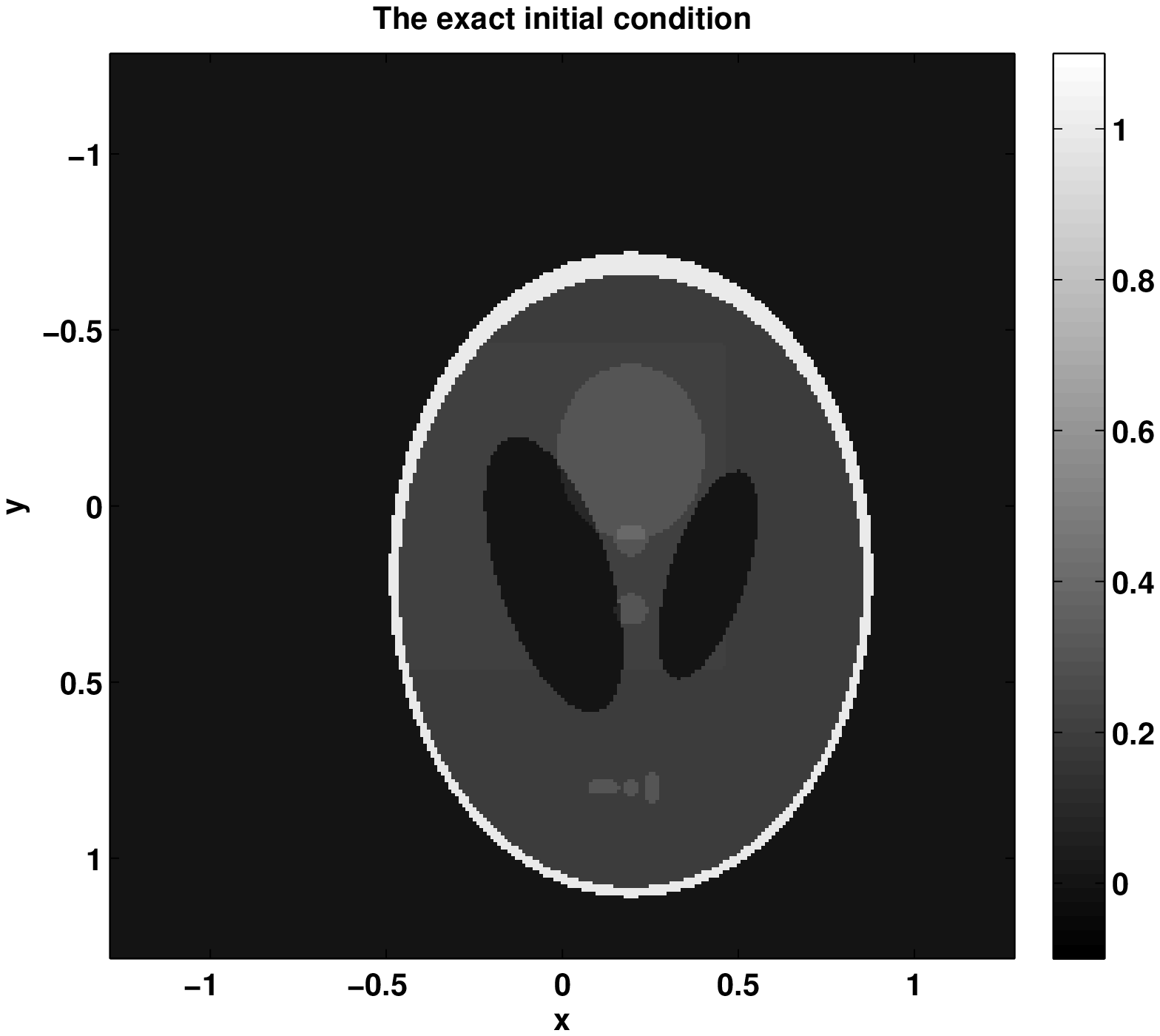,height=5.5cm }\\
(c)\epsfig{figure=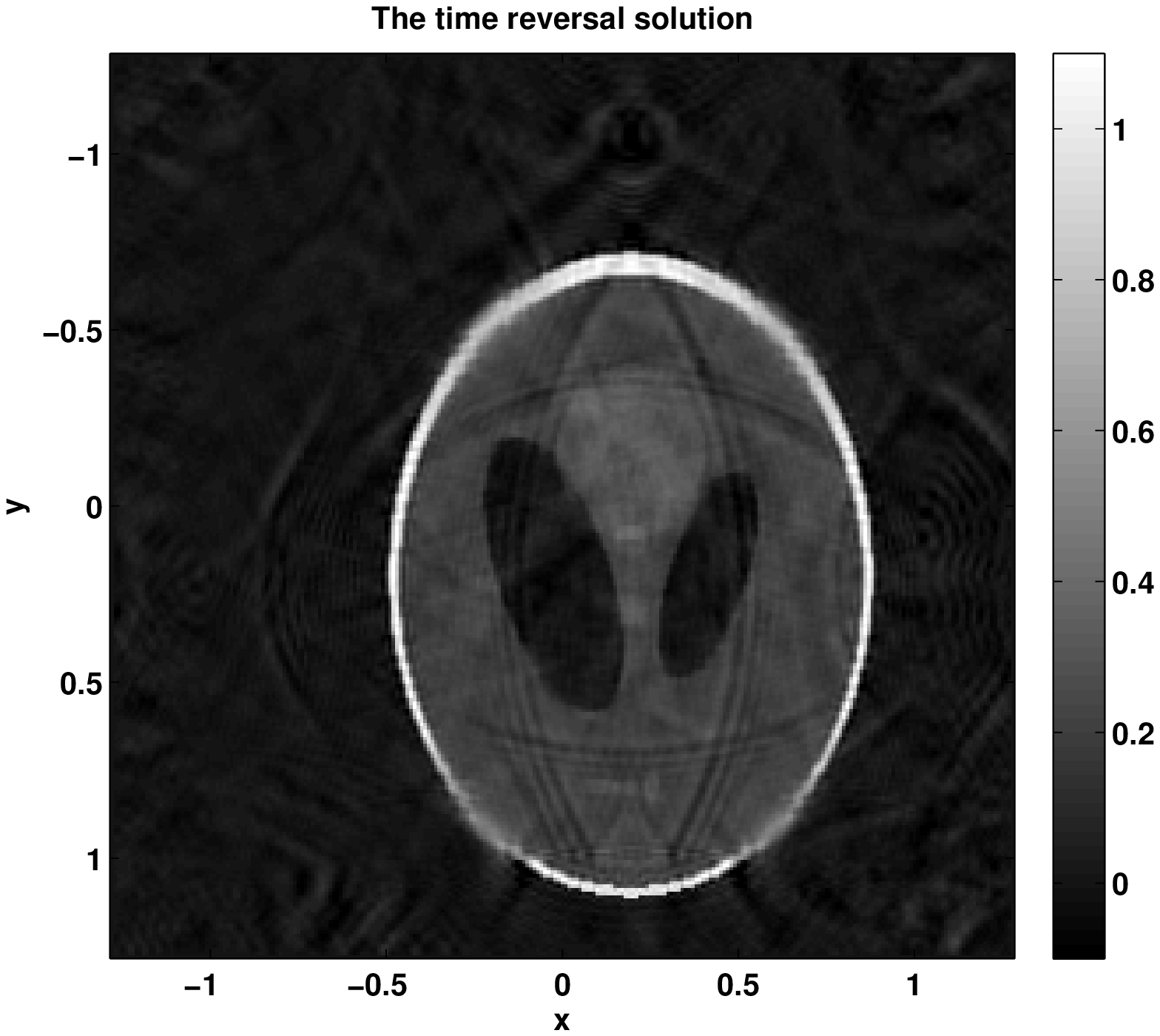,height=5.5cm }
(d)\epsfig{figure=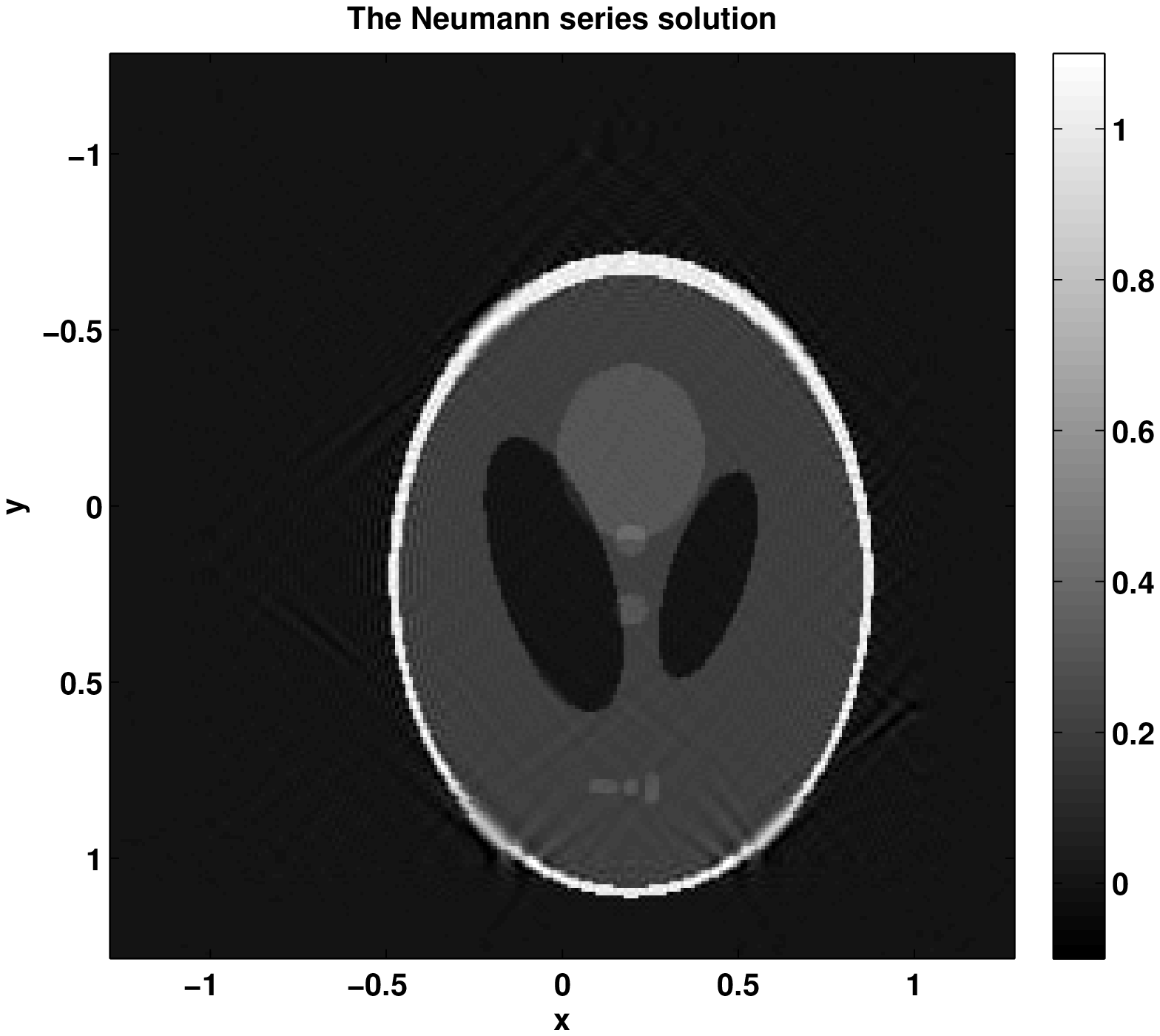,height=5.5cm }
\end{center}
\caption{Example 3 with the discontinuous sound speed $c_4$. $T=4T_0$.   
(a): the boundary distance map. 
(b): the exact initial condition. 
(c): the time reversal solution. 
(d): the Neumann series solution. 
 }
\label{Fig:2dZebraDis1}
\end{figure}

\textbf{Figure \ref{Fig:2dZebraDis1}:} 
The sound speed is $c_4$ given by Figure \ref{Fig:VelocitiesDis}(a).  
We have $T_0\approx 1.50$ with $T=4T_0$. The artifacts in the TR image are quite strong and can be explained by the way that singularities propagate in this case. In Figure~\ref{fig:skull2}, for example, the ray that exits on top carries a fraction of the energy only. When we reverse the time, at the first contact with the outer boundary of the ``skull'', this ray will create a reflected one (not seen in the figure) together with the transmitted one, shown there. The reflected one is not part of the actual graph that we are trying to invert. It will reflect off $\bo$, and go back to the interior of $\Omega$, creating more artificial rays, etc. If we had an infinite time $T$, then those artificial rays will be canceled by other such artificial rays, leading to an exact reconstruction, as $T\to\infty$ (at a very slow logarithmic rate for smooth $f$); this explains the artifacts in the TR image. The NS series expansion creates very few artifacts, and the few that are seen are mostly due to singularities near or on the ``skull''. 

\subsubsection{Discontinuous speed $c_5$}
\begin{figure}
\begin{center}
(a)\epsfig{figure=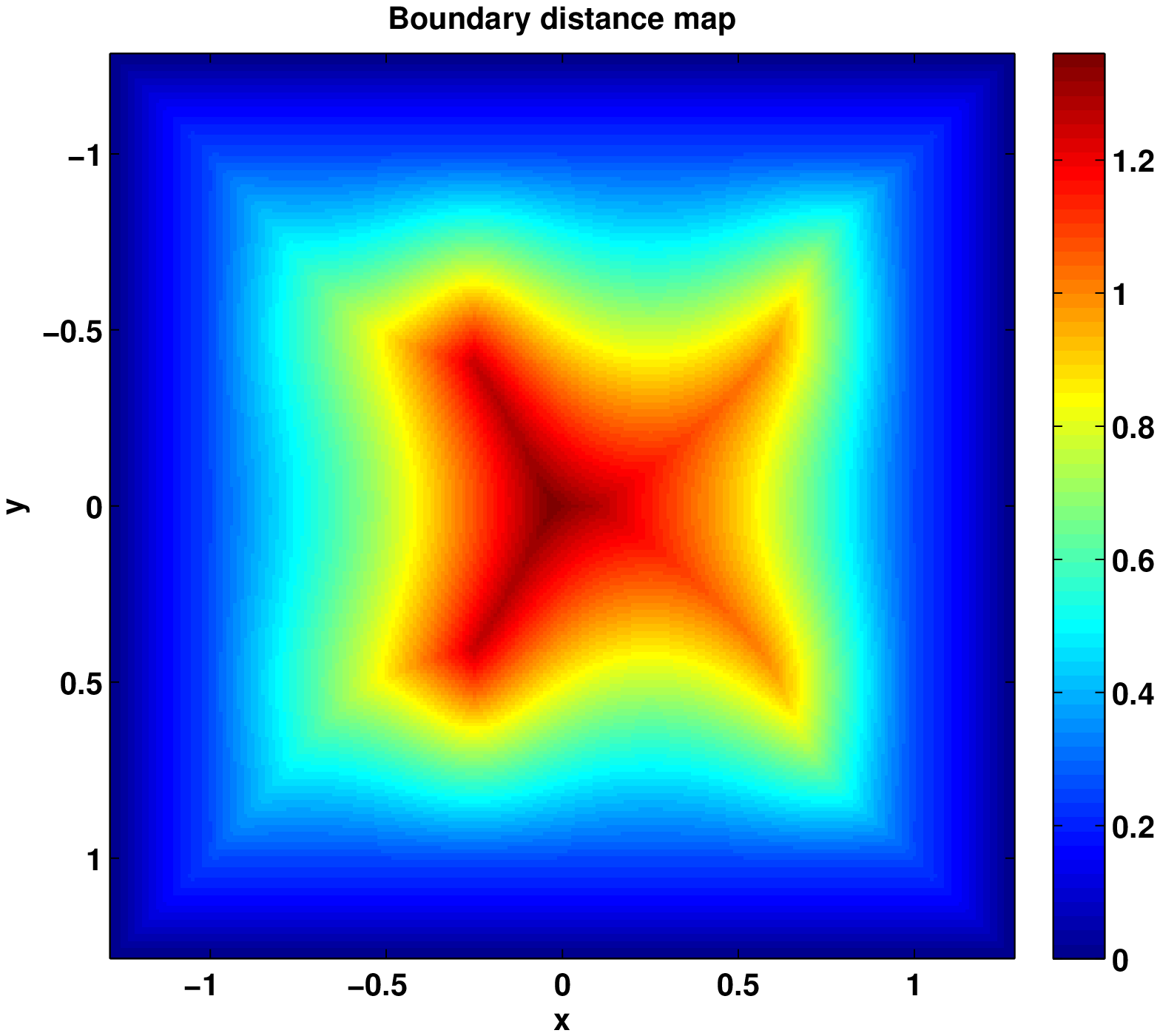,height=5.5cm }
(b)\epsfig{figure=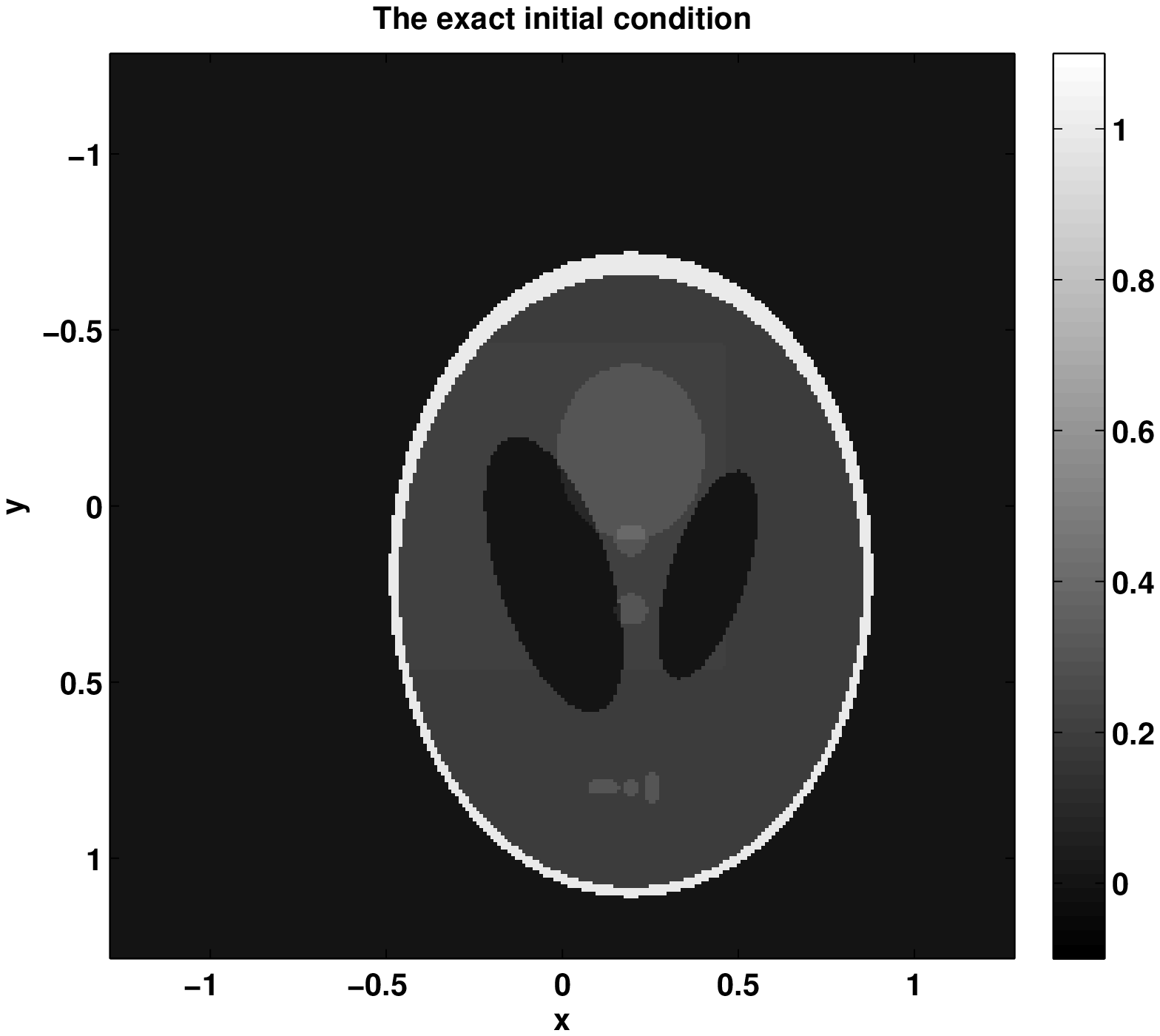,height=5.5cm }\\
(c)\epsfig{figure=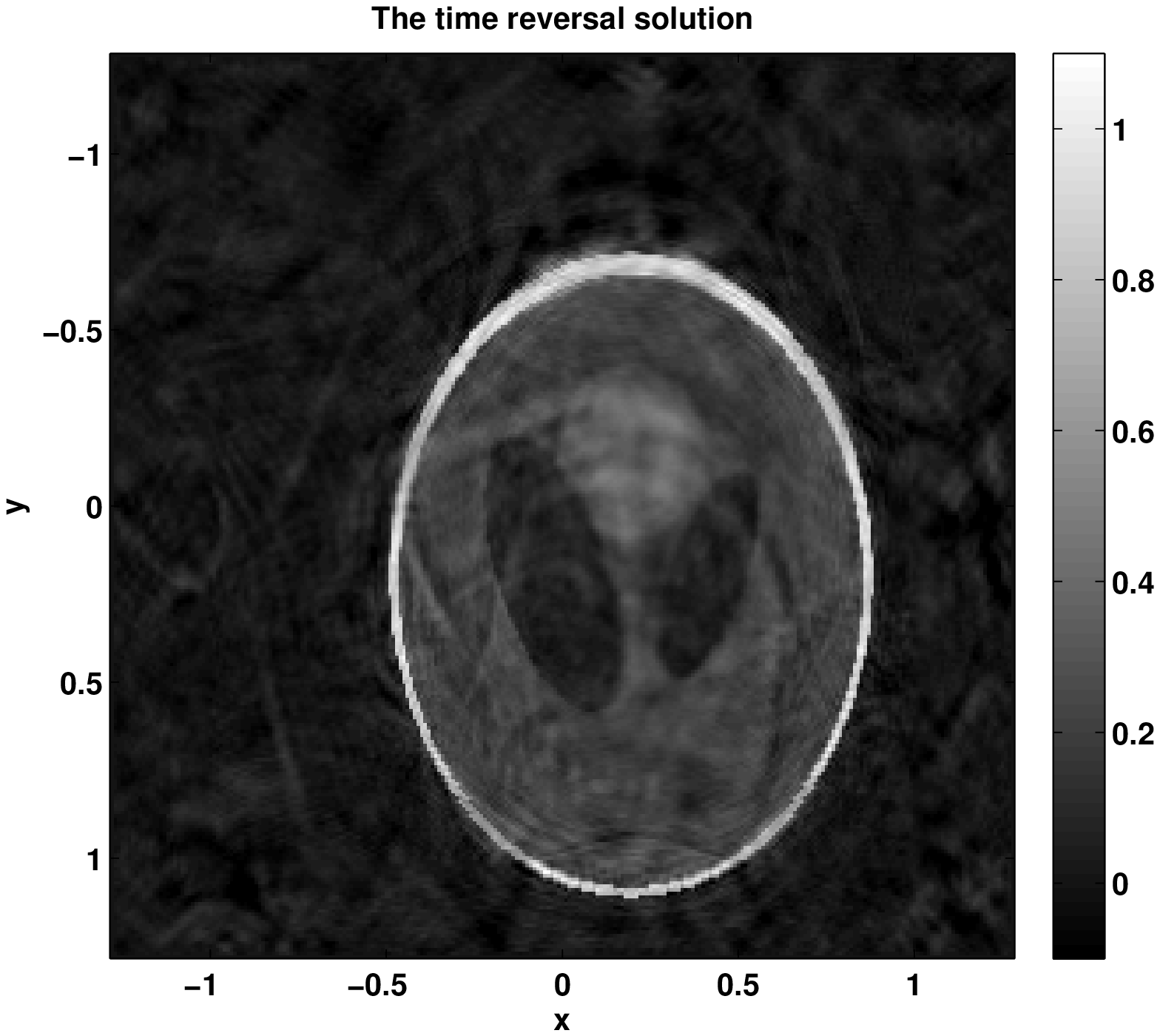,height=5.5cm }
(d)\epsfig{figure=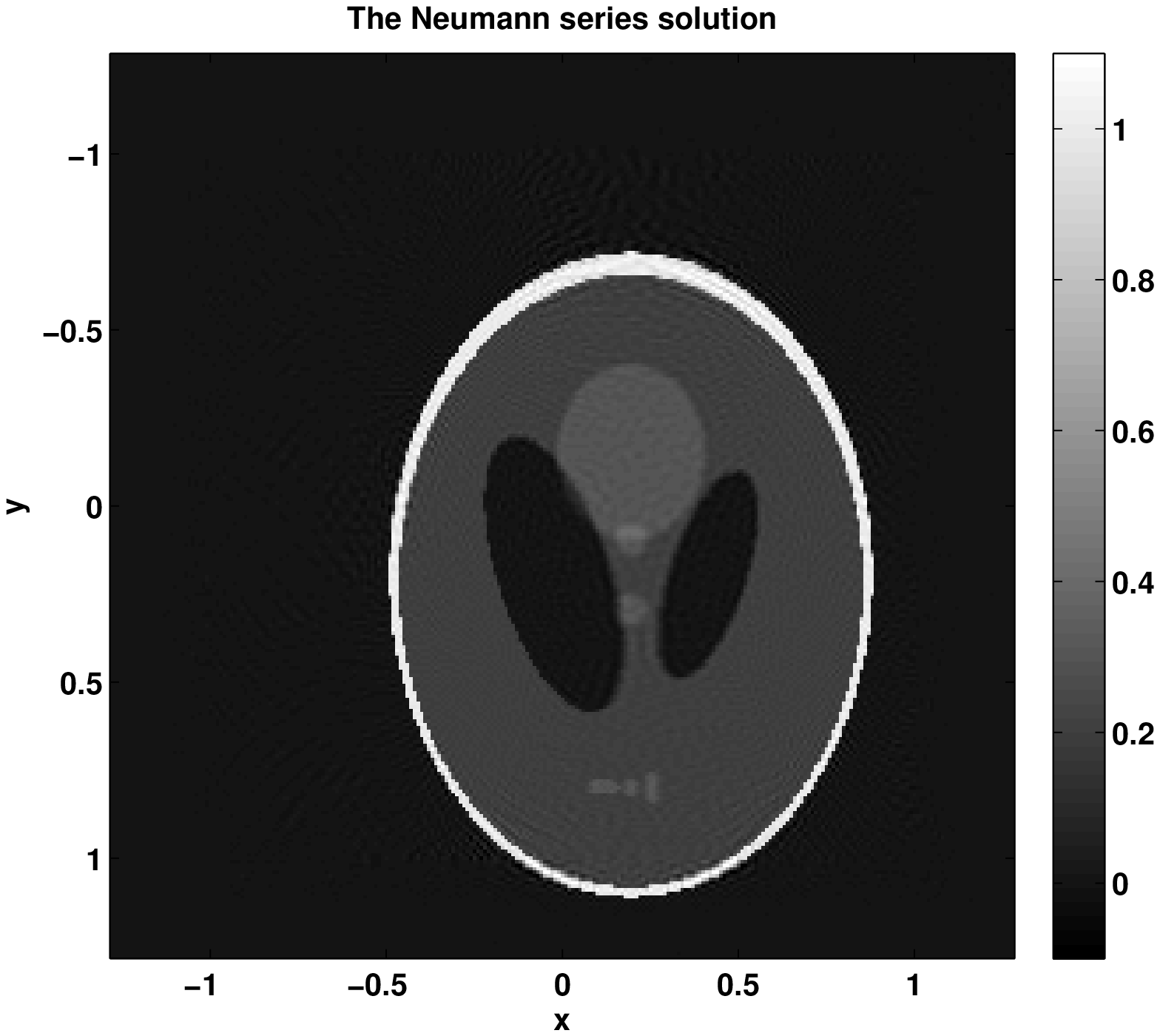,height=5.5cm }
\end{center}
\caption{Example 3 with the discontinuous sound speed $c_5$.  
(a): the boundary distance map. 
(b): the exact initial condition. 
(c): the time reversal solution. 
(d): the Neumann series solution. 
 }
\label{Fig:2dZebraDis2}
\end{figure}

\textbf{Figure \ref{Fig:2dZebraDis2}:} The sound speed $c_5$ is given by Figure \ref{Fig:VelocitiesDis}(b).  Here, 
$T_0\approx 1.36$,  $T=4T_0$. The speed is not so symmetric anymore and the artifacts in the TR image are still strong but more random. The variable speed inside improves the NS image. 

\begin{figure}
\begin{center}
(a)\epsfig{figure=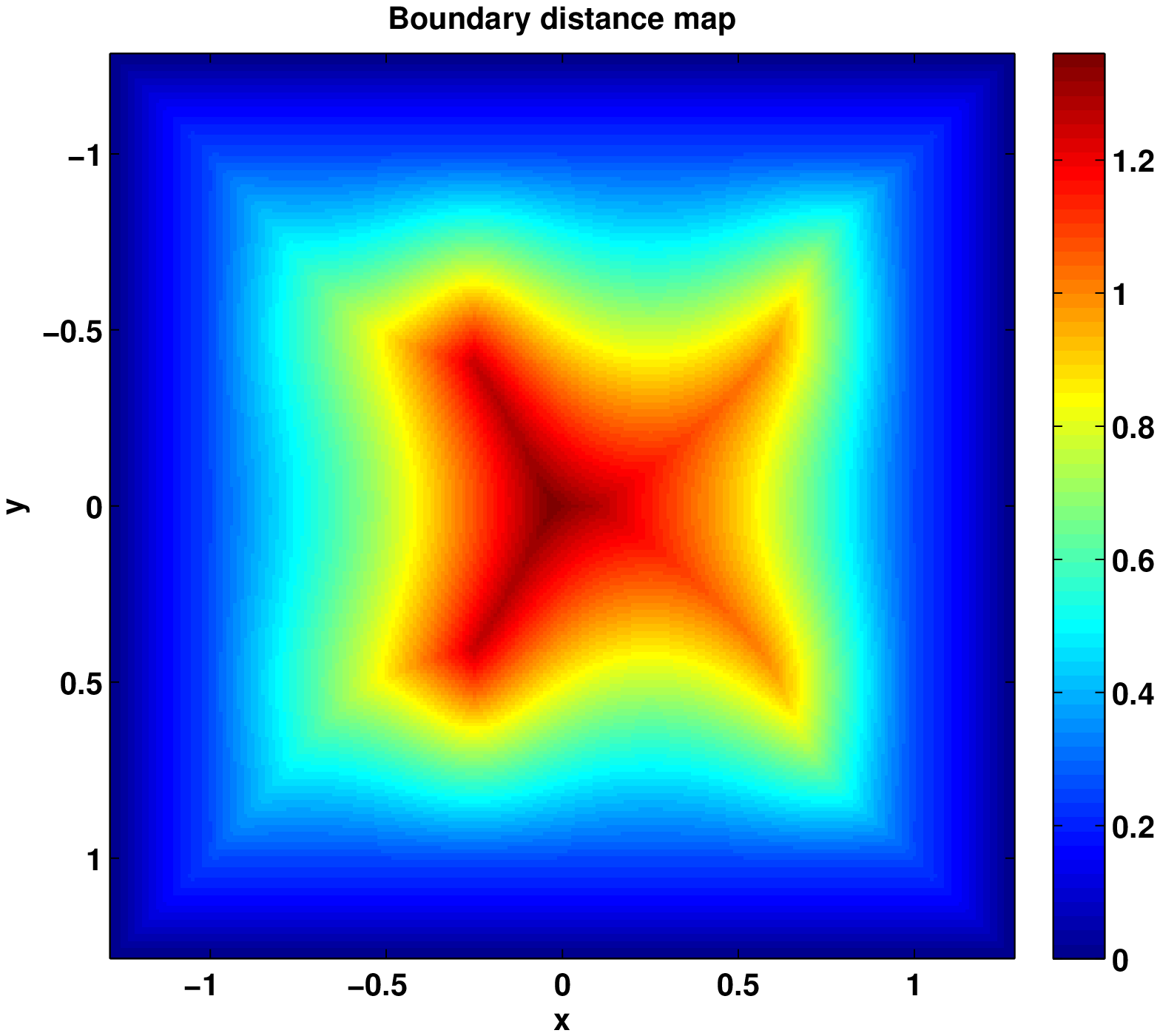,height=5.5cm }
(b)\epsfig{figure=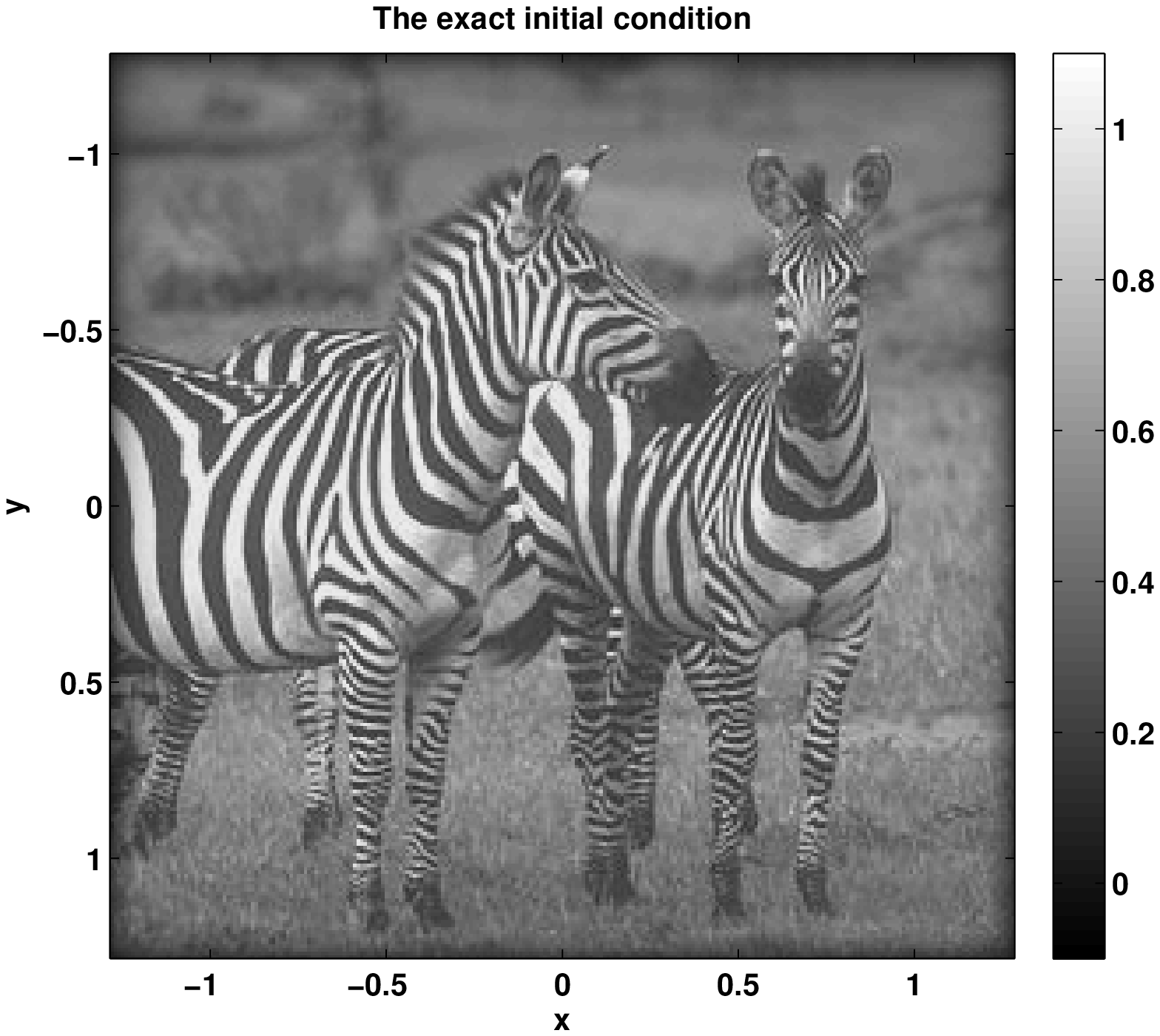,height=5.5cm }\\
(c)\epsfig{figure=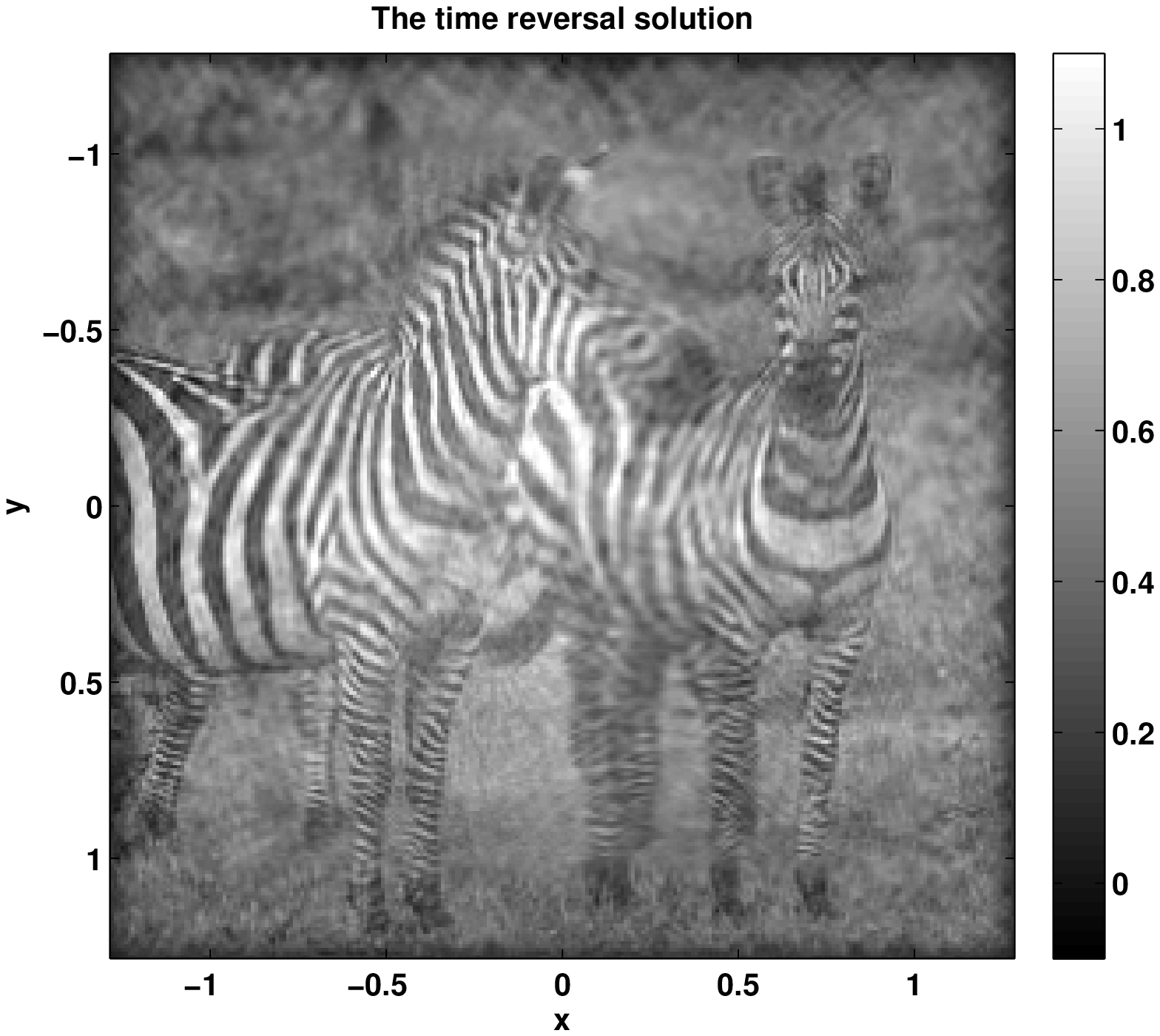,height=5.5cm }
(d)\epsfig{figure=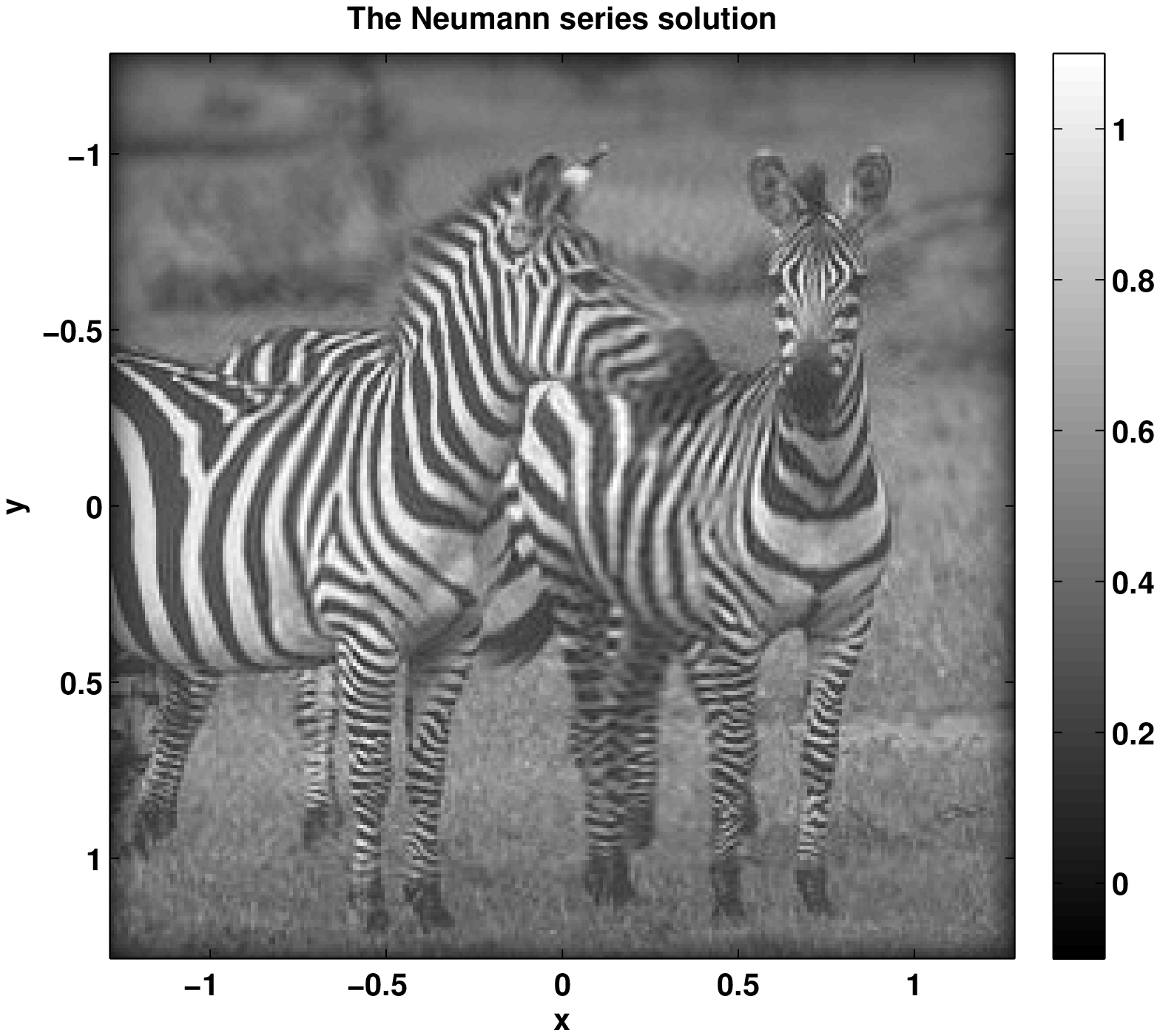,height=5.5cm }
\end{center}
\caption{Example 4 with the discontinuous speed $c_5$: Figure \ref{Fig:VelocitiesDis}(b). 
(a): the boundary distance map. 
(b): the exact initial condition. 
(c): the time reversal solution. 
(d): the Neumann series solution. 
}
\label{Fig:2dZebraDis3}
\end{figure}

\subsection{Example 4: Zebras, Figure~\ref{Fig:2dZebraDis3}}
\textbf{Figure~\ref{Fig:2dZebraDis3}:} 
The sound speed $c_5$ is given by Figure \ref{Fig:VelocitiesDis}(b).
Here, 
$T_0\approx 1.36$, and we take  $T=4T_0$. As in Figure~\ref{Fig:2dZebraDis2}, the TR reconstruction (error $21.83\%$) contains many wrong singularities, while  the NS image (error $8\%$, $k=16$) is very clean.

\section{Numerical results: partial data, Figures~\ref{Fig:2dwaveT4pNotrapPartialA}--\ref{Fig:2dZebra4T0trapPartialOne}}
\label{Sec:NumericalPartial}

We present here numerical examples with data given on three or two sides only. In the first case, we remove data on the right hand side of the square; in the second one, we remove data on the right hand side and on the bottom of the square. We also use a smooth cutoff. We present examples with the non-trapping and the trapping speeds that we used before. Note that the notion of trapping changes with partial data since the stability depends whether all singularities can reach $\Gamma$ for times $\le T$. Still, removing the data on some of the sides in addition to the much longer times that signals need to reach $\Gamma$ produces worse reconstructions when the speed is trapping. 

\subsection{Partial data, non-trapping speed $c_1$, Figure~\ref{Fig:2dwaveT4pNotrapPartialA}}
\begin{figure} 
\begin{center}
(a)\epsfig{figure=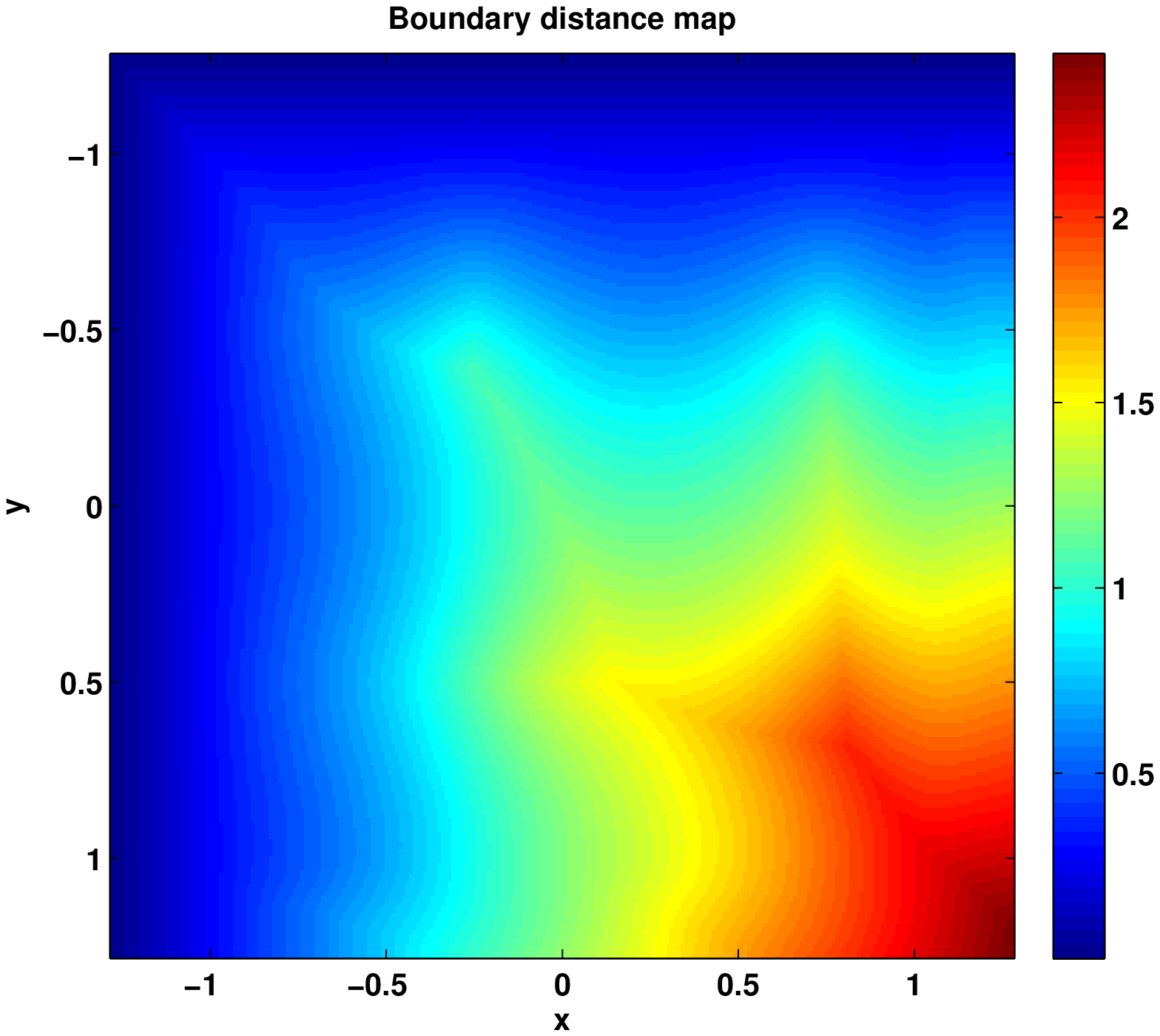,height=5.5cm }
(b)\epsfig{figure= 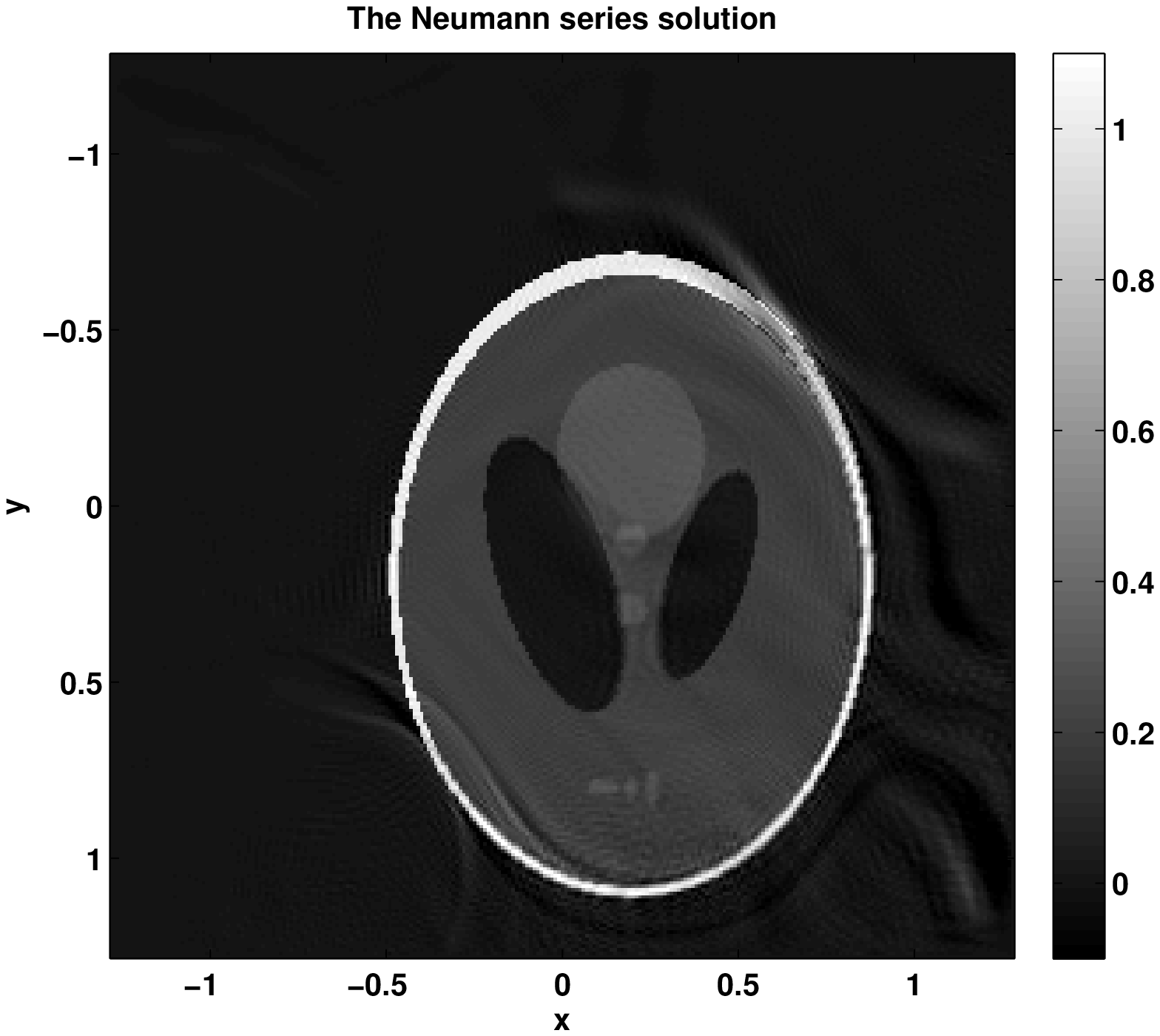,height=5.50cm }
\end{center}
\caption{Non-trapping speed $c_1$ with partial data on two adjacent sides. $T=4.7$. 
(a): the boundary distance map, data on two sides. 
(b): The Shepp-Logan phantom reconstruction. 
}
\label{Fig:2dwaveT4pNotrapPartialA}
\end{figure}

The speed $c_1$ is given by equation \eqref{key8}. In all three cases, $T=4.7$. 
\subsubsection{The Shepp-Logan phantom, Figure~\ref{Fig:2dwaveT4pNotrapPartialA} (a), (b)} 
We estimate $T_0$ to be $T_0\approx 2.5$. It is greater than before because $\Gamma$ is not the whole $\bo$ anymore. 
We have data on two adjacent sides. If the speed were constant, the singularities below the lower-left-to-upper-right diagonal that would exit in the sides without measurements are invisible. That would affect jumps across surfaces in that triangle with normals parallel or nearly parallel to that diagonal, roughly speaking. The speed is variable but not too far from constant, and we see the expected behavior. Note that this is a unstable case, and our analysis does not exclude an even exponential divergence (in the low frequency part) but the error gets smaller up to $k=8$ when the computation is stopped. 

\subsubsection{Zebras, Figure~\ref{Fig:2dwaveT4pNotrapPartialB} (a), (b)}
We have data on three sides.  In this case, $T_1\approx 1.35$, and we would have stability if the speed were constant. However, the speed is not constant but in some sense, not too far from a constant one. There might be a small set of geodesics that enter through the right-hand side and exit there as well, thus creating instability. The reconstruction is very good, however, with an error of $6.11\%$, $k=10$. The chosen $T$ is slightly greater than what would be the stability time, but the result is computed without the contribution of those geodesics. 

\begin{figure} 
\begin{center}
(a) \hspace{-0.12in} \epsfig{figure=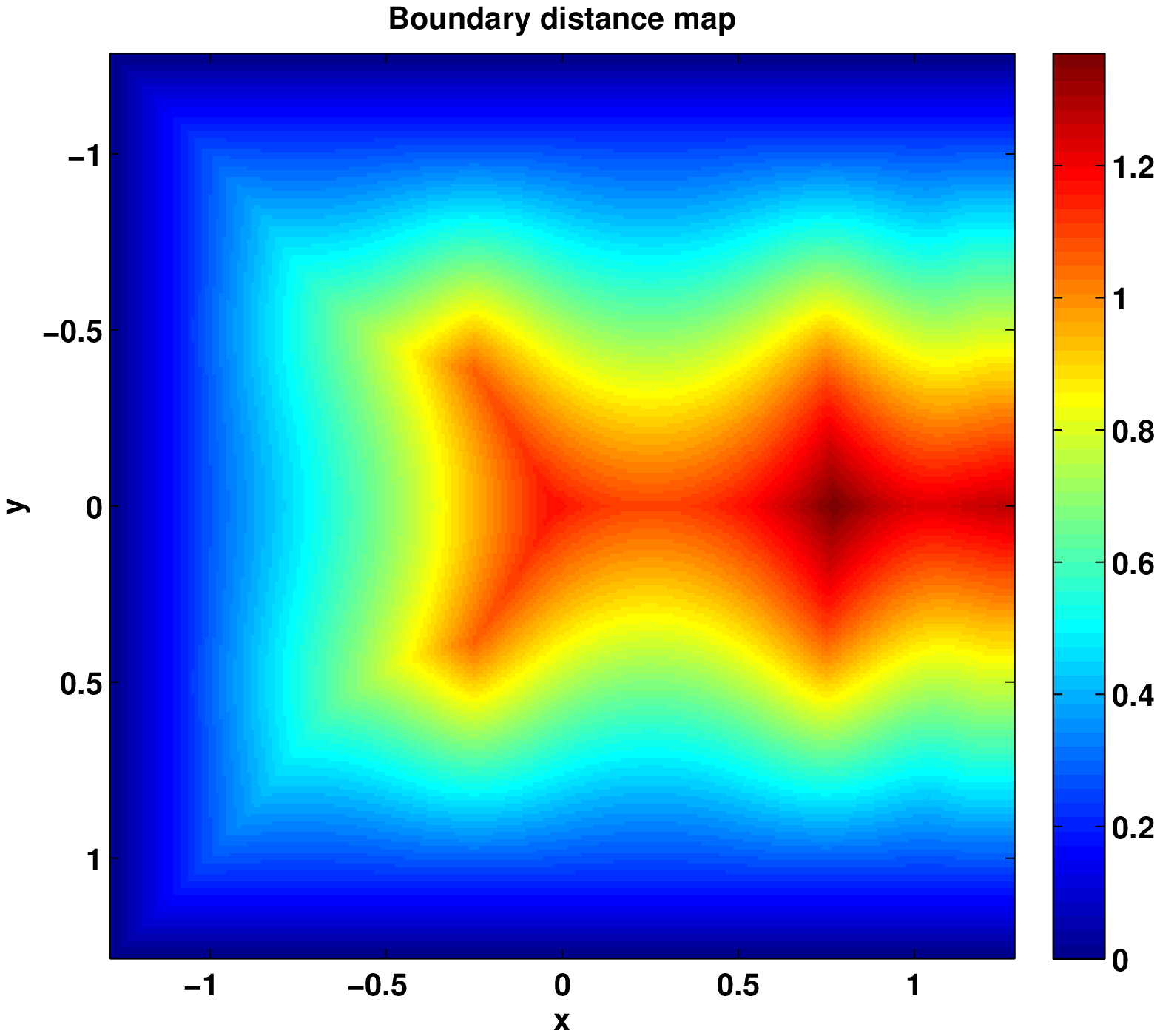,height=5.5cm }
(b)\hspace{-0.07in} \epsfig{figure= 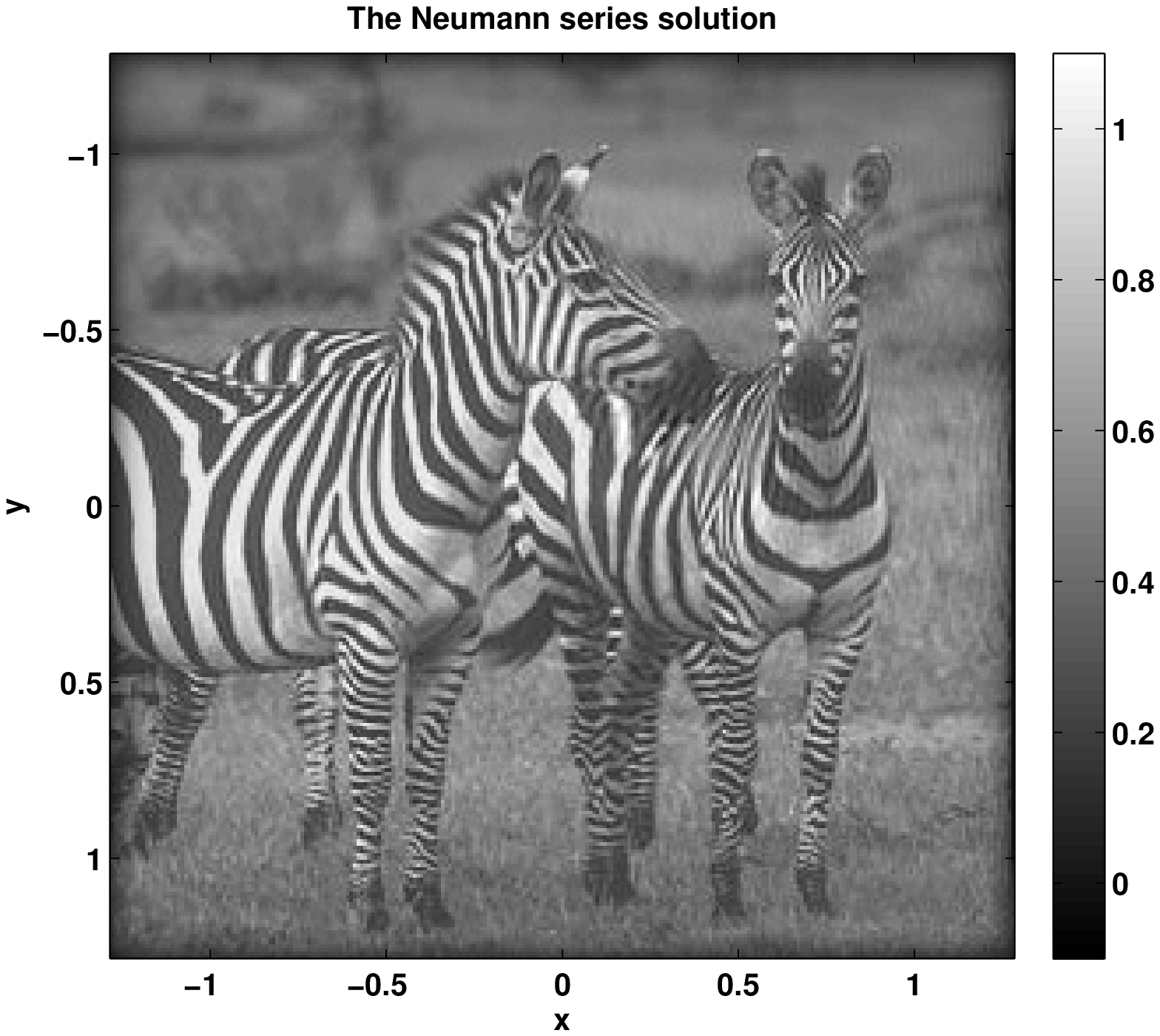,height=5.5cm }\\
(c)\epsfig{figure=Fig/Key8Keyini6341Sides105Traveltime.ps,height=5.5cm }
(d)\epsfig{figure=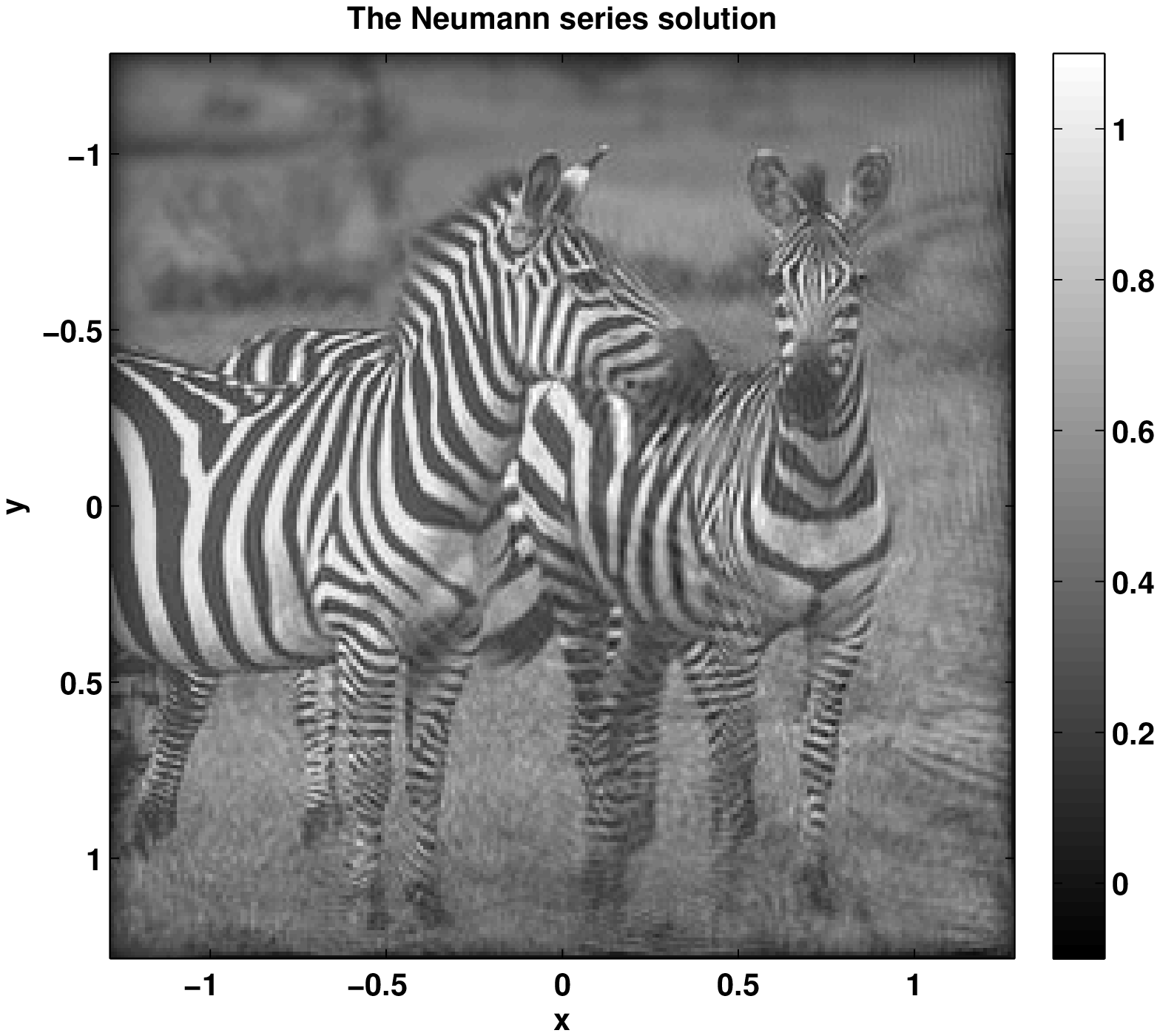,height=5.5cm }
\end{center}
\caption{A modified example with the non-trapping speed $c_1$, partial data.  $T=4.7$.
(a): the boundary distance map, data on three sides. 
(b): the reconstructed ``zebras'' image.  
(c): the boundary distance map, data on two sides. 
(d): the reconstructed ``zebras'' image. 
}
\label{Fig:2dwaveT4pNotrapPartialB}
\end{figure}

\subsubsection{Zebras, Figure~\ref{Fig:2dwaveT4pNotrapPartialB} (c), (d)}
We use  data on two  sides.   This is a very unstable case and one can see artifacts where the invisible singularities lie --- in the lower right-hand triangle, with slopes close to $1$, roughly speaking (jumps across curves with slopes in a neighborhood of $-1$). The error is $10.6\%$, and $k=10$. 

\begin{figure} 
\begin{center}
(a) \hspace{-0.12in} \epsfig{figure=  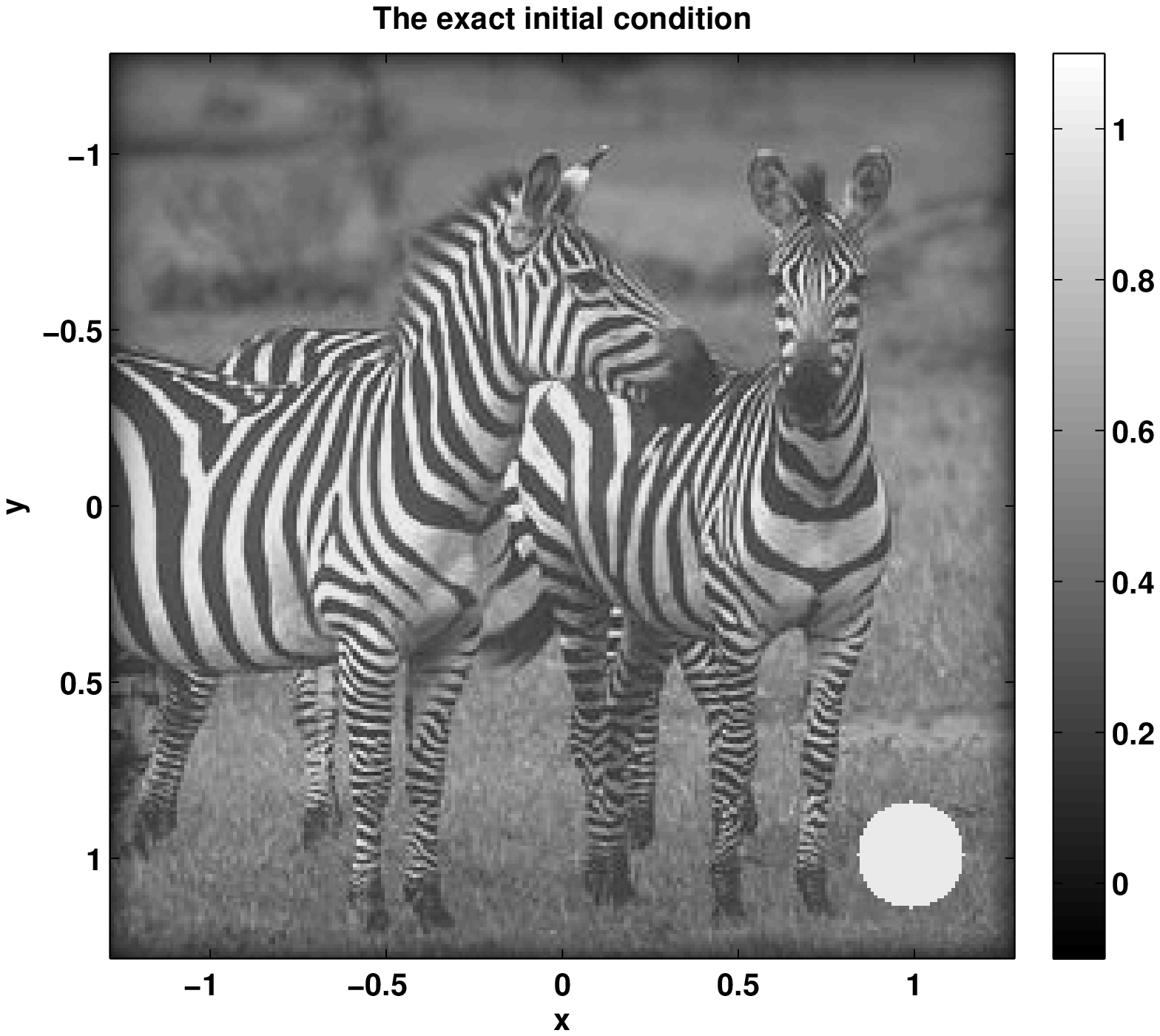,height=5.5cm }
(b)\hspace{-0.07in} \epsfig{figure=  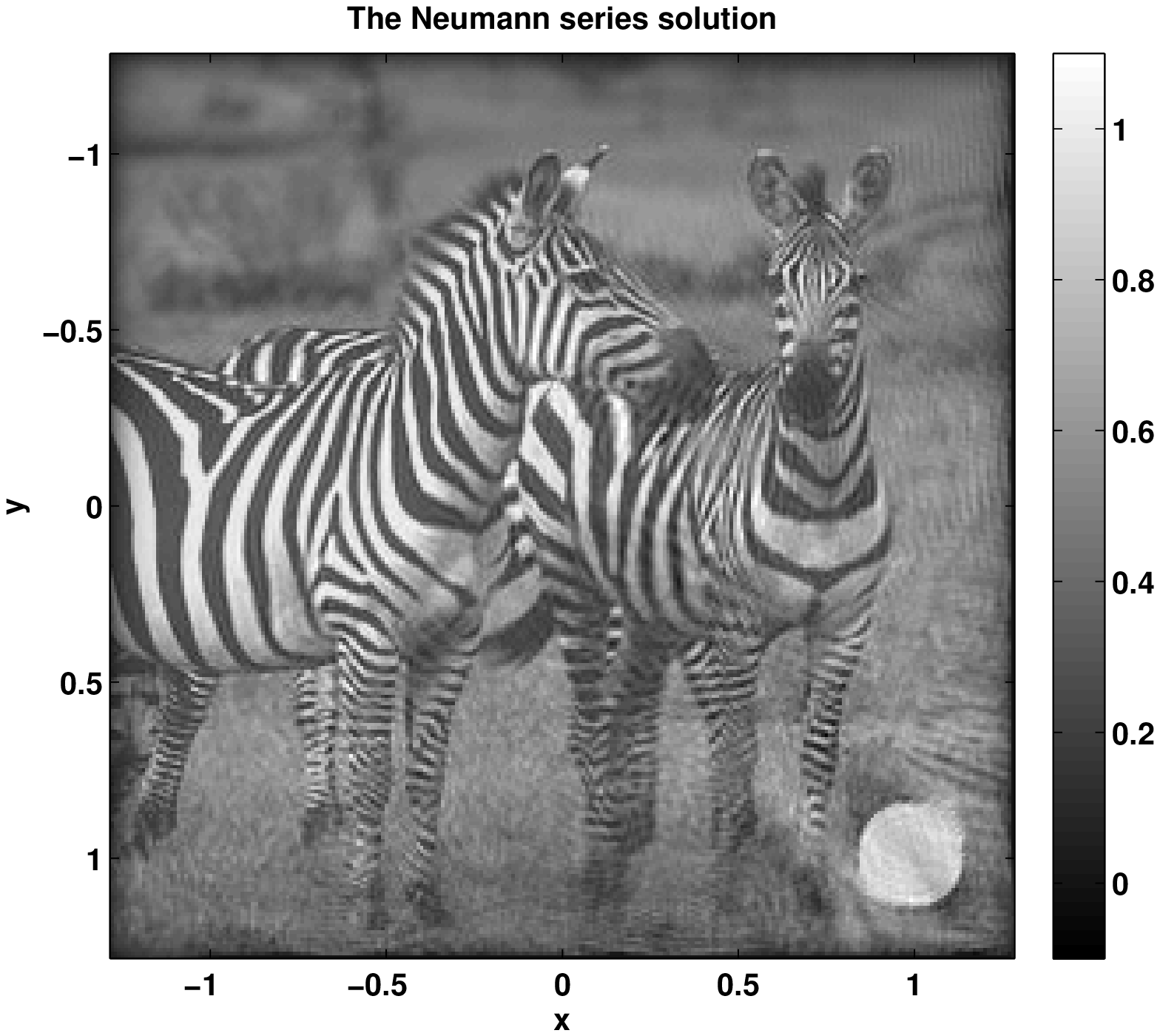,height=5.5cm }
\end{center}
\caption{Non-trapping speed $c_1$ examples, partial data.  $T=4.7$.
(a): the exact initial condition. 
(b): the reconstructed ``zebras'' image.  
}
\label{Fig:2dwaveT4pNotrapPartialC}
\end{figure}

\subsubsection{Zebras, Figure~\ref{Fig:2dwaveT4pNotrapPartialC} (a), (b)}
The zebras image reconstructed in Figure~\ref{Fig:2dwaveT4pNotrapPartialB} (d)  does not have so many invisible singularities in this particular setup (data on two sides), however.  For this reason, we present another example, Figure~\ref{Fig:2dwaveT4pNotrapPartialB} (a)(b), with a modified image that shows the expected behavior of the visible and invisible singularities.

\subsection{Partial data, trapping sound speeds $c_2$ and $c_2$, Figure~\ref{Fig:2dZebra4T0trapPartialOne}} 
The sound speed is $c_3$ (the first two),  and $c_2$ (the third example). 


\subsubsection{Figure~\ref{Fig:2dZebra4T0trapPartialOne} (a), (b)}
Here,  $T=4.93>T_0\approx 1.3$,  with data on three sides. The ``chaotic'' trapping speed $c_3$ makes the reconstruction worse than before. This is an unstable case because the trapping speed leaves many singularities invisible. As expected, the worst part is near the side with no observations due to geodesics that enter and exit through that side. There are invisible singularities everywhere, as well, due to the speed. 

\subsubsection{Figure~\ref{Fig:2dZebra4T0trapPartialOne} (c), (d)} 
Here,  $T=4.93>T_0\approx 2.1$,  with data on two sides, the same speed as above. As expected, the reconstruction is quite bad near the sides with no data. 

\begin{figure} 
\begin{center}
(a)\epsfig{figure=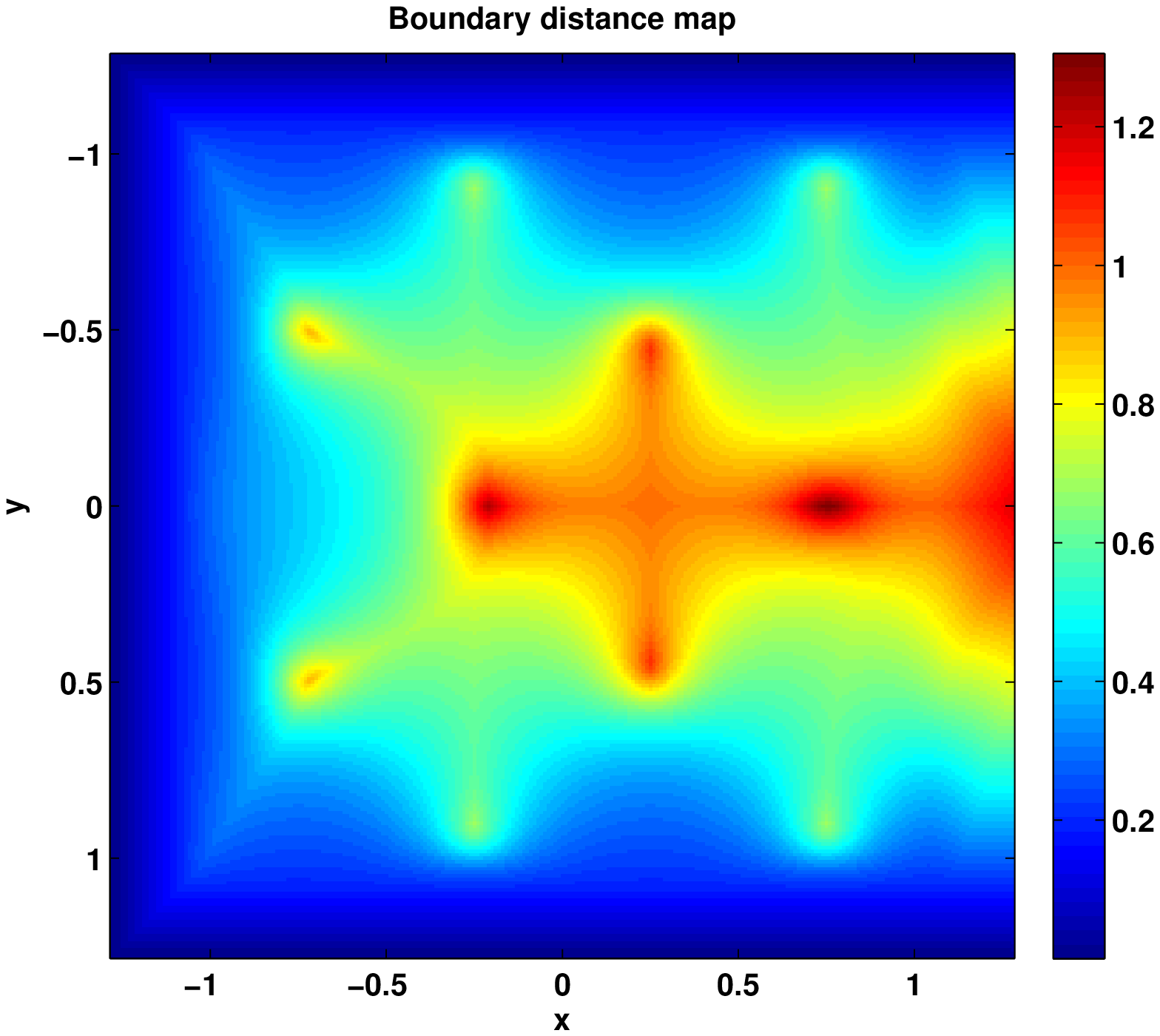,height=5.5cm }
(b)\epsfig{figure=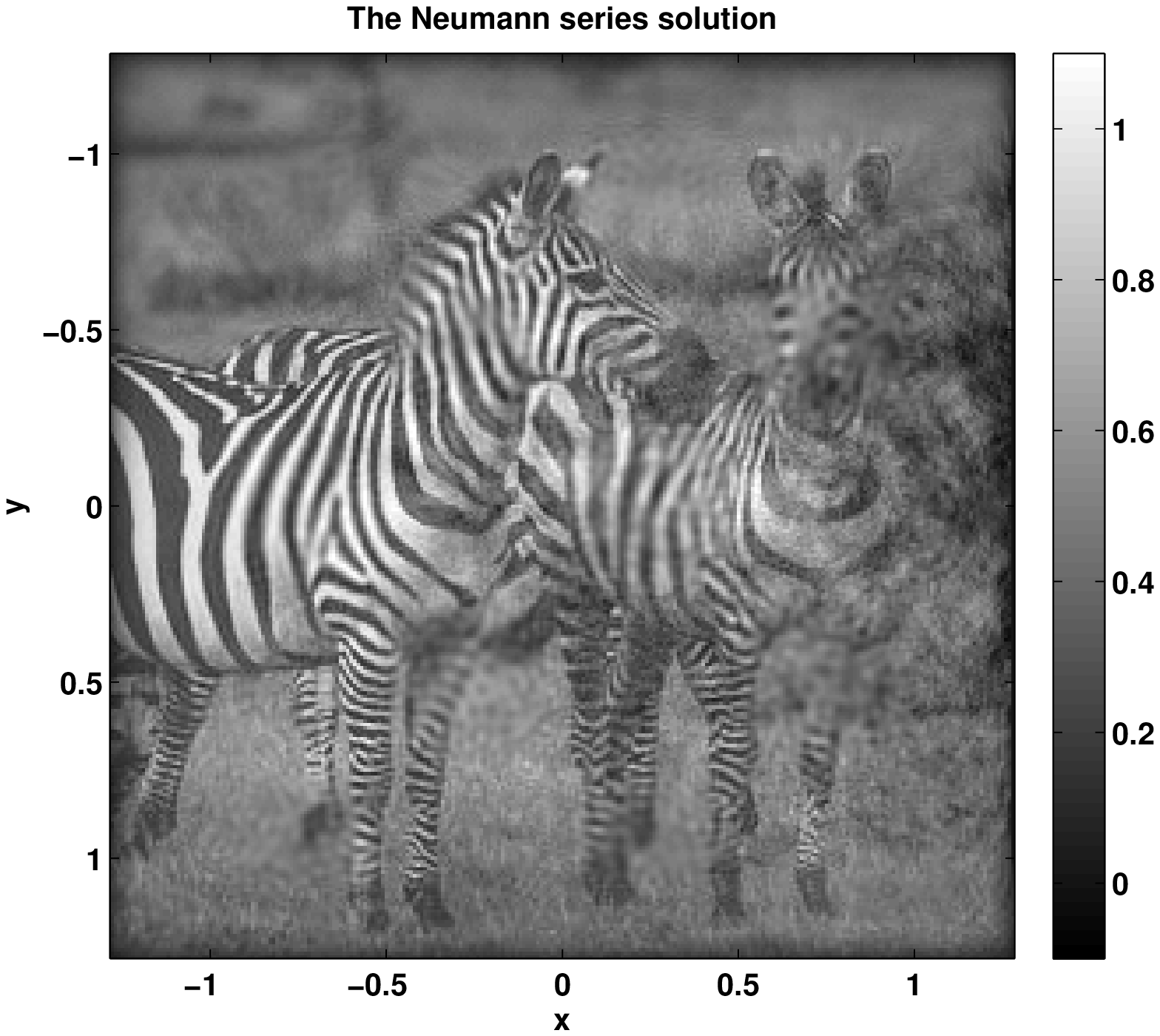
,height=5.5cm }\\ 
(c)\epsfig{figure=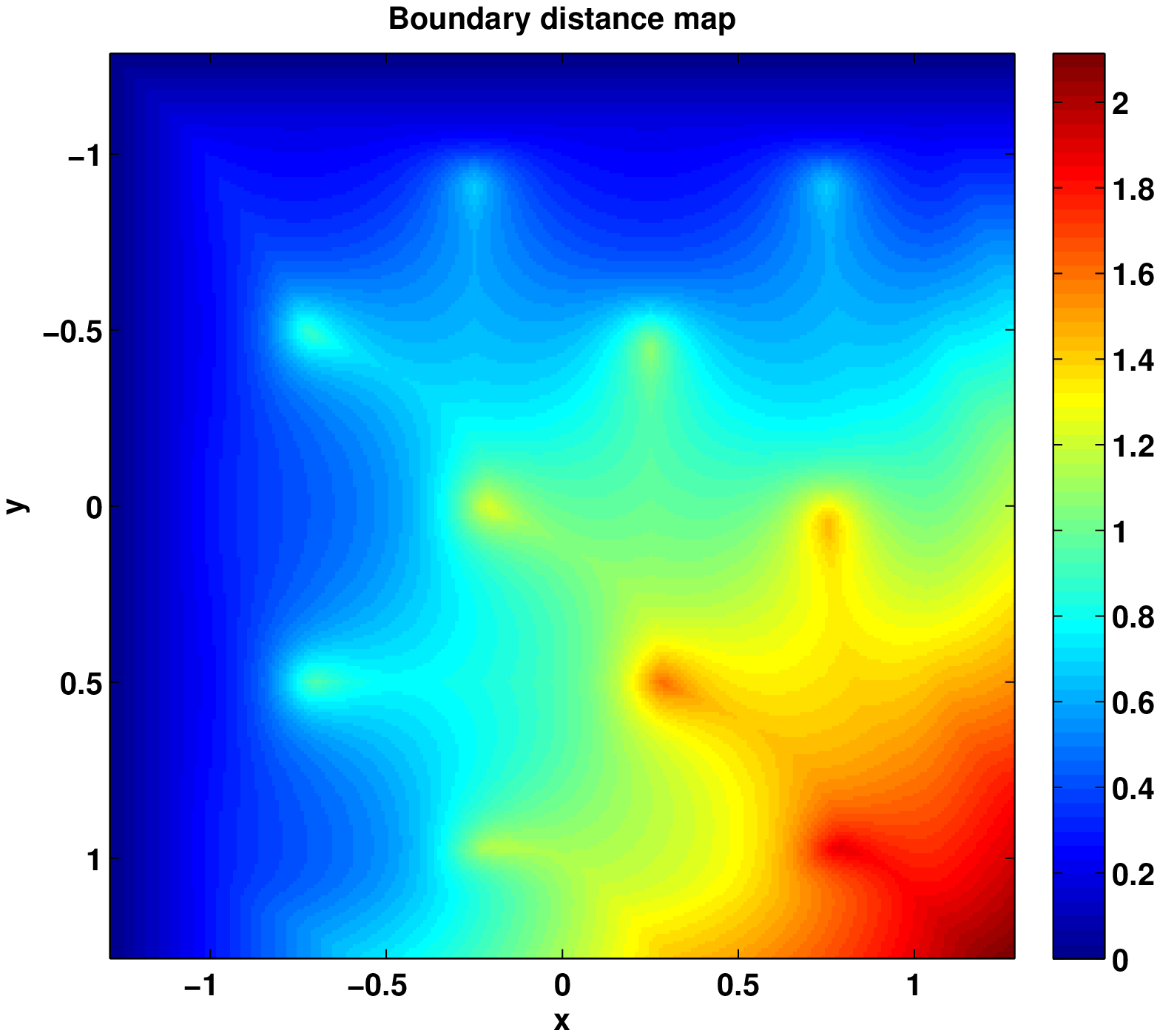,height=5.5cm }
(d)\epsfig{figure=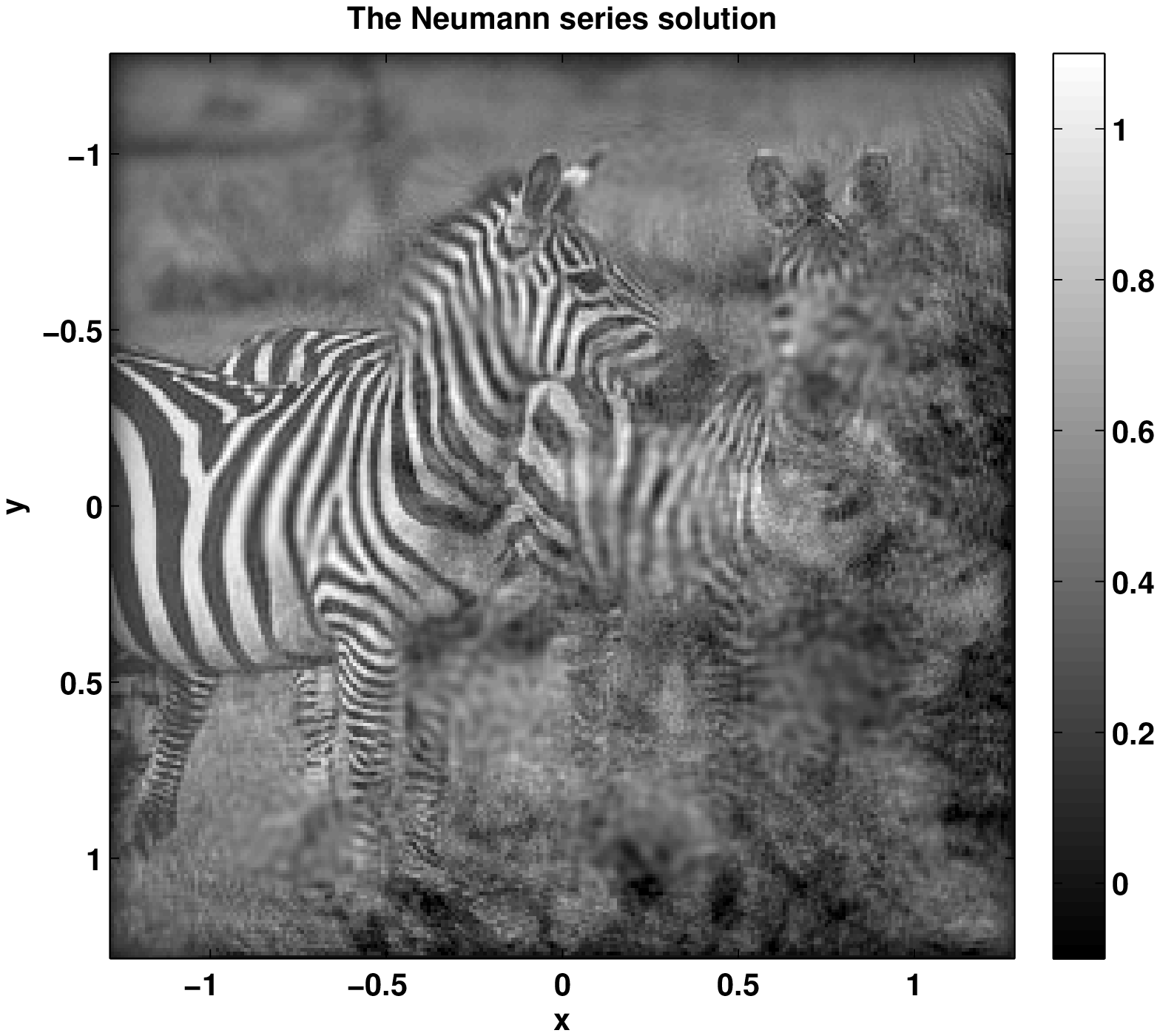,height=5.5cm }\\
(e)\epsfig{figure=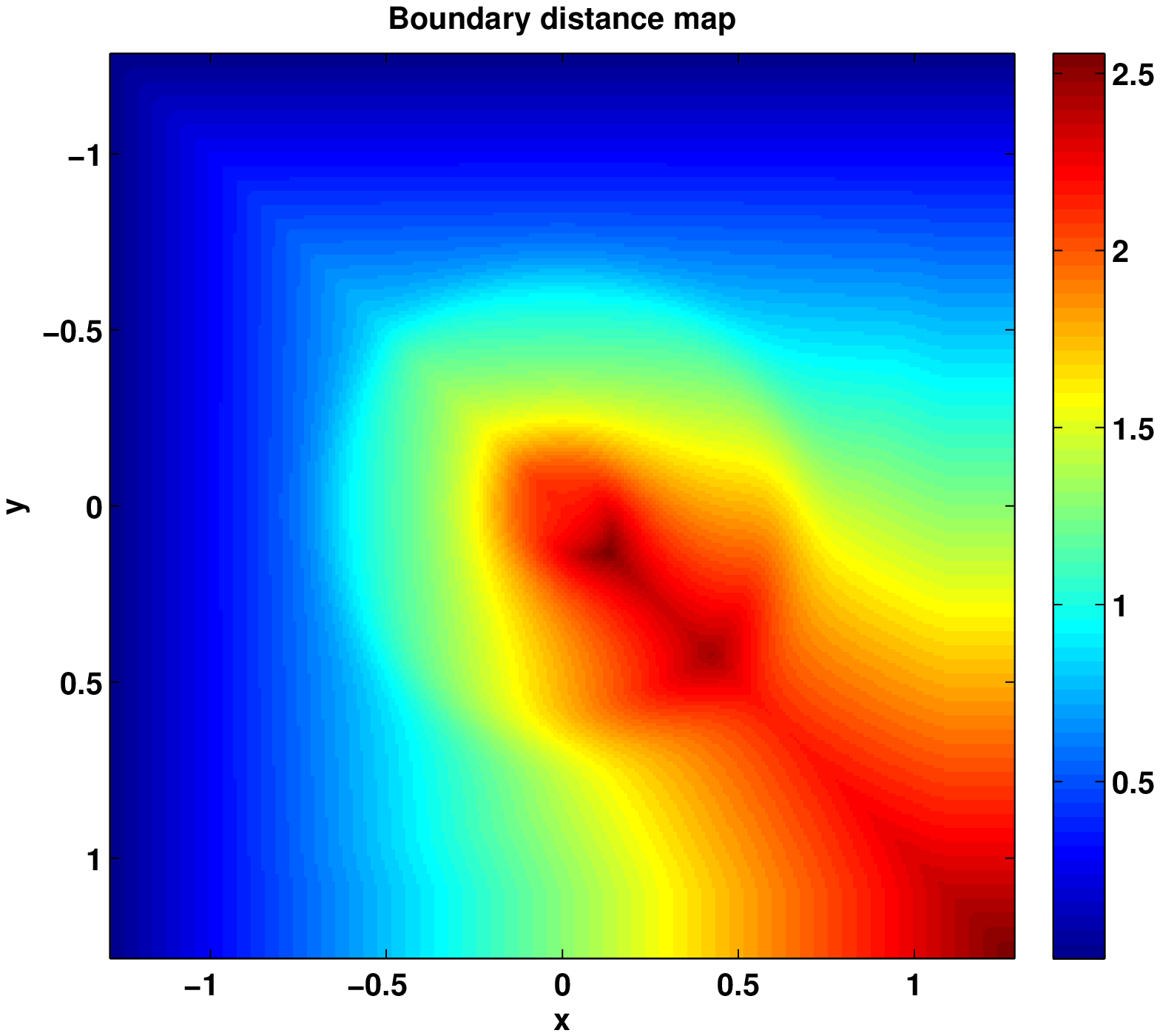,height=5.5cm }
(f)\epsfig{figure=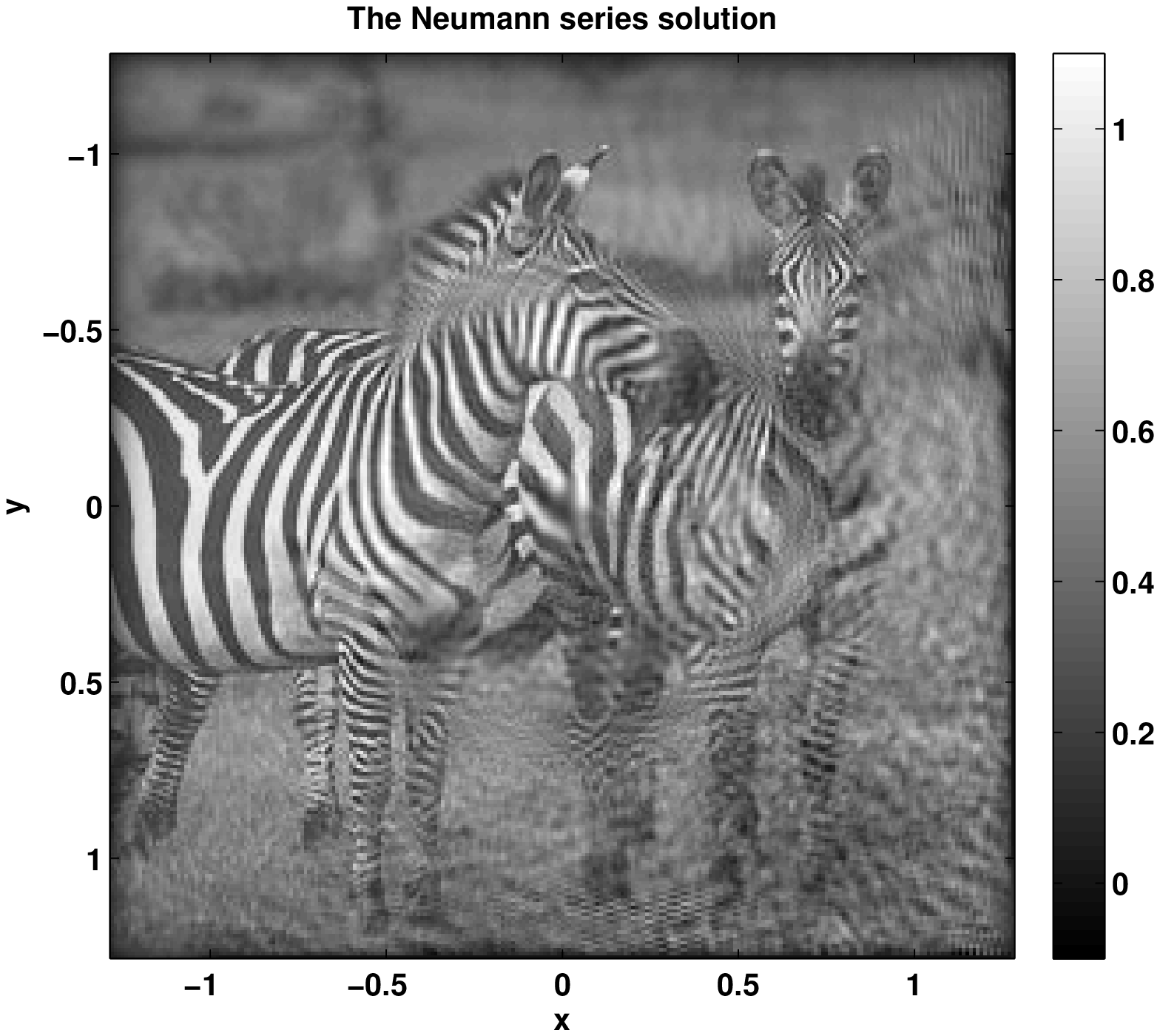,height=5.5cm }
\end{center}
\caption{Examples with the trapping speeds $c_3$ (the first two rows), and $c_2$  (the last row). 
(a): the boundary distance map, data on three sides. 
(b): the reconstructed ``zebras'' image, $T=4.92$.
(c): the boundary distance map, data on two sides. 
(d): the reconstructed ``zebras'' image, $T=4.92$, $k=16$.
(e): the boundary distance map, data on two sides. 
(f): the reconstructed ``zebras'' image, $T=8.62$, $k=16$. 
}
\label{Fig:2dZebra4T0trapPartialOne}
\end{figure}

\subsubsection{Figure~\ref{Fig:2dZebra4T0trapPartialOne} (e), (f)} 
Here,  $T=8.61>T_0\approx 2.6$  with data on two sides, and the speed is $c_2$. The time $T$ is larger than above but $T_0$ is slightly larger as well. The reconstruction is better due to the larger time and (probably) due to the fact that the trapping region of this speed is farther away from the sides where no observations are done. Experiments with times $T$ closer to that in the two examples above, not shown, still yield a better reconstruction with this speed. 

%
%

\section{Conclusion}
We present new algorithms for reconstructing an unknown source in TAT and PAT based on the recent advances in understanding 
the theoretical nature of the problem. We work with variable sound speeds that might be also discontinuous across some surface. 
The latter problem arises in brain imaging. The new algorithm is based on an explicit formula in the form of a Neumann series. 
We present numerical examples with non-trapping, trapping and piecewise smooth speeds, as well as examples with data on a part 
of the boundary. These numerical examples demonstrate the robust performance of the new algorithm. 


\section*{Acknowledgement}
Qian is partially supported by NSF 0810104 and NSF 0830161. Stefanov partially supported by NSF  Grant DMS-0800428. Uhlmann partially supported by NSF, a Chancellor Professorship at UC Berkeley and a Senior 
Clay Award. Zhao is partially supported by NSF 0811254.


\bibliographystyle{abbrv}
\bibliography{myref,myreferences,TAT}
\end{document}